  %
  %

  \font \eightbf = cmbx8
  \font \eightit = cmti8
  \font \eightrm = cmr8
  \font \eightsl = cmsl8
  \font \eighttt = cmtt8 \hyphenchar\eighttt = -1
  
  \font \tensc = cmcsc10
  
  \font \titlefont = cmr7 scaled \magstep4
  \font \eighti = cmmi8 \skewchar \eighti = '177

  \font \eightsy = cmsy8 \skewchar \eightsy = '60
  
  \font \ninebf = cmbx9
  \font \ninei = cmmi9 \skewchar \ninei = '177
  \font \nineit = cmti9
  \font \ninerm = cmr9
  \font \ninesl = cmsl9
  \font \ninesy = cmsy9 \skewchar \ninesy = '60
  \font \ninett = cmtt9 \hyphenchar \ninett = -1
  \font \sixbf = cmbx6
  \font \sixi = cmmi6 \skewchar \sixi = '177
  \font \sixrm = cmr6
  \font \sixsy = cmsy6 \skewchar \sixsy = '60
  
   \hyphenchar \tentt = -1

  \def \Headlines #1#2{\nopagenumbers 
    \voffset = 2\baselineskip 
    \advance \vsize by -\voffset 
    \headline {\ifnum \pageno = 1 \hfil 
    \else \ifodd \pageno \tensc \hfil \lowercase {#1} \hfil \folio 
    \else \tensc \folio \hfil \lowercase {#2} \hfil 
    \fi \fi }}

  \def \Title #1{\noindent {\titlefont #1}}
  
  \def \Date #1 {\footnote {}{\eightit Date: #1.}}

  \def \Note #1{\footnote {}{\eightrm #1}}

  \def \Authors #1{\bigskip \bigskip \noindent #1}

  \long \def \Addresses #1{\begingroup \eightpoint \parindent0pt 
\medskip #1\par \par \endgroup }

  \long \def \Abstract #1{\begingroup \eightpoint 
  \bigskip \bigskip \noindent 
  {\sc ABSTRACT.} #1\par \par \endgroup }

  \def \Contents {\bigskip \centerline {CONTENTS} \medskip }

  \newcount \ccnt \ccnt = 0
  \newdimen \ldim 
  \newdimen \rdim 
  \newdimen \ndim 
  \newdimen \mdim 
  \ldim = 0.1\hsize 
  \rdim = 0.1\hsize 
  \ndim = 15pt
  \mdim = \hsize 
  \advance \mdim by -\ldim 
  \advance \mdim by -\rdim 
  \advance \mdim by -\ndim 
  \def \triplebox #1#2#3{\noindent \hbox {\hbox to \ldim {\hfil
#1}\hbox to \mdim { #2\dotfill }\hbox to \rdim {\hbox to \ndim {\hfil 
#3}\hfil }}\par }

  \def \pag #1 ... #2{\triplebox {\global \advance \ccnt by 1 \number 
\ccnt .}{#2}{#1}}
  \def \ppg #1 ... #2{\triplebox {}{#2}{#1}}

  \def \vg #1{\ifx #1\null \null \else 
    \ifx #1\ { }\else 
    \ifx #1,,\else 
    \ifx #1..\else 
    \ifx #1;;\else 
    \ifx #1::\else 
    \ifx #1''\else 
    \ifx #1--\else 
    \ifx #1))\else 
    { }#1\fi \fi \fi \fi \fi \fi \fi \fi \fi }

  \def \goodbreak {\vskip0pt plus.1\vsize \penalty -250 \vskip0pt 
plus-.1\vsize }

  \newcount \secno \secno = 0
  \newcount \stno 

  \def \seqnumbering {\global \advance \stno by 1
    \number \secno .\number \stno }

  \def \label #1{\def\localsystemvariable {\number \secno 
    \ifnum \number \stno = 0\else .\number \stno \fi }
    \global \edef #1{\localsystemvariable }}

  \def \section #1{\stno = 0 \global \advance \secno by 1
    \bigskip \bigskip \goodbreak
    \noindent {\bf \number \secno .\enspace #1.}
    \bigskip}

  \def \References {
    \begingroup
    \bigskip \bigskip \goodbreak
    \eightpoint
    \centerline {\tensc References}
    \nobreak \medskip \frenchspacing }

  \long \def \sysstate #1#2#3{\medbreak \noindent {\bf 
\seqnumbering .\enspace #1.\enspace }{#2#3\vskip 0pt}\medbreak }
  \def \state #1 #2\par {\sysstate {#1}{\sl }{#2}}
  \def \definition #1\par {\sysstate {Definition}{\rm }{#1}}

  \def \proof {\medbreak \noindent {\it Proof.\enspace }}
  \def \proofend {\ifmmode \eqno \square \else \hfill \square 
\looseness = -1 \medbreak \fi }

  \newcount \zitemno \zitemno = 0
  \def \stdcount {\global \advance \zitemno by 1 \ifcase \zitemno \or 
i\or ii\or iii\or iv\or v\or vi\or vii\or viii\or ix\or x\or xi\or 
xii\fi }
  \def \iItem #1{\medskip \item {#1}}
  \def \Item #1{\smallskip \item {#1}}

  \def \izitem {\zitemno = 0\iItem {{\rm (\stdcount )}}}
  \def \zitem {\Item {{\rm (\stdcount )}}}

  \def \iBitem {\zitemno = 0\iItem {$\bullet $}}
  \def \Bitem {\Item {$\bullet $}}

  \newdimen \itemwidth \itemwidth = 40pt \newdimen \itemcontents 
  \def \witem #1#2{\medskip 
  \itemcontents = \hsize \advance \itemcontents by -\itemwidth 
  \hbox {\hbox to \itemwidth {\hfill #1\ }\vtop {\noindent \parshape 1
0pt \itemcontents \ #2}}\par }

  \def \$#1{#1 $$$$ #1}

  \def \({\left (}
  \def \){\right )}
  \def \[{\left \Vert }
  \def \]{\right \Vert }
  \def \*{\otimes }
  \def \+{\oplus }
  \def \:{\colon }
  \def \<{\langle }
  \def \>{\rangle }
  \def \text #1{\hbox {\rm #1}}
  \def \and {\hbox {\quad and \quad }}
  \def \arw {\rightarrow }
  \def \calcat #1{\,{\vrule height8pt depth4pt}_{\,#1}}
  
  \def \crossproduct {\hbox to 1.8ex{$\times \kern -.45ex\vrule 
height1.1ex depth0pt width0.45truept$\hfill }}
  \def \cstar {$C^*$}
  \def \for #1{,\quad #1}
  \def \inv {^{-1}}
  
  \def \square {\hbox {$\sqcap \!\!\!\!\sqcup $}}
  \def \stress #1{{\it #1}\/}
  
  \def \x {\times }
  \def \|{\Vert }
  \def \inv {^{-1}}

  \newskip \ttglue 
  \def \tenpoint {\def \rm {\fam0 \tenrm }%
  \textfont0 = \tenrm 
  \scriptfont0 = \sevenrm \scriptscriptfont0 = \fiverm 
  \textfont1 = \teni 
  \scriptfont1 = \seveni \scriptscriptfont1 = \fivei 
  \textfont2 = \tensy 
  \scriptfont2 = \sevensy \scriptscriptfont2 = \fivesy 
  \textfont3 = \tenex 
  \scriptfont3 = \tenex \scriptscriptfont3 = \tenex 
  \def \it {\fam \itfam \tenit }%
  \textfont \itfam = \tenit 
  \def \sl {\fam \slfam \tensl }%
  \textfont \slfam = \tensl 
  \def \bf {\fam \bffam \tenbf }%
  \textfont \bffam = \tenbf 
  \scriptfont \bffam = \sevenbf 
  \scriptscriptfont \bffam = \fivebf 
  \def \tt {\fam \ttfam \tentt }%
  \textfont \ttfam = \tentt 
  \tt \ttglue = .5em plus.25em minus.15em
  \normalbaselineskip = 12pt
  \def \MF {{\manual META}\-{\manual FONT}}%
  \let \sc = \eightrm 
  \let \big = \tenbig 
  \setbox \strutbox = \hbox {\vrule height8.5pt depth3.5pt width 0pt}%
  \normalbaselines \rm }

  \def \ninepoint {\def \rm {\fam0 \ninerm }%
  \textfont0 = \ninerm 
  \scriptfont0 = \sixrm \scriptscriptfont0 = \fiverm 
  \textfont1 = \ninei 
  \scriptfont1 = \sixi \scriptscriptfont1 = \fivei 
  \textfont2 = \ninesy 
  \scriptfont2 = \sixsy \scriptscriptfont2 = \fivesy 
  \textfont3 = \tenex 
  \scriptfont3 = \tenex \scriptscriptfont3 = \tenex 
  \def \it {\fam \itfam \nineit }%
  \textfont \itfam = \nineit 
  \def \sl {\fam \slfam \ninesl }%
  \textfont \slfam = \ninesl 
  \def \bf {\fam \bffam \ninebf }%
  \textfont \bffam = \ninebf 
  \scriptfont \bffam = \sixbf 
  \scriptscriptfont \bffam = \fivebf 
  \def \tt {\fam \ttfam \ninett }%
  \textfont \ttfam = \ninett 
  \tt \ttglue = .5em plus.25em minus.15em
  \normalbaselineskip = 11pt
  \def \MF {{\manual hijk}\-{\manual lmnj}}%
  \let \sc = \sevenrm 
  \let \big = \ninebig 
  \setbox \strutbox = \hbox {\vrule height8pt depth3pt width0pt}%
  \normalbaselines \rm }

  \def \eightpoint {\def \rm {\fam0 \eightrm }%
  \textfont0 = \eightrm 
  \scriptfont0 = \sixrm \scriptscriptfont0 = \fiverm 
  \textfont1 = \eighti 
  \scriptfont1 = \sixi \scriptscriptfont1 = \fivei 
  \textfont2 = \eightsy 
  \scriptfont2 = \sixsy \scriptscriptfont2 = \fivesy 
  \textfont3 = \tenex 
  \scriptfont3 = \tenex \scriptscriptfont3 = \tenex 
  \def \it {\fam \itfam \eightit }%
  \textfont \itfam = \eightit 
  \def \sl {\fam \slfam \eightsl }%
  \textfont \slfam = \eightsl 
  \def \bf {\fam \bffam \eightbf }%
  \textfont \bffam = \eightbf 
  \scriptfont \bffam = \sixbf 
  \scriptscriptfont \bffam = \fivebf 
  \def \tt {\fam \ttfam \eighttt }%
  \textfont \ttfam = \eighttt 
  \tt \ttglue = .5em plus.25em minus.15em
  \normalbaselineskip = 9pt
  \def \MF {{\manual opqr}\-{\manual stuq}}%
  \let \sc = \sixrm 
  \let \big = \eightbig 
  \setbox \strutbox = \hbox {\vrule height7pt depth2pt width0pt}%
  \normalbaselines \rm }


  \newcount \bibno \bibno =0
  \def \newbib #1{\global \advance \bibno by 1 \edef #1{\number \bibno }}
  \def \cite #1{{\rm [\bf #1\rm ]}}
  \def \scite #1#2{[{\bf #1}:#2]}
  \def \lcite #1{#1}

  \def \se #1#2#3#4{\def \a {#1}\def \b {#2}\ifx \a \b #3\else #4\fi }
  \def \setem #1{\se {#1}{}{}{, #1}}
  \def \index #1{\smallskip \item {[#1]}}
  \def \tit #1{``#1''}

  \def \Article #1{\zarticle #1 xyzzy }
  \def \zarticle #1, auth = #2,
   title = #3,
   journal = #4,
   year = #5,
   volume = #6,
   pages = #7,
   NULL#8 xyzzy {\index {#1} #2, \tit {#3}, {\sl #4\/} {\bf #6} (#5), #7.}

  \def \ztechreport #1,
    auth = #2,
    title = #3,
    institution = #4,
    year = #5,
    type = #6,
    note = #7,
    NULL#8
    xyzzy
    {\index {#1} #2, \tit {#3}\setem {#6}\setem {#4}\setem {#5}\setem {#7}.}

  \def \Unpublished #1{\zunpublished #1 xyzzy }
  \def \zunpublished #1,
    auth = #2,
    title = #3,
    institution = #4,
    year = #5,
    type = #6,
    note = #7,
    NULL#8
    xyzzy
    {\index {#1} #2, \tit {#3}\setem {#6}\setem {#4}\setem {#5}\setem {#7}.}

  \def \Book #1{\zbook #1 xyzzy }
  \def \zbook #1,
    auth = #2,
    title = #3,
    publisher = #4,
    year = #5,
    volume = #6,
    series = #7,
    NULL#8
    xyzzy
    {\index {#1} #2, \tit {#3}\setem {#7}\se {#6}{}{}{ vol. #6}, #4, #5.}

  \def \zmasterthesis #1,
    auth = #2,
    title = #3,
    school = #4,
    year = #5,
    type = #6,
    NULL#7
    xyzzy
    {\index {#1} #2, \tit {#3}, #6, #4, #5.}

  \def \zinproceedings #1,
    auth = #2,
    title = #3,
    booktitle = #4,
    year = #5,
    pages = #6,
    organization = #7,
    note = #8,
    NULL#9
    xyzzy
    {\index {#1} #2, \tit {#3}\se {#4}{}{}{, In
    {\sl #4}}\setem {#5}\setem {#6}\setem {#7}\setem {#8}.}

  \def \zphdthesis #1,
    auth = #2,
    title = #3,
    school = #4,
    year = #5,
    type = #6,
    NULL#7
    xyzzy
    {\index {#1} #2, \tit {#3}, #6, #4, #5.}

  \def \zbooklet #1,
    auth = #2,
    title = #3,
    howpublished = #4,
    year = #5,
    NULL#6
    xyzzy
    {\index {#1} #2, \tit {#3}, #4, #5.}

  \def \zmisc #1,
    auth = #2,
    title = #3,
    note = #4,
    howpublished = #5,
    year = #6,
    NULL#7
    xyzzy
    {\index {#1} #2, \tit {#3}\setem {#4}\setem {#5}\setem {#6}.}



  \def\fcite#1{\lcite{#1}}
  \def\Graph{{\rm Gr}(A)}
  \def\col#1{c_{#1}}
  \def\STree#1#2{F(#1/#2)}
  \def\trunc#1#2{#1|_{#2}}
  \def\root#1#2{R_{#2}(#1)}
  \def\orb#1{{\rm Orb}(#1)}
  \def\stem#1{\sigma(#1)}
  \def\Gen{{\cal G}}
  \def\F{{I\!\!F}}
  \def\Pos{\F_+}
  \def\RelTCK{\Rel_A^\tau}
  \def\RelCK{\Rel_A}
  \def\LA{{\cal T\kern-2ptO}_A}
  \def\OA{{\cal O}_A}
  \def\uOA{\widetilde{\cal O}_A}
  \def\OmegaE{\Omega^e}
  \def\GenSpec#1{\Omega(#1)}
  \def\specTCK{\Omega_A^\tau}
  \def\specUCK{\widetilde{\Omega}_A}
  \def\specCK{\Omega_A}
  \def\ha{\h^A}
  \def\Xa{\Delta^A}
  \def\PR#1{{P\kern-1ptR}$_{#1}$}
  \def\TCKCond#1{T\kern-1ptC\kern-1ptK$_{#1}$}
  \def\CKcond#1{CK$_{#1}$}
  \def\UAlg#1#2{C^*(#1,#2)}
  \def\W#1#2#3{W_{#1,#2,#3}}
  \def\A#1#2{A(#1,#2)}
  \def\N{{\bf N}}
  \def\B{{\cal B}}
  \def\E{{\cal E}}
  \def\e{{\rm e}}
  \def\a{\alpha}
  \def\b{\beta}
  \def\d{\delta}
  \def\g{\gamma}
  
  \def\w{\omega}
  \def\last#1{#1_{|#1|}}
  \def\piso{partial isometry\vg}
  \def\pisos{partial isometries\vg}
  \def\pr{u}
  \def\Rel{{\cal R}}
  \def\h{h}
  \def\X{\Delta}
  \def\prep{partial representation\vg}
  \def\preps{partial representations\vg}
  \def\ss{semi-saturated\vg}
  \def\claim#1 #2\par{\medskip\noindent{\tensc claim} #1: #2\medskip}
  \def\commute#1#2{#1 and #2 commute}
  \def\imply{\Rightarrow}
  \def\bimply{\Leftrightarrow}
  \def\IN{\!\in\!}
  \def\NOTIN{\!\notin\!}
  \def\INVIN{\inv\!\IN}
  \def\INVNOTIN{\inv\!\NOTIN}
  \def\set#1{[\![#1]\!]}
  \def\0{\emptyset}
  \def\labar#1{\buildrel\hbox{$#1$} \over \longrightarrow}
  \def\frac#1#2{{#1 \over #2}}
  \def\tfrac#1#2{\frac{#1}{#2}}
  \def\bvec{\varepsilon}


  \newbib\FAbadie
  \newbib\ArchSpiel
  \newbib\Blackadar
  \newbib\ChoiEffros
  \newbib\CuntzOInfinite
  \newbib\Cuntz
  \newbib\CKbib
  \newbib\Evans
  \newbib\newpim
  \newbib\Inverse
  \newbib\TPA
  \newbib\Amena
  \newbib\Ortho
  \newbib\ELQ
  \newbib\FD
  \newbib\Glf
  \newbib\HR
  \newbib\KPRR
  \newbib\KPR
  \newbib\Laca
  \newbib\LR
  \newbib\LS
  \newbib\McCl
  \newbib\Nica
  \newbib\QR
  \newbib\Renault

  \Headlines
  {Cuntz--Krieger Algebras for Infinite Matrices} 
  {Ruy Exel and Marcelo Laca} 

  \Date {January 28, 1997}

  \Title
  {Cuntz--Krieger Algebras for Infinite Matrices}\footnote
  {$\null^\dagger$}
  {\eightrm Research partially supported by FAPESP.}

  \Authors
  {Ruy Exel\footnote
    {*}{\eightrm Partially supported by CNPq.}
  and Marcelo Laca\footnote
    {**}{\eightrm Supported by the Australian Research Council.}}

  \Note {Note: This version incorporates a correction in Theorem 9.1.}

  \Addresses
  {Departamento de Matem\'atica;
  Universidade de S\~ao Paulo;
  Rua do Mat\~ao, 1010;
  05508-900 S\~ao Paulo -- Brazil, and
  \par
  Department of Mathematics;
  University of Newcastle;
  NSW 2308 -- Australia.}

  \Abstract
  {Given an arbitrary infinite matrix $A = \{\A ij\}_{i,j\in \Gen}$ with
entries in $\{0,1\}$ and having no identically zero rows, we define an
algebra $\OA$ as the universal \cstar-algebra generated by partial
isometries subject to conditions that generalize, to the infinite case,
those introduced by Cuntz and Krieger for finite matrices.  We realize
$\OA$ as the crossed product algebra for a partial dynamical system and,
based on this description, we extend to the infinite case some of the
main results known to hold in the finite case, namely the uniqueness
theorem, the classification of ideals, and the simplicity criteria.
$\OA$ is always nuclear and we obtain conditions for it to be unital and
purely infinite.}

  \Contents

  \pag 2 ... {Introduction}
  \pag 4 ... {Preliminaries}
  \pag 7 ... {Toeplitz--Cuntz--Krieger partial isometries}
  \pag 10 ... {The Toeplitz--Cuntz--Krieger algebra}
  \pag 13 ... {Stems and roots}
  \pag 16 ... {The topology of $\specTCK$}
  \pag 18 ... {The unital Cuntz--Krieger algebra}
  \pag 23 ... {Cuntz--Krieger algebras for infinite matrices}
  \pag 26 ... {A faithful representation}
  \pag 28 ... {Invariant sets}
  \pag 31 ... {Fixed points}
  \pag 32 ... {Topological freeness}
  \pag 34 ... {Uniqueness of $\OA$}
  \pag 35 ... {The simplicity criteria}
  \pag 35 ... {Classification of ideals}
  \pag 36 ... {Pure infiniteness}
  \ppg 39 ... {Index of symbols}
  \ppg 40 ... {References}

  \section{Introduction}
  \label\Introduction
  Let $\Gen$ be a \stress{finite} set and let $A = \{\A ij\}_{i,j\in
\Gen}$ be a matrix with entries in the set $\{0,1\}$, having no
identically zero rows.  The Cuntz--Krieger algebra $\OA$, introduced and
studied in \cite{\CKbib} in the special case where $A$ satisfies a
certain condition (I) (see the bottom of page 254 in \cite{\CKbib}), may
also be defined as the universal \cstar-algebra generated by a family of
partial isometries $\{S_i\}_{i\in\Gen}$ subject to the conditions
  \itemwidth = 40pt
  \witem{\CKcond1)}{$1=\sum_{j\in\Gen} S_j S_j^*$, and}
  \witem{\CKcond2)}{$S_i^* S_i = \sum_{j\in\Gen} \A ij S_j S_j^*$, for
all $i\in\Gen$.}

  \medskip\noindent

If the index set $\Gen$ is \stress{infinite}, any attempt at
  generalizing the theory of Cuntz and Krieger will face a serious
stumbling block at the first onset, the problem being that the series in
\lcite{\CKcond1} and \lcite{\CKcond2} will not converge in norm.

A natural attempt to deal with these infinite sums is to work with the
strong operator topology \scite{\CKbib}{Remark 2.15} since its terms may
be assumed to be pairwise orthogonal projections, in analogy with the
finite case.  However, this leads one away from the standard \cstar{}
concepts because the theory of universal \cstar-algebras for sets of
generators and relations \cite{\Blackadar} breaks down when the
relations involve the strong topology.

One way to side-step the problem above (and to remain in the category of
\cstar-algebras) is to assume that $A$ is row-finite, i.e., that there
are only finitely many ones in each row of $A$.  In this case one may
still consider families of partial isometries $\{S_i\}_{i\in \Gen}$
satisfying \lcite{\CKcond2} because the sum involved will only have a
finite number of nonzero terms.  Although condition \lcite{\CKcond1} no
longer makes sense in the infinite case, as it involves a divergent sum,
one is justified in abandoning it, since, in the finite case, it is
actually unnecessary: the right hand side of \lcite{\CKcond1} ends up
behaving like a unit \scite{\CKbib}{Section 2} even if \lcite{\CKcond1}
is not required, as long as one retains the requirement that the ranges
of the $S_i$ be pairwise orthogonal.  The row-finite case is the subject
of \cite{\KPRR} and \cite{\KPR}, where Renault's theory of groupoid
\cstar-algebras \cite{\Renault} was used to extend some of the main
results originally proved by Cuntz and Krieger \cite{\CKbib,\Cuntz}.

Beyond the row-finite case there is one and only one known example: it
is widely recognized that the algebra ${\cal O}_\infty$ introduced in
\cite{\CuntzOInfinite} is a model for $\OA$ when $A$ is the (countably)
infinite matrix whose entries are all equal to 1.  However, due to the
lack of a general theory, this recognition is merely based on the formal
analogy with the finite case.  In truth ${\cal O}_\infty$ is but a
beacon, signalling towards a hitherto elusive theory of Cuntz--Krieger
algebras for genuinely infinite matrices.

It is the goal of this work to give a definition for $\OA$ and to
develop its theory entirely within the context of \cstar-algebras for an
arbitrary infinite 0--1 matrix having no identically zero rows.

To motivate our definition suppose that the index set $\Gen$ is infinite
but the $i^{th}$ row of $A$ contains only a finite number of zeros.
Ignoring the convergence problems in \lcite{\CKcond1} and
\lcite{\CKcond2}, subtract \lcite{\CKcond2} from \lcite{\CKcond1} to
derive
  $$
  1 - S_i^* S_i = \sum_{j\in\Gen} (1-\A ij) S_j S_j^*,
  \eqno{(\seqnumbering)}
  \label \ComplementCondition
  $$
  where the sum is now finite, since only a finite number of
coefficients $(1-\A ij)$ do not vanish.

There is another meaningful conclusion to be inferred from
\lcite{\CKcond2}: suppose that $i$ and $i'$ are labels for two rows of
$A$ such that, for all but finitely many $j$'s, either $\A ij$ or $\A
{i'}j$ vanishes.  An alternative way to say this is that the set
$\{j\in\Gen : \A ij \A{i'}j = 1\}$ is finite.  If we then write down
\lcite{\CKcond2} for both $i$ and $i'$ and, again ignoring the
convergence problems, multiply them term-wise, we get
  $$
  S_i^* S_i ~ S_{i'}^* S_{i'} = \sum_{j\in\Gen} \A ij \A{i'}j S_j S_j^*,
  $$
  where we have also assumed that the $S_j S_j^*$ are pairwise
orthogonal projections.  Again the sum happens to be finite.

These two ways of extracting meaningful conclusions from the
  Cuntz--Krieger equations above may be combined as follows: suppose
  $X$ and $Y$ are finite subsets of $\Gen$ such that
  $$
  A(X,Y,j) :=
  \prod_{x\in X} \A xj
  \prod_{y\in Y} (1-\A yj)
  \eqno{(\seqnumbering)}
  \label \AXYj
  $$
  is nonzero only for a finite number of $j$'s.  If we write
\lcite{\CKcond2} for each $x$ in $X$ and \lcite{(\ComplementCondition)}
for each $y$ in $Y$ and then multiply term-wise the resulting equations
we will arrive at
  $$
  \prod_{x\in X} S_x^* S_x \prod_{y\in Y} (1 - S_y^* S_y) =
  \sum_{j\in\Gen} A(X,Y,j) S_j S_j^*,
  \eqno{(\seqnumbering)}
  \label \ELCond
  $$
  which, once more, involves a finite sum.

We will then define $\OA$ (up to a subtle question of units) as the
universal \cstar-algebra generated by a family of partial isometries
$\{S_i\}_{i\in\Gen}$ satisfying:
  \izitem Condition \lcite{(\ELCond)}, for each pair of finite subsets
$X$ and $Y$ of $\Gen$ such that $A(X,Y,j)$ vanishes for all but a finite
number of $j$'s,
  \zitem \commute{$S_i^* S_i$}{$S_j^*S_j$}, for all $i,j$,
  \zitem $S_i^* S_j = 0$, if $i\neq j$, and
  \zitem $(S_i^* S_i) S_j = \A ij S_j $, for all $i,j$.

  \medskip
  In the case of finite matrices \lcite{(\ELCond)} is easily seen to be
equivalent to \lcite{\CKcond 1} and \lcite{\CKcond 2} which, in turn,
imply (ii--iv).  However, we will see that there are infinite matrices
$A$ with respect to which one cannot find a single pair of sets $X$ and
$Y$ such that $A(X,Y,j)$ is finitely supported in $j$. Condition (i) is
then vacuously satisfied, suggesting that we still need to include
(ii--iv) explicitly.  Moreover, we will see that conditions (ii--iv)
by themselves lead to an interesting theory as well as to an important
\cstar-algebra which we will call the Toeplitz--Cuntz--Krieger algebra.

As in \cite{\KPRR}, we associate to any given $A$ a graph called
$\Graph$ by taking the index set $\Gen$ as the set of vertices and by
joining a pair of vertices $(i,j)$ with an oriented edge from $i$ to $j$
precisely when $\A ij=1$.  The requirement that no row of $A$ is
identically zero says that there exists at least one edge out of each
vertex of $\Graph$.  However, since we do not require $A$ to be
row-finite, $\Graph$ may have an infinite number of edges in and out of
each vertex.

We then show that certain properties of $\Graph$ imply important
\cstar-algebraic properties of $\OA$.  Among these is the famous
uniqueness theorem of Cuntz and Krieger \scite{\CKbib}{Theorem 2.13}
which we show to follow, in the general infinite case, from the absence
of \stress{terminal circuits} in $\Graph$.  By a terminal circuit we
mean a circuit from which there is no exit (see below for the formal
definitions).  The property of having no terminal circuits is a
generalization, to the infinite case, of condition (I) of \cite{\CKbib},
as proven in \scite{\KPR}{Lemma 3.3}, so that our result implies the
classical uniqueness theorem when specialized to the finite case.

The next problem treated is the question of simplicity for $\OA$.
Precisely, we show that $\OA$ is simple whenever $\Graph$ is transitive
and not a finite cycle, generalizing the original simplicity result of
Cuntz and Krieger \scite{\CKbib}{Theorem 2.14}.

In the row-finite case simplicity is related to a property of $\Graph$
referred to as \stress{co-finality} in \cite{\KPRR} and \cite{\KPR}
--- a directed graph is said to be co-final if from each vertex one
can intercept any predetermined infinite path.  However, we will give
an example to show that the effectiveness of co-finality in proving
the simplicity of $\OA$ is restricted to the row-finite situation.
This turns out to be one of the surprises encountered when
venturing out of the row-finite situation.

We then discuss the classification of ideals of $\OA$ under another
property of $\Graph$, namely the absence of \stress{transitory
circuits}.  We say that a circuit in $\Graph$ is transitory if, upon
exiting it, one cannot come back to it by following the oriented edges
of $\Graph$.  This property is referred to as condition (K) in
\cite{\KPRR} and it is the natural generalization of Cuntz's condition
(II) of \cite{\Cuntz} to infinite matrices, as observed in
\scite{\KPRR}{5.1}.

Assuming the inexistence of transitory circuits in $\Graph$, we show
that the ideals of $\OA$ correspond bijectively to certain open sets in
a topological space $\specCK$ that are invariant under a partial action
of the free group $\F$ generated by $\Gen$.  Here we encounter another
surprise presented by the non-row-finite case, as the structure of
ideals of $\OA$ is shown to be more complicated than suggested by the
finite case \scite{\Cuntz}{2.5} or the row-finite one
\scite{\KPRR}{6.6}.

The method employed to obtain most of our results is the theory of
\stress{partial dynamical systems}, introduced in \cite{\newpim} for the
case of the group of integers and subsequently extended to more general
situations in \cite{\McCl}, \cite{\TPA}, and \cite{\Inverse}.  By a
partial dynamical system we mean a topological space $\Omega$ where a
group $G$ acts by means of partial homeomorphisms, that is,
homeomorphisms between open subsets of $\Omega$.  The reader will find
the complete definitions and a rough summary of the theory of partial
dynamical systems in the following section.

Given any partial dynamical system there is an associated
\stress{covariance} \cstar-algebra, also referred to as its
\stress{crossed product}.  One of our main technical results is that we
are able to exhibit $\OA$ as resulting from such a construction for a
partial dynamical system in which the acting group is the free group on
the set $\Gen$ of indices of $A$.

Once that description is obtained, we use the tools of the theory of
partial dynamical systems to obtain our main results.  The most
important among these tools are \stress{amenability}
\cite{\Amena,\Ortho} and \stress{topological freeness} \cite{\ELQ}.  The
concept of amenability is related to the equality between the reduced
and full crossed products.  On the other hand, a partial dynamical
system is said to be topologically free if the set of fixed points for
the partial homeomorphism associated to each non-trivial group element
has empty interior.

When combined, these properties imply that any non-trivial ideal of the
crossed product \cstar-algebra has a non-trivial intersection with the
subalgebra of continuous functions on the given topological space
\scite{\ELQ}{2.6}.  This result is the
cornerstone for most of our main results, namely the uniqueness theorem,
the simplicity criteria, and the classification of ideals.

One of the most concrete ways to describe a \cstar-algebra is by
exhibiting a faithful representation of it on a Hilbert space.  In the
case of $\OA$ we provide such a representation on the Hilbert space
$\ell^2(P_A)\otimes \ell^2(\F)$, where $P_A$ denotes the set of all
infinite paths in $\Graph$.  This representation can be described
quite simply by specifying that each generating partial isometry $S_x$
is sent to the operator $L_x \otimes \lambda_x$, where $\lambda_x$
refers, as usual, to the left regular representation of $\F$ and
  $L_x$ is given on the canonical basis
  $\{\bvec_\w: \w \in P_A \}$ of $\ell^2(P_A)$ by
  $$
  L_x (\bvec_\w) =
  \left\{\matrix{
  \bvec_{x\w} & \text{ if $x\w$ is an admissible path,}\cr
  0 \hfill & \text{ otherwise.}\hfill}
  \right. 
  $$

In our final section we present a necessary condition for $\OA$ to be
purely infinite.  There we show that if a circuit can be reached from
every vertex of $\Graph$ and if terminal circuits do not exist then
$\OA$ is purely infinite.

A recent result of F.~Abadie \cite{\FAbadie} shows that the
cross-sectional \cstar-algebra of a Fell bundle is nuclear provided the
bundle satisfies the approximation property of \cite{\Amena} and the
unit fiber algebra is itself nuclear.  This result can be used in
combination with our crossed product characterization of $\OA$ to show
that $\OA$ is a nuclear \cstar-algebra for any infinite 0--1 matrix
having no identically zero rows.

  \section{Preliminaries}
  \label \Preliminaries
  In this section we will conduct a very brief survey on \stress{partial
actions} and \stress{partial group representations}
\cite{\newpim,\McCl,\Inverse,\TPA,\QR,\ELQ} and \stress{amenability of
Fell bundles} \cite{\Amena,\Ortho} which play crucial roles in our study
of Cuntz--Krieger algebras for infinite matrices.

A \stress{partial action} of a group $G$ on a set $\Omega$ is, by
definition \scite{\Inverse}{1.2}, a pair
  $\Theta = \(\{\X_t\}_{t\in G}, \{\h_t\}_{t\in G}\)$,
  where, for each $t$ in $G$, $\X_t$ is a subset of $\Omega$ and
  $\h_t : \X_{t\inv} \arw \X_t$ is a bijection satisfying, for all $t$
and $s$ in $G$,
  \iBitem $\X_e = \Omega$ and $\h_e$ is the identity map on $\Omega$
      (here, as always, $e$ denotes the identity element of $G$),
  \Bitem $\h_t(\X_{t\inv}\cap \X_s) = \X_t\cap \X_{ts}$, and
  \Bitem $\h_t(\h_s(x)) = \h_{ts}(x)
       \for{x\in \X_{s\inv} \cap \X_{s\inv t\inv}}$.

  \medskip\noindent
  If $\Omega$ happens to be a topological space we will require, in
addition, that
  \iBitem each $\X_t$ is an open subset of $\Omega$, and that
  \Bitem $\h_t$ is a homeomorphism from $\X_{t\inv}$ onto $\X_t$.

  \medskip\noindent
  If $\Theta = \(\{D_t\}_{t\in G}, \{\theta_t\}_{t\in G}\)$ is a partial
action of $G$ on a \cstar-algebra $A$ we shall require, instead, that
  \iBitem each $D_t$ is a closed two-sided ideal of $A$, and that
  \Bitem $\theta_t$ is a *-isomorphism from
  $D_{t\inv}$ onto $D_t$.

  \medskip\noindent
  Given a partial action
  $\(\{\X_t\}_{t\in G}, \{\h_t\}_{t\in G}\)$
  of $G$ on a locally compact topological space $\Omega$, let
$D_t=C_0(\X_t)$ be identified, in the usual way, with the ideal of
functions in $C_0(\Omega)$ vanishing off $\X_t$.  Then the map
  $$
  \theta_t : f\in D_{t\inv} \mapsto f\cdot\h_t\inv \in D_t,
  $$ is easily seen to be a *-isomorphism from
  $D_{t\inv}$ onto $D_t$.  Moreover, we have that the pair
  $\(\{D_t\}_{t\in G}, \{\theta_t\}_{t\in G}\)$ is a partial action of
$G$ on $C_0(\Omega)$, which takes into account the requirements of the
definition of partial actions on \cstar-algebras.

Given any partial action of a group $G$ on a \cstar-algebra $A$, one may
construct a Fell bundle
  $\B = \{B_t\}_{t\in G}$ (see \cite{\FD} for an extensive study of Fell
bundles, also referred to as \cstar-algebraic bundles), called the
\stress{semi-direct product bundle} \scite{\TPA}{2.8} of $A$ by $G$, as
follows: each $B_t$ is taken to be $D_t$ and an element $a\in D_t$ is
denoted by $a\d_t$, when viewed as an element of $B_t$.  The
multiplication and involution operations
  $$
  B_t \times B_s {\buildrel\cdot \over \arw} B_{ts} \and
  B_t {\buildrel* \over \arw} B_{t\inv},
  $$
  are defined by
  $$
  (a_t\d_t)\cdot(b_s\d_s) = \theta_t\(\theta_t\inv(a_t)b_s\)\d_{ts},
  \and
  (a_t\d_t)^* = \theta_t\inv(a_t^*)\d_{t\inv},
  $$
  for $a_t$ in $D_t$ and $b_s$ in $D_s$.

The \stress{crossed product} of $A$ by $G$, denoted by
$A\crossproduct_\Theta G$, or simply $A\crossproduct G$ if $\Theta$ is
understood, is then defined to be the cross-sectional \cstar-algebra
$C^*(\B)$ of $\B$ \scite{\FD}{VIII.17.2}.
  On the other hand the \stress{reduced crossed product}, denoted by
  $A\crossproduct_{\Theta,r} G$, or simply
  $A\crossproduct_r G$, is the reduced cross-sectional algebra
$C_r^*(\B)$ of $\B$ \scite{\Amena}{2.3}.

There is a canonical epimorphism
  $$
  \Lambda: C^*(\B) \arw C_r^*(\B)
  $$
  which may or may not be an isomorphism.  When $\Lambda$ happens to be
isomorphic we say that $\B$ is an \stress{amenable Fell bundle}
  \scite{\Amena}{4.1},
  or that $\Theta$ is an \stress{amenable partial action}.  When $G$ is
an amenable group \cite{\Glf}, any Fell bundle over $G$ is amenable
\scite{\Amena}{4.7}.

A Fell bundle $\B = \{B_t\}_{t\in G}$ is said to satisfy the
\stress{approximation property} \scite{\Amena}{4.4} if there exists a
net
  $\( a_i \)_{i\in I}$
  of functions
  $a_i\: G \arw B_e$,
  which is uniformly bounded in the sense that
  $$
  \sup_{i\in I} \left\| \sum_{t\in G} a_i(t)^* a_i(t) \right\| < \infty,
  $$
  and such that, for all $b_t$ in each $B_t$, one has
  $$
  \lim_{i \rightarrow \infty} \sum_{r\in G} a_i(tr)^* b_t a_i(r)=
  b_t.
  $$

The relevance of the approximation property lies in the fact that all
Fell bundles that satisfy this property are necessarily amenable
\scite{\Amena}{4.6}.

Even if $G$ is non-amenable as a group there may be interesting Fell
bundles over $G$ which are amenable.  An important class of such
examples is given by \cite{\Ortho}, where it is shown that any Fell
bundle $\B$ over a (finitely or infinitely generated) free group $\F$ is
amenable provided that
  $\B$ is \stress{semi-saturated} in the sense that
  \Bitem $B_{ts} = B_tB_s$ (closed linear span), whenever
$|ts|=|t|+|s|$,
  \smallskip\noindent and \stress{orthogonal}, that is,
  \Bitem $B_x^*B_y=\{0\}$ if $x$ and $y$ are distinct generators of
$\F$.

  \medskip\noindent
  Here and elsewhere we denote by $|t|$ the \stress{length} of $t$, that
is, the number of generators appearing in the reduced decomposition of
$t$.

A fundamental source of examples of Fell bundles is given by partial
group representations.  Recall from \cite{\Inverse} that a
\stress{partial representation} of a group $G$ on a Hilbert space $H$ is
a map
  $\pr : G \rightarrow \B(H)$
  such that
  \itemwidth = 40pt
  \witem{\PR1)}{$\pr(t)\pr(s)\pr(s\inv) = \pr(ts)\pr(s\inv)$,}
  \witem{\PR2)}{$\pr(t\inv) = \pr(t)^*$, and}
  \witem{\PR3)}{$\pr(e) = I$,}

  \medskip \noindent
  for all $t,s\in G$, where $I$ is the identity operator on $H$.

Given a \prep{} $\pr$, one has that each $\pr(t)$ is a partial isometry
and hence that
  $$
  \e(t) := \pr(t)\pr(t)^*
  $$
  is a projection.  These projections commute with each other
\cite{\Inverse} and hence generate a commutative \cstar-algebra of
operators on $H$, which we will denote by $B_e$.  For each $t$ in $G$,
let us also consider the set $B_t=B_e\pr(t)$, which is easily seen to be
a closed linear subspace of $\B(H)$.  It is not hard to show that, with
the induced operations from $\B(H)$, the collection $\B=\{B_t\}_{t\in
G}$ forms a Fell bundle \scite{\Amena}{6.1}.

Interestingly enough, this Fell bundle also arises as the semi-direct
product bundle for a partial action of $G$ on $B_e$: the ideals $D_t$
are defined by
  $D_t=B_t B_{t\inv}$ and the partial automorphisms are given by
  $$
  \theta_t(a) = \pr(t)a\pr(t\inv) \for a\in D_{t\inv}.
  $$
  The proof of this fact is elementary and is based on the isometry
between $D_t$ and $B_t$ given by
  $x\mapsto x\pr(t)$.

Let $\F$ be the free group on a set $\Gen$ of arbitrary cardinality.
Then a \prep{} $\pr$ of
  $\F$ is said to be \stress{semi-saturated} if
  \Bitem $\pr(ts) = \pr(t) \pr(s)$, whenever
  $|ts| = |t|+|s|$.

  \smallskip\noindent Also, $\pr$ is said to be \stress{orthogonal} if
  \Bitem $\pr(x)^*\pr(y)=0$, when $x$ and $y$ are in $\Gen$ and $x\neq
y$.

  \medskip
  It is not hard to see that if $\pr$ possesses the last two properties
mentioned then the associated bundle has the corresponding
properties of semi-saturatedness and orthogonality and hence is
amenable, as mentioned above.  In fact, for the special case of Fell
bundles arising from semi-saturated orthogonal \preps of $\F$, one has
that the approximation property holds \scite{\Ortho}{4.1}.

One is often interested in studying \preps of a group $G$ such that the
corresponding $\e(t)$ satisfy previously specified relations.  Let
$\Rel$ be a given set of such relations.  One may then consider the
\stress{universal} \prep of $G$ under $\Rel$, namely the direct sum of
all \preps of $G$ whose associated $\e(t)$'s obey the conditions listed
in $\Rel$ (some care needs to be taken here to avoid set theoretical
problems caused by considering the \stress{set of all partial
representations}).  One may then consider the \cstar-algebra
  of operators generated by the image of this \prep.  This algebra,
which we denote by $\UAlg{G}{\Rel}$, has an obvious universal property
according to which its non-degenerate *-representations are in a natural
one--to--one correspondence with the \preps of $G$ satisfying $\Rel$.

By \scite{\ELQ}{4.4}, $\UAlg{G}{\Rel}$ is isomorphic
to the crossed product of a commutative \cstar-algebra by a partial
action of $G$ which we now briefly describe.

Let $\{0,1\}^G$, also denoted by $2^G$, be given the product topology.
When convenient, we shall identify $2^G$ with the power set of $G$, that
is, the set whose elements are all the subsets of $G$.  Let $\OmegaE$ be
the compact subset of $2^G$ defined by
  $$
  \OmegaE = \{\xi \in 2^G : e\in\xi\}.
  $$

There is a canonical action of $G$ by means of partial homeomorphisms on
$\OmegaE$ defined as follows: for each $t$ in $G$, let
  $$
  \X_t = \{\xi \in \OmegaE: t\in\xi\},
  $$
  and put
  $$
  \h_t : \xi \in \X_{t\inv} \mapsto t\xi\in \X_t,
  $$
  where we interpret $t\xi$ as the set $t\xi = \{ts:s\in\xi\}$.  Then
each $\X_t$ is an open subset of $\OmegaE$ (it is also closed) and each
$\h_t$ is a homeomorphism. Also, $\(\{\X_t\}_{t\in G},\{\h_t\}_{t\in
G}\)$ is a partial action of $G$ on $\OmegaE$ according to the
definition given earlier in this section.

Let
  $C(\OmegaE)$
  be the commutative \cstar-algebra of all continuous complex valued
functions on $\OmegaE$.  For each $t$ in $G$, let $D_t$ be the ideal of
$C(\OmegaE)$ consisting of the functions vanishing off $\X_t$.  Then, by
dualization, we get the *-isomorphisms
  $$
  \theta_t : f\in D_{t\inv} \mapsto f\cdot\h_t\inv \in D_t,
  $$
  thus providing a partial action of $G$ on
  $C(\OmegaE)$.  The corresponding crossed product algebra, denoted
  $C^*_p(G)$,
  is called the \stress{partial group \cstar-algebra} of $G$
\scite{\Inverse}{6.4} and it is a model for $\UAlg{G}{\Rel}$ when $\Rel$
is the empty set of relations \scite{\Inverse}{6.5}.

In general, $\Rel$ should consist of relations of the form
  $$
  \sum \lambda_{r_1,r_2,\ldots,r_n} \e(r_1) \e(r_2) \ldots \e(r_n) = 0,
  \eqno{(\seqnumbering)}
  \label \TipicalRel
  $$
  where the sum is finite and the $\lambda$'s are complex scalars.  The
case $n=0$ is also allowed, in which case $\e(r_1) \e(r_2) \ldots
\e(r_n)$ should be interpreted as the identity operator.

Denote by $1_t$ the characteristic function
 of $\X_t$, so that $1_t\in C(\OmegaE)$.  In fact, $1_t$ belongs to
$D_t$ and it is actually the unit for this ideal.  Consider the
element $f$ of $C(\OmegaE)$ obtained by replacing each $\e(t)$ in the
left hand side of \lcite{\TipicalRel} by the corresponding $1_t$.
  If this is done for each relation in $\Rel$ we get a set of functions
on $\OmegaE$, which we denote by $\Rel'$.  The closed
subset $\GenSpec\Rel$ of $\OmegaE$, defined by
  $$
  \GenSpec\Rel = \{\xi\in\OmegaE : f(t\inv\xi) = 0, \hbox{ for all }
t\in\xi \hbox{ and } f\in\Rel'\},
  $$
  will be referred to as the \stress{spectrum} of $\Rel$.

By \scite{\ELQ}{4.4}, $\GenSpec\Rel$ is invariant by
the above partial action of $G$ on $\OmegaE$ and the crossed product of
$C(\GenSpec\Rel)$ by the restricted partial action of $G$ is isomorphic
to
  $\UAlg{G}{\Rel}$.  See \scite{\ELQ}{Section 4} for more details.

  \section{Toeplitz--Cuntz--Krieger partial isometries}
  \label \PCKPisos
  We shall now study sets of \pisos satisfying certain conditions which
are strictly weaker than the ones studied by Cuntz and Krieger in
\cite{\CKbib}.  These conditions form part of the set of conditions that
will later be used to define
  $\OA$.  Since this smaller set of conditions implies some of the
crucial properties of $\OA$, their study is relevant in itself.

Let $\Gen$ be a (possibly infinite) set and let $A$ be a function
defined on $\Gen\x\Gen$ and taking values in the set $\{0,1\}$.  We
shall also view $A$ as a matrix of zeros and ones, that is
  $A = \{\A ij\}_{i,j\in \Gen}$.  For the time being we shall make no
extra assumptions on $A$.  Later, however, we will suppose that no row
of $A$ is identically zero in order to avoid trivialities.

Our attention will be focused on families $\{S_i\}_{i\in \Gen}$ of
\pisos on a Hilbert space $H$ such that their
  initial projections $Q_i=S_i^* S_i$
  and final projections $P_j=S_j S_j^*$
  satisfy the following conditions, for all $i$ and $j$ in $\Gen$,
  \itemwidth = 40pt
  \witem{\TCKCond1)}{\commute{$Q_i$}{$Q_j$},}
  \witem{\TCKCond2)}{$P_i\perp P_j$, if $i\neq j$, and}
  \witem{\TCKCond3)}{$Q_iP_j = \A ij P_j$.}

  \medskip
  Condition \TCKCond2 can be rephrased in terms of the partial
isometries themselves:
  \itemwidth = 45pt
  \witem{\TCKCond2$'$)}{$S_i^* S_j = 0$, for $i\neq j$.}

  \medskip
  Note that \TCKCond3 means that $Q_i \geq P_j$ when $\A ij=1$ and that
$Q_i \perp P_j$ when $\A ij=0$.  Also, it is easy to see that \TCKCond3
is equivalent to
  \itemwidth = 45pt
  \witem{\TCKCond3$'$)}{$Q_i S_j = \A ij S_j$.}

  \medskip
  The reader acquainted with the work of Cuntz and Krieger will notice
that \TCKCond1--\TCKCond3 are satisfied by the partial isometries
studied in \cite{\CKbib}.  However, we will see that our conditions are
strictly weaker than the ones considered there.

Throughout this section we shall fix a family $\{S_i\}_{i\in \Gen}$ of
\pisos satisfying the above conditions.

Let $\F$ be the free group on $\Gen$ and consider the map $\pr$ from
$\F$ into the algebra $\B(H)$ of bounded operators on $H$ defined as
follows (cf. \scite{\Amena}{Section 5}): if $x$ is in $\Gen$, put
$\pr(x)=S_x$ and $\pr(x\inv)=S_x^*$.  For a general $t\in\F$, write
  $$
  t = x_1 x_2 \cdots x_n,
  $$
  in reduced form, that is,
  each $x_k\in \Gen \cup \Gen\inv$ and $x_{k+1} \neq x_k\inv$ and set
  $$
  \pr(t)=\pr(x_1)\cdots\pr(x_n).
  $$

Let us denote by $\Pos$ the \stress{positive cone} of $\F$, that is, the
unital sub-semigroup of $\F$ generated by $\Gen$.  A simple consequence
of \TCKCond2$'$ is:

  \state Proposition
  \label \PosNeg
  For every $t$ in $\F$ which is not of the form $\mu\nu\inv$, with
$\mu$ and $\nu$ in $\Pos$, we have that $\pr(t) = 0$.

Our next result is a generalization of \scite{\Amena}{5.2}.  See also
\scite{\Ortho}{5.4}.

  \state Proposition
  \label \Suffering
  The map $\pr$ defined above is an orthogonal \ss \prep of $\F$.

  \proof
  Although the present hypothesis are somewhat weaker than those of
\scite{\Amena}{5.2}, the proof follows essentially the same steps.  By
construction, it is clear that $\pr(t)\pr(s)=\pr(ts)$ if $|ts|=|t|+|s|$.
The property characterizing orthogonal \preps is also easy to check
based on \TCKCond2$'$.  Properties \lcite{\PR 2} and \lcite{\PR 3}
are trivial, so we concentrate on proving \lcite{\PR 1}.  For this we
make a series of claims, beginning with:

  \claim 1 If $n\geq1$,
  $x_1,\ldots,x_n\in\Gen$, and
  $\a=x_1\ldots x_n$ then
    $\pr(\a)^* \pr(\a) = \delta \pr(x_n)^* \pr(x_n)$,
  where
  $$
  \delta = \prod_{k=1}^{n-1}\A{x_k}{x_{k+1}}.
  $$

This is obvious if $n=1$.  Proceeding by induction, assume that $n\geq2$
and let
  $\b= x_1\ldots x_{n-1}$ and
  $\delta' = \prod_{k=1}^{n-2}\A{x_k}{x_{k+1}}$.
  Using \TCKCond3$'$, we have
  $$
  \pr(\a)^* \pr(\a) =
  \pr(\b x_n )^* \pr(\b x_n ) =
  \pr(x_n)^* \pr(\b )^* \pr(\b ) \pr(x_n) \$=
  \delta' \pr(x_n)^* \pr(x_{n-1})^* \pr(x_{n-1}) \pr(x_n) =
  \delta' \pr(x_n)^* Q_{x_{n-1}} S_{x_n} \$=
  \delta' \pr(x_n)^* \A {x_{n-1}}{x_n} S_{x_n} =
  \delta \pr(x_n)^* \pr(x_n),
  $$
  proving our claim.  Since $\delta$ is either zero or one, this also
shows that for $\a\in \Pos$ one has that $\pr(\a)^*\pr(\a)$ is
idempotent and hence that $\pr(\a)$ is a partial isometry.

  \claim 2 For every $\a$ and $\b$ in $\Pos$, if $|\a|=|\b|$, but $\a
\neq \b$ then
  $\pr(\a)^*\pr(\b) = 0$.

Let $m=|\a|=|\b|$.  If $m=1$ then $\a,\b\in\Gen$ and so, by \TCKCond2$'$
we have that
  $\pr(\a)^*\pr(\b) = 0$.  If $m > 1$, write
  $\a = \tilde \a x$ and
  $\b = \tilde \b y$ with
  $\tilde \a, \tilde \b \in \Pos$ and
  $x,y \in \Gen$.  Then
  $$
  \pr(\a)^* \pr(\b) =
  \pr(\tilde \a x)^* \pr(\tilde \b y) =
  \pr(x)^* \pr(\tilde \a)^* \pr(\tilde \b) \pr(y).
  $$
  Arguing by contradiction, if the above does not vanish then
  $ \pr(\tilde \a)^* \pr(\tilde \b) \neq 0$
  and, by induction, we have
  $\tilde \a = \tilde \b$. By claim (1) it follows that
  $\pr(\tilde \a)^* \pr(\tilde \a) = \delta \pr(z)^* \pr(z)$
  for some $z\in \Gen$, where $\delta$ is either zero or one.  Therefore
  $$
  0 \neq
  \pr(\a)^* \pr(\b) =
  \delta \pr(x)^* \pr(z)^* \pr(z) \pr(y) =
  \delta \A zy \pr(x)^* \pr(y),
  $$
  which implies that $x=y$ and hence that $\a=\b$, a contradiction.

  \claim 3 For every $\a$ and $\b$ in $\Pos$, if
  $\pr(\a)^*\pr(\b) \neq 0$ then $\a\inv \b \in \Pos \cup \Pos\inv$.

Without loss of generality assume
  $|\a| \leq |\b|$
  and write
  $\b = \tilde \b \g$ with
  $|\tilde \b| = |\a|$ and $\tilde \b, \g \in \Pos$.  Then
  $$
  0 \neq \pr(\a)^* \pr(\tilde \b \g) =
  \pr(\a)^* \pr(\tilde \b) \pr(\g),
  $$
  which implies that
  $\pr(\a)^* \pr(\tilde \b) \neq 0$ and hence, by claim (2), that
  $\a = \tilde \b$.  So $\a\inv \b = \g \in \Pos \subseteq \Pos \cup
\Pos\inv$.

Given any $t$ in $\F$, let
  $\e(t):= \pr(t) \pr(t)^*$.
  Since we already know that $\pr(\a)$ is a partial isometry for $\a\in
\Pos$, we see that $\e(\a)$ is a self-adjoint idempotent operator.

  \claim 4 For all $\a$ and $\b$ in $\Pos$ the operators
  \commute{$\e(\a)$}{$\e(\b)$}.

In the case that $\a\inv\b \not \in \Pos \cup \Pos\inv$ we have
  $$
  \e(\a) \e(\b) = \pr(\a) \pr(\a)^* \pr(\b) \pr(\b)^* = 0,
  $$
  by claim (3),
  and similarly
  $\e(\b) \e(\a) = 0$.

If, on the other hand,
  $\a\inv\b \in \Pos \cup \Pos\inv$, without loss of generality write
  $\a\inv \b = \g \in \Pos$ and note that
  $$
  \e(\a) \e(\b) =
  \pr(\a) \pr(\a)^* \pr(\a \g) \pr(\a \g)^* =
  \pr(\a) \pr(\a)^* \pr(\a) \pr(\g) \pr(\g)^* \pr(\a)^*
  \$=
  \pr(\a) \pr(\g) \pr(\g)^* \pr(\a)^* =
  \pr(\a) \pr(\g) \pr(\g)^* \pr(\a)^* \pr(\a) \pr(\a)^*
  \$=
  \pr(\b) \pr(\b)^* \pr(\a) \pr(\a)^* =
  \e(\b) \e(\a).
  $$

  \claim 5 For every $x\in\Gen$ and $\a\in\Pos$ the operators
  \commute{$Q_x$}{$\pr(\a)Q\pr(\a)^*$}, where $Q$ is either the identity
operator or the initial projection of one of the isometries $S_i$.

In case $|\a|=0$ this is a consequence of \TCKCond1.  Otherwise write
$\a=y\tilde\a$ with $y\in\Gen$ and $\tilde\a\in\Pos$ and observe that
  $$
  Q_x \pr(\a) =
  Q_x S_y \pr(\tilde\a) =
  \A xy S_y \pr(\tilde\a) =
  \A xy \pr(\a),
  $$
  which also gives $\pr(\a)^* Q_x = \A xy\pr(\a)^*$.  Therefore
  $$
  Q_x \pr(\a)Q\pr(\a)^* =
  \A xy\pr(\a)Q\pr(\a)^* =
  \pr(\a)Q\pr(\a)^* Q_x,
  $$
  proving our claim.

  \claim 6 For all $t$ and $s$ in $\F$ the operators $\e(t)$ and $\e(s)$
commute.

By \lcite{\PosNeg}, we may assume that $t=\mu\nu\inv$ and $s=\a\b\inv$,
where $\mu,\nu,\a$, and $\b$ are positive.
  We may clearly also assume that
  $|t|=|\mu| + |\nu|$
  and that
  $|s|=|\a| + |\b|$.
  Then
  $$
  \e(t)=\pr(\mu)\pr(\nu)^*\pr(\nu)\pr(\mu)^*
  \and
  \e(s)=\pr(\a)\pr(\b)^*\pr(\b)\pr(\a)^*.
  $$
  Now, observe that the product of $\e(t)$ and $\e(s)$, in either order,
will vanish if $\pr(\mu)^*\pr(\a)=0$, which is the case, unless
  $\mu\inv\a \in \Pos \cup \Pos\inv$, according to claim (3).  So, let
us suppose, without loss of generality, that $\mu\inv\a=\g$ for some
$\g\in\Pos$.

Let us use the notation
  $Q=\pr(\nu)^*\pr(\nu)$, observing that, in case $|\nu|>0$, $Q$ is
either zero, or one of the $Q_x$'s, by claim (1), or, if $\nu=e$,
$Q=I$. The same observation applies to
  $Q'=\pr(\b)^*\pr(\b)$ and
  $Q''=\pr(\mu)^*\pr(\mu)$.

We then have
  $$
  \e(t)\e(s) =
  \pr(\mu)\pr(\nu)^*\pr(\nu)\pr(\mu)^*
  \pr(\mu)\pr(\g)\pr(\b)^*\pr(\b)\pr(\g)^*\pr(\mu^*) \$=
  \pr(\mu)QQ''\pr(\g)Q'\pr(\g)^*\pr(\mu)^*.
  $$
  Using \TCKCond1 and claim (5), we conclude that this coincides with
  $$
  \pr(\mu)
  \pr(\g)Q'\pr(\g)^*
  Q''
  Q
  \pr(\mu)^*
  \$=
  \pr(\mu)
  \pr(\g) \pr(\b)^*\pr(\b) \pr(\g)^*
  \pr(\mu)^*\pr(\mu)
  \pr(\nu)^*\pr(\nu)
  \pr(\mu)^*
  =
  \e(s)\e(t).
  $$

This proves our last claim and we are now ready to prove \PR 1, that is
  $$
  \pr(t) \pr(s) \pr(s)^* = \pr(ts) \pr(s)^*,
  $$
  for $t,s \in \F$.
  To do this we use induction on
  $|t| + |s|$.  If either $|t|$ or $|s|$ is zero, there is nothing to
prove.  So, write
  $t = \tilde t x$ and
  $s = y \tilde s$,
  where $x,y\in \Gen \cup \Gen\inv$ and, moreover,
  $|t|=|\tilde t|+1$ and
  $|s|=|\tilde s|+1$.

In case
  $x\inv \neq y$ we have
  $|ts|=|t| + |s|$ and hence $\pr(ts) = \pr(t) \pr(s)$.
  If, on the other hand,
  $x\inv = y$, we have
  $$
  \pr(t) \pr(s) \pr(s)^* =
  \pr(\tilde t x) \pr(x\inv \tilde s) \pr(\tilde s\inv x) =
  \pr(\tilde t) \pr(x) \pr(x)^* \pr(\tilde s) \pr(\tilde s)^* \pr(x).
  $$
  With an application of claim (6) and the use of the induction
hypothesis we conclude that the above equals
  $$
  \pr(\tilde t) \pr(\tilde s) \pr(\tilde s)^* \pr(x) \pr(x)^* \pr(x) =
  \pr(\tilde t \tilde s) \pr(\tilde s)^* \pr(x)
  \$=
  \pr(ts) \pr(\tilde s\inv x) =
  \pr(ts) \pr(s)^*.
  \proofend
  $$

  \section{The Toeplitz--Cuntz--Krieger algebra}
  \label \TCKAlgSec
  Let $\LA$ be the universal unital \cstar-algebra generated by a family
$\{S_i\}_{i\in \Gen}$ of \pisos satisfying \TCKCond1--\TCKCond3 of
Section {\PCKPisos}.  See \cite{\Blackadar} for a definition of
universal \cstar-algebras given by generators and relations.

  \state Proposition
  \label \ParcRepAlg
  The non-degenerate *-representations of $\LA$ are in a one--to--one
correspondence with the \preps $\pr$ of $\F$ such that
  \izitem $\e(x)\e(y) = 0$, for $x\neq y$ in $\Gen$,
  \zitem $\e(x\inv) \e(y) = \A xy \e(y)$, for all $x$ and $y$ in $\Gen$,
and
  \zitem $\e(ts) \leq \e(t)$, whenever $t$ and $s$ are in $\F$ and
satisfy
  $|ts|=|t|+|s|$,
  \smallskip\noindent where, as usual, $\e(t)=\pr(t)\pr(t)^*$.

  \proof
  Given any partial representation $\pr$ of $\F$ satisfying the
conditions above let, for each $x\in\Gen$, $S_x=\pr(x)$.  Then, since
the final projections $\e(t)$ associated to any \prep commute with each
other, we see that $\{S_x\}_{x\in\Gen}$ satisfies \lcite{\TCKCond1}.
Conditions \lcite{\TCKCond2} and \lcite{\TCKCond3} follow from (i) and
(ii), respectively.  Therefore we get a non-degenerate representation
$\pi$ of $\LA$ by the universal property.

  Conversely, given a representation $\pi$ of $\LA$, consider the
collection of partial isometries $\{S_x\}_{x\in\Gen}$ obtained as the
image, under $\pi$, of the canonical ones in $\LA$.  Using
\lcite{\Suffering} we get a partial representation $\pr$ of $\F$ such
that $\pr(x)=S_x$.  It is clear that $\pr$ will satisfy (i) and (ii).
Also, since $\pr$ is semi-saturated, (iii) will hold by
\scite{\Amena}{5.4}.
  It is now easy to see that the correspondences described are
inverses of each other.
  \proofend

We therefore see that $\LA$ is isomorphic to $\UAlg{\F}{\RelTCK}$, for
the set $\RelTCK$ consisting of the relations
  \izitem $\e(x)\e(y) = 0$, for $x\neq y$ in $\Gen$,
  \zitem $\e(x\inv) \e(y) - \A xy \e(y) = 0$, for $x, y$ in $\Gen$, and
  \zitem $\e(ts)\e(t) - \e(ts) = 0$, for all $t,s\in\F$ such that
$|ts|=|t|+|s|$.

  \medskip\noindent
  The reason we use the superscript in $\RelTCK$ is that we will later
introduce a larger set of relations associated to $A$.

By \scite{\ELQ}{4.4} we conclude that $\LA$ is
isomorphic to the crossed product $C(\GenSpec{\RelTCK}) \crossproduct
\F$.  Our next major goal is to describe the partial dynamical system
giving rise to $\LA$ and, in particular, to characterize the spectrum
$\GenSpec{\RelTCK}$ of these relations, which we denote simply by
$\specTCK$.

  \medskip\noindent
  The corresponding collection ${\RelTCK}'$ of functions on $\OmegaE$,
as described in Section \lcite{\Preliminaries}, is the union of
  \medskip
  $$\displaystyle
  \matrix{
  \Rel'_1 &=& \{1_x1_y: x,y\in\Gen, x\neq y\},\hfill\cr\cr
  \Rel'_2 &=& \{1_{x\inv} 1_y - \A xy 1_y: x,y\in\Gen\},
  \hbox{ and}\hfill\cr\cr
  \Rel'_3 &=& \{1_{ts}1_t - 1_{ts}:t,s\in\F,\quad |ts|=|t|+|s|\}.
\hfill
  }
  $$

  \medskip
  We would now like to describe, in terms as explicit as possible, the
properties that a given $\xi$ in $\OmegaE$ should have in order to
belong to $\specTCK$.  It is convenient to study these properties
separately for $\Rel'_1$, $\Rel'_2$ and $\Rel'_3$.  Accordingly, let
  $$
  \Omega_i = \{\xi\in\OmegaE : f(t\inv\xi) = 0, \hbox{ for all } t\in\xi
\hbox{ and } f\in\Rel_i'\} \for i=1,2,3,
  $$
  and observe that
  $\specTCK = \Omega_1\cap\Omega_2\cap\Omega_3$.

  \state Proposition
  \label \CondOne
  Let $\xi\in\OmegaE$.  A necessary and sufficient condition for
$\xi$ to belong to $\Omega_1$ is that, whenever $t\in\xi$, there exists
at most one $x\in\Gen$ such that $tx\in\xi$.

  \proof
  Since $1_t$ is the characteristic function of
  $\X_t = \{\xi \in \OmegaE: t\in\xi\},$ we have that
  $$
  1_t(\xi) = [t\in\xi],
  $$
  where the brackets correspond to the obvious boolean valued function,
taking values in the set $\{0,1\}$, seen as a subset of complex numbers.
Thus, if
  $f=1_x1_y$
  and $t\in\xi$, we have
  $$
  f(t\inv\xi) =
  1_x1_y(t\inv\xi) =
  [x\in t\inv\xi] [y\in t\inv\xi] =
  [tx\in \xi] [ty\in \xi].
  $$
  So, to say that $f(t\inv\xi)=0$ is to say that either $tx\notin\xi$ or
$ty\notin\xi$.
  \proofend

  \state Proposition
  \label \CondTwo
  Let $\xi\in\OmegaE$.  A necessary and sufficient condition for $\xi$
to belong to $\Omega_2$ is that, whenever
  $t\in\F$ and $y\in\Gen$ are such that both $t$ and $ty$ belong to
$\xi$ then
  \iBitem $\xi$ contains
  $tx\inv$ for all $x\in\Gen$ for which $\A xy = 1$, and
  \Bitem $\xi$ does not contain those
  $tx\inv$ for which $\A xy = 0$.

  \proof
  The condition mentioned is equivalent to
  $$
  \forall x,y\in\Gen, ~\forall t\in\xi, ~
  ty\in\xi \imply (tx\inv\in\xi \bimply \A xy = 1).
  \eqno{(\dagger)}
  $$
  Let $\phi$ and $\psi$ be logical propositions.  Then the boolean value
of ``$\phi\imply\psi$'' is given by
  $$
  [\phi\imply\psi] =
  [\neg \phi \vee \psi] =
  1 - [\phi \wedge \neg \psi] =
  1 - [\phi] (1 - [\psi]).
  $$
  Therefore
  $$
  [\phi\bimply\psi] =
  [(\phi\imply\psi)\wedge(\psi\imply\phi)] =
  [\phi\imply\psi] [\psi\imply\phi] \$=
  \Big(1 - [\phi] (1 - [\psi]) \Big)
  \Big(1 - [\psi] (1 - [\phi]) \Big) =
  1 - \Big( [\phi]-[\psi] \Big)^2.
  $$
  So, the value of the logical proposition after the quantifiers in
($\dagger$) is
  $$
  \Big[ty\IN \xi \imply ( tx\INVIN\xi \bimply \A xy=1)\Big] =
  1 - [ty\IN\xi] \Big( 1- \Big[tx\INVIN\xi \bimply \A xy=1\Big] \Big)
\$=
  1 - [ty\IN\xi] \Big( [tx\INVIN\xi] - \A xy \Big)^2.
  \eqno{(\ddagger)}
  $$
  On the other hand,
  if $f=1_{x\inv} 1_y - \A xy 1_y$ and $t\in\xi$ then
  $$
  f(t\inv\xi) =
  [x\INVIN t\inv\xi] [y\IN t\inv\xi] - \A xy [y\IN t\inv\xi] =
  [ty\IN\xi] \Big( [tx\INVIN\xi] - \A xy \Big).
  $$
  It is now easy to see that
  $f(t\inv\xi) = 0$ if and only if ($\ddagger$) equals one.
  \proofend

It is often convenient to view the free group $\F$ as the set of
vertices of its Cayley graph, in which one draws an (unoriented) edge
{}from $t$ to $s$ if and only if $t\inv s \in \Gen \cup \Gen\inv$.  Thus a
typical edge joins $t$ to $tx$ if $x\in\Gen\cup\Gen\inv$.
  Since the Cayley graph of $\F$ is a connected tree, given $t,s\in\F$,
there exists a unique shortest path joining $t$ and $s$.  It is easy to
see that the length of this path (number of edges) is given by $|t\inv
s|$.  In fact, if
  $t\inv s=x_1 x_2\ldots x_n$ is in reduced form then this path is given
by
  $$
  (t\ ,\ tx_1\ ,\ tx_1x_2\ ,\ \ldots\ ,\ tx_1 x_2\ldots x_n=s).
  $$
  It is easy to see that a vertex $r$ belongs to this path if and only
if
  $
  |t\inv s| = |t\inv r| + |r\inv s|.
  $

  \definition
  A subset $\xi\subseteq\F$ is said to be \stress{convex} if, whenever
$t,s\in\xi$, the whole shortest path joining $t$ and $s$ lies in $\xi$.
Equivalently, $\xi$ is convex if and only if
  $$(t,s\in\xi)\ \wedge\ \(|t\inv s| = |t\inv r| + |r\inv s|\)
  \ \Longrightarrow\ r\in\xi.
  $$

  \state Proposition
  \label \CondSS
  Let $\xi\in\OmegaE$.  A necessary and sufficient condition for $\xi$
to belong to $\Omega_3$ is that $\xi$ be convex.

  \proof
  Let $s,r\in\F$ with $|sr|=|s|+|r|$ and let
  $f=1_{sr}1_s - 1_{sr}$.  Then, for $t\in\xi$,
  $$
  f(t\inv\xi) =
  [sr\IN t\inv\xi] [s\IN t\inv\xi] - [sr\IN t\inv\xi] =
  [tsr\IN\xi]([ts\IN\xi]-1) \$=
  1 - [tsr\IN\xi](1-[ts\IN\xi]) - 1 =
  [tsr\IN\xi \imply ts\IN\xi] - 1.
  $$
  Thus we conclude that $f(t\inv\xi) = 0$ if and only if $tsr\in\xi
\imply ts\in\xi$.  This holds for all $t\in\xi$ and $s,r$ satisfying
$|sr|=|s|+|r|$ if and only if $\xi$ is convex.
  \proofend

We thus see that a given $\xi\in\OmegaE$ belongs to $\specTCK$ if and
only if it satisfies the conditions described in our last three
Propositions, namely
  \lcite{\CondOne}, \lcite{\CondTwo} and \lcite{\CondSS}.
  This concludes the description of $\specTCK$ and hence we obtain, from
\scite{\ELQ}{4.4}, the following characterization
of $\LA$:

  \state Theorem
  \label \TheoremOnLA
  Let
  $\Gen$ be any set,
  $A = \{\A ij\}_{i,j\in \Gen}$ be a 0--1 matrix and
  $\F$ be the free group on $\Gen$.  Let $\LA$ be the universal unital
\cstar-algebra generated by a family $\{S_i\}_{i\in \Gen}$ of \pisos,
satisfying the conditions
  \izitem \commute{$S_i^*S_i$}{$S_j^*S_j$},
  \zitem $S_i^*S_j=0$, for $i\neq j$, and
  \zitem $S_i^*S_iS_j = \A ij S_j$,
  \medskip \noindent
  for all $i$ and $j$ in $\Gen$.
  Let $\specTCK$ be the compact subspace of $2^\F$ given by
  $$
  \specTCK =
  \left\{
  \matrix{ \xi\in 2^\F : & e\in\xi,
  \hbox{ $\xi$ is convex,}\hfill\cr
  &\hbox{if $t\in\xi$ there is at most one $x\in\Gen$ such that
$tx\in\xi$, and}\hfill\cr
  &\hbox{if $t\in\xi$, $y\in\Gen$, and $ty\in\xi$ then
  $tx\inv\in\xi \bimply \A xy = 1$}\hfill\cr
  }\right\},
  $$
  and consider the partial action of $\F$ on
  $C(\specTCK)$ induced by the partial homeomorphisms
  $$
  \h_t:\xi\in\X_{t\inv} \ \mapsto\ t\xi\in \X_t \for t\in\F,
  $$
  where
  $\X_t=\{\xi\in\specTCK:t\in\xi\}$.
  For $x\in\Gen$ let $1_x$ be the characteristic function of $\X_x$.
Then the correspondence
  $S_x \mapsto 1_x\delta_x$
  extends to an isomorphism of $\LA$ onto
  $C(\specTCK)\crossproduct \F$.

  \section{Stems and roots}
  The characterization of $\specTCK$ obtained in the last section does
not lend itself to the study of the dynamical properties of our partial
action.  For this reason we will need to obtain a more convenient one.
The first important ingredient of the new characterization will be the
intersection of an element $\xi\in\specTCK$ with the positive cone
$\Pos$.  This will provide a link with the familiar path space.  The
second ingredient is introduced in \fcite{5.6} and reveals an
unexpected subtlety concerning finite paths.

  \definition
  \label \WordDef
  \Bitem By a \stress{word} we shall mean a finite or infinite sequence
$\w = (x_1,x_2,\ldots)$, where each $x_i\in\Gen$.
  \Bitem The \stress{length} of $\w$, denoted $|\w|$, is the (possibly
infinite) number of coordinates of $\w$.
  \Bitem For each $i$, we shall denote by $\w_i$ the $i^{th}$ coordinate
of $\w$, so that
  $\w = (\w_1,\w_2,\ldots)$.
  \Bitem Given a word $\w$ and an integer $n$, with $0\leq n < |\w|$, we
shall denote by
  $\trunc{\w}{n}$ the \stress{sub-word} of $\w$ defined by
  $\trunc{\w}{n} = (\w_1,\w_2,\ldots,\w_n)$. Note that only initial
segments are considered.

If $\w$ is a finite word we will often identify $\w$ with the product of
its coordinates, namely
  $\w_1 \w_2\ldots \w_{|\w|} \in \Pos$.
  The empty word will be identified with the unit group element $e$.  In
this way, the set of finite words may be identified with $\Pos$.  This
identification will be made, from now on, without further warning.

  \definition
  \label \SetWDef
  Given a word $\w$ we shall denote by $\set\w$ the subset of $\Pos$
consisting of the group elements associated to the finite sub-words of
$\w$.  In symbols,
  $\set\w=\{e,\w_1,\w_1\w_2,\w_1\w_2\w_3,\ldots\}$.

There is a natural order structure in $\F$ defined by
  $$
  t\leq s \bimply t\inv s\in\Pos.
  $$
  The subsets of the form $\set\w$ are precisely
  the hereditary directed subsets of $\Pos$ (cf. \cite{\Nica,\Laca}).
In the specific case of a finite word $\w$, we have that
$\set\w=\{t\in\F: e\leq t\leq \w\}$.

{}From now on
  $A = \{\A ij\}_{i,j\in\Gen}$
  will be a fixed 0--1 matrix indexed over an arbitrary set $\Gen$.

  \definition
  \label \AdmissibleWord
  A word $\w$ is said to be \stress{admissible}
  if
  $\A{\w_i}{\w_{i+1}}=1$, for all $i$.  If $|\w|<2$ then $\w$ is
admissible by default.

The following is an important tool in the characterization of
$\specTCK$:

  \state Proposition
  \label \TheWord
  Let $\xi\in\specTCK$.  Then there exists a unique admissible word
$\stem\xi$ such that
  $$
  \xi\cap\Pos =
  \set{\stem\xi}.
  $$

  \proof
  Follows easily from the fact that $\xi$ satisfies the conditions of
  \lcite{\CondOne}, \lcite{\CondTwo} and \lcite{\CondSS}.
  \proofend

  \definition
  \label \BoundedDef
  Let $\xi\in\specTCK$.  The word $\stem\xi$ mentioned in
\lcite{\TheWord} will be called the \stress{stem} of $\xi$.  We will say
that $\xi$ is \stress{bounded} if its stem $\stem\xi$ is of finite
length.  Otherwise we will say that $\xi$ is \stress{unbounded}.

  \definition
  \label \RootDef
  Let $\xi\in\specTCK$ and let $t\in\xi$. The \stress{root of $t$
relative to $\xi$}, denoted $\root t\xi$, is the subset of $\Gen$
defined by
  $
  \root t\xi = \{x\in\Gen : tx\inv\in\xi \}.
  $

  \sysstate{Remark}{\rm}
  {\label \RootIsColumn
  Suppose that $\xi$ is in $\specTCK$ and let $t\in\xi$.  We already
know from \lcite{\CondOne} that there is at most one $y\in\Gen$ such
that $ty\in\xi$.  Suppose, for a moment, that such a $y$ exists.  Then,
by \lcite{\CondTwo}, we have that
  $tx\inv\in\xi$ if and only if $\A xy = 1$.  Therefore the root of $t$
relative to $\xi$ is given simply by $\root t\xi = \{x\in\Gen: \A
xy=1\}$.  However, if no such $y$ exists, we will see shortly that
$\root t\xi$ can be almost any subset of $\Gen$.}

The relevance of stems and roots is that, together, they characterize
$\xi$, as we shall now show.

  \state Proposition
  \label \MembershipCriteria
  Each $\xi\in\specTCK$ consists, precisely, of the elements $t$ in $\F$
of the form $t=\a\b\inv$ such that
  \iBitem $|t|=|\a|+|\b|$,
  \Bitem $\a$ is a finite sub-word of $\stem\xi$, and
  \Bitem $\b$ is a finite admissible word such that, either $\b=e$, or
$\last\b\in\root\a\xi$.


  \proof
  Let $t\in\xi$ and write $t = x_1 x_2 \cdots x_n$ in reduced form.
Since every $\xi$ in $\specTCK$ is convex, we have that
  $x_1 x_2 \cdots x_i \in \xi,$
  for all $i=1,\ldots,n$.

We claim that there is no $i$ for which $x_i\in\Gen\inv$ and
$x_{i+1}\in\Gen$.  Otherwise, letting
  $s = x_1 x_2 \cdots x_i$, $x=x_i\inv$, and $y=x_{i+1}$,
  we would have that $s\in\xi$, that $x,y\in\Gen$, and that both $sx$
and $sy$ belong to $\xi$, contradicting \lcite{\CondOne}.

This proves our claim and hence that, in the reduced form of $t$, all
elements from $\Gen$ must be to the left of the elements from
$\Gen\inv$.  That is,
  $t = x_1 \ldots x_k y_\ell\inv \ldots y_1\inv$, where every
  $x_i$ and $y_j$ belong to $\Gen$.  Putting $\a=x_1 \ldots x_k$ and
  $\b=y_1\ldots y_\ell$ we have that $t=\a\b\inv$ and that
$|t|=|\a|+|\b|$.

Since $\a$ is in $\xi\cap\Pos$, it is a sub-word of $\stem\xi$, by
\lcite{\TheWord}.

If $|\b|\geq 2$, using the convexity of $\xi$, we have, for all $i=
2,\ldots,\ell$, that
  $s:= \a y_\ell\inv \ldots y_i\inv$ belongs to $\xi$, as well as
  $s y_{i-1}\inv$
  and
  $s y_i$.  Therefore, by \lcite{\CondTwo}, we have that
$\A{y_{i-1}}{y_i} = 1$.  This shows that $\b$ is an admissible word.
Finally, since $\a y_\ell\inv\in\xi$, we have that
  $\last\b = y_\ell\in\root\a\xi$.  This concludes the proof that every
$t$ in $\xi$ is of the form mentioned in the statement.

Conversely, suppose that $t=\a\b\inv$ is as above. Then
$\a\in\set{\stem\xi} = \xi\cap\Pos$ and hence $\a\in\xi$.

Assuming now that $\b\neq e$, write $\b=y_1\ldots y_\ell$ with
$y_i\in\Gen$.  Then $\last\b=y_\ell$ and hence, by hypothesis, $\a
y_\ell\inv\in\xi$.
  We shall now prove that $\a y_\ell\inv \ldots y_k\inv \in\xi$, for all
$k$, using (backwards) induction on $k$.  We thus suppose that $k \leq
\ell$ and that
  $s:=\a y_\ell\inv \ldots y_k\inv \in\xi$. Observing that $s y_k\in
\xi$, we then conclude, by \lcite{\CondTwo},
  that $sz\inv \in \xi$ for all $z\in\Gen$ for which $\A z{y_k} =1$.
Since $\b$ is an admissible word, and hence $\A {y_{k-1}}{y_k} =1$, we
then obtain the desired induction step, that is, $\a y_\ell\inv \ldots
y_{k-1}\inv \in\xi$. It follows that $\a\b\inv\in\xi$.
  \proofend

The following consequence of our last result gives a simple
characterization of the bounded elements:

  \state Corollary
  \label \FormOfBounded
  Let $\xi\in\specTCK$ be bounded and put $\w=\stem\xi$ and
  $R = \root\w\xi$.  Then
  $$
  \xi = \{
  \w \b\inv : \b=e,
  \hbox{\rm\ or $\b$ is a finite admissible word such that }
  \last\b\in R\}.$$

  \proof
  Left to the reader.
  \proofend

Let us again consider the free group $\F$ as the set of vertices of its
Cayley tree $T$.  Recall that the edges emanating from a given $t$ in
$\F$ are exactly those of the form $(t,s)$ where
  $t\inv s\in\Gen\cup\Gen\inv$.  If we remove such an edge from $T$ we
are left with two connected components, one of which contains $t$, the
other containing $s$.  Let us temporarily denote these components by
$T_t$ and $T_s$, respectively.

  \definition
  \label \Subtree
  Let $t,s\in\F$ be such that $t\inv s\in\Gen\cup\Gen\inv$.  The set of
elements of $\F$ corresponding to the vertices of $T_t$ will be denoted
by $\STree{t}{s}$.

If $r\in\STree{t}{s}$ then the shortest path joining $r$ and $s$ cannot
stay within $T_t$ and hence must pass through the edge $(t,s)$.  It
follows that the distance from $r$ to $s$, namely $|r\inv s|$, is bigger
than $|r\inv t|$ by one.  In other words
  $$
  \STree{t}{s} = \{r\in\F: |r\inv s| = |r\inv t| +1 \}.
  $$
  We also have that $\F$ is the disjoint union of $\STree ts$ and
$\STree st$.

  \state Proposition
  \label \Edge
  Let $t\in\F$ and $x\in\Gen$.  Also let $\xi_1$ and $\xi_2$ be in
$\specTCK$ and suppose that both $t$ and $tx$ belong to both $\xi_1$ and
$\xi_2$.  Then
  $\xi_1\cap \STree{t}{tx} = \xi_2\cap \STree{t}{tx}$.

  \proof
  Replacing $\xi_i$ by $t\inv\xi_i$ we may assume that $t=e$.  We begin
by showing that
  $\xi_1\cap \STree{e}{x} \subseteq \xi_2\cap \STree{e}{x}$. Let
$s\in\xi_1\cap \STree{e}{x}$ and, using \lcite{\MembershipCriteria},
write $s=\a\b\inv$ where $\a$ is a sub-word of $\stem{\xi_1}$, $\b$ is
an admissible word, and $\b=e$, or $\last\b\in\root\a{\xi_1}$.  Now,
note that, since $x\in\xi_1$, we have that
  $\stem{\xi_1}$ must start with $x$.  Therefore, either $\a=e$ or $\a_1
= x$.  However, this last possibility is ruled out by the fact that
$s\in \STree{e}{x}$, so we conclude that $\a=e$.  Since $x\in\xi_1$
then, by \lcite{\CondTwo}, we have that
  $\root e{\xi_1} = \{y\in\Gen:\A yx=1\}
  = \root e{\xi_2}$ and we see that
  $\A{\last\b}x =1$.

We conclude that $s=\b\inv$, where either $\b=e$ or $\b$ is an
admissible word such that $\A{\last\b}x =1$.  Now, under the light of
\lcite{\MembershipCriteria}, this implies that $s\in\xi_2$.  By
symmetry, we also have that
  $\xi_2\cap \STree{e}{x} \subseteq \xi_1$ and hence the proof is
complete.
  \proofend

We now present a uniqueness criteria for elements of $\specTCK$ in terms
of stems and roots.

  \state Proposition
  \label \UniqueWordRoot
  Suppose $\xi_1$ and $\xi_2$ are elements of $\specTCK$ such that
$\stem{\xi_1} = \stem{\xi_2}$.  Then
  \iBitem if $\stem{\xi_1}$ is infinite then $\xi_1=\xi_2$,
  \Bitem if $\stem{\xi_1}$ is finite and
  $\root{\stem{\xi_1}}{\xi_1} = \root{\stem{\xi_2}}{\xi_2}$ then
$\xi_1=\xi_2$.

  \proof
  Let us start with the infinite case, denoting by $\w$ the common value
of $\stem{\xi_1}$ and $\stem{\xi_2}$.  For each positive integer $n$
consider the subset
  $F_n = \STree{\trunc{\w}{n}}{\trunc{\w}{n+1}}$ of $\F$, as in
\lcite{\Subtree}.  It is easy to see that
  $\F = \bigcup_{n=1}^\infty F_n$ and hence that
  $$
  \xi_i = \bigcup_{n=1}^\infty F_n\cap\xi_i.
  $$
  Since both $\xi_1$ and $\xi_2$ contain $\trunc{\w}{n}$ as well as
  $\trunc{\w}{n+1} = \trunc{\w}{n}\w_{n+1}$
  we have, by \lcite{\Edge}, that
  $F_n\cap\xi_1 = F_n\cap\xi_2$ and hence that $\xi_1=\xi_2$.
  The finite case follows immediately from Corollary
\lcite{\FormOfBounded}.
  \proofend

This takes care of uniqueness.  As for existence, let us begin with the
infinite case.

  \state Proposition
  \label \ExistenceInfiniteWord
  Let $\w$ be an infinite admissible word.  Then there exists a
(necessarily unique) $\xi$ in $\specTCK$ such that $\stem\xi=\w$.

  \proof
  Define $\xi$ to be the subset of\/ $\F$ formed by the elements $t$ of
the form $t=\a\b\inv$, where $|t|=|\a|+|\b|$, $\a$ is a finite sub-word
of $\w$ and $\b$ is a finite admissible word such that $\b=e$ or
  $\A{\last\b}{\w_{|\a|+1}} = 1$. Clearly $e\in\xi$ and hence
$\xi\in\OmegaE.$

Recall that
  $\specTCK = \Omega_1 \cap \Omega_2 \cap \Omega_3$,
  where the $\Omega_i$ were defined right before \lcite{\CondOne}.  In
order to prove that $\xi\in\specTCK$ we must therefore verify that $\xi$
belongs to each $\Omega_i$.

Obviously $\xi$ is convex and hence $\xi\in\Omega_3$ by \lcite{\CondSS}.
Suppose that $t\in\xi$ and write $t=\a\b\inv$ as above. Let $x\in\Gen$.
Assuming that $\b\neq e$, observe that, in order for $tx=\a\b\inv x$ to
belong to $\xi$, $x$ must coincide with $\b_1$.  If $\b=e$, that is, if
$t=\a$ then clearly $x =\w_{|\a|+1}$ is the only element in $\Gen$ for
which $tx\in\xi$.  Hence $\xi\in\Omega_1$ by \lcite{\CondOne}.

Suppose $t\in\xi$ and let $y\in\Gen$ be such that $ty\in\xi$.  Again,
write $t=\a\b\inv$.

If $\b\neq e$ we have already seen that $y=\b_1$.  So, for every
$x\in\Gen$ such that $\A xy = 1$, we have that $x\b$ is admissible and
hence
  $tx\inv = \a(x\b)\inv$ lies in $\xi$.  When $\A xy =0$ the reduced
form of $tx\inv$ is
  $\a\b\inv x\inv$, which shows that $tx\inv$ is not in $\xi$, since
$x\b$ is not admissible.

    If $\b=e$ then $ty=\a y$, so that $\a y$ is also a sub-word of $\w$
and thus $y=\w_{|\a|+1}$.  We must now ask ourselves for which $x$ in
$\Gen$ will we have $\a x\inv\in\xi$.  By definition of $\xi$ this is so
exactly when $\A xy = 1$ and hence we see that $\xi\in\Omega_2$ by
\lcite{\CondTwo}.
  Finally, it is easy to see that $\stem\xi=\w$.
  \proofend

The existence of $\xi \in\specTCK$ with a prescribed finite stem is
studied next.

  \state Proposition
  \label \ExistenceFiniteWord
  Let $\w$ be a finite admissible word and let $R$ be a subset of $\Gen$
which, in case $|\w|>0$, contains $\last\w$.  Then there exists a
(necessarily unique) $\xi$ in $\specTCK$ such that $\stem\xi=\w$ and
$\root\w\xi = R$.

  \proof
  Define
  $$
  \xi = \{
  \w \b\inv : \b=e,
  \hbox{\rm\ or $\b$ is a finite admissible word such that }
  \last\b\in R\}.$$
  Much as in the previous Proposition, it follows that $\xi$ is in
$\specTCK$, that $\stem\xi=\w$ and that $\root\w\xi=R$.
  \proofend


The new characterization of $\specTCK$, promised at the beginning of this
section, may then be summarized as follows:

  \iBitem
  The correspondence $\xi\mapsto\stem\xi$
  is a bijection between the sets

  \newdimen \ldim
  \newdimen \rdim
  \newdimen \mdim
  \mdim = 2em
  \rdim = 0.55\hsize
  \ldim = \hsize
  \advance \ldim by -\mdim
  \advance \ldim by -\rdim
  \def\fix#1#2{\bigskip\hbox{
               \hbox to \ldim{\hfil$#1$}
               \hbox to \mdim{\hfil$\longrightarrow$\hfil}
               \hbox to \rdim{$#2$\hfil}}\bigskip}

  \fix{\Big\{\hbox{unbounded elements of $\specTCK$}\Big\}}
      {\Big\{\hbox{infinite admissible words}\Big\}.}
  \Bitem
  The correspondence
  $\xi\mapsto \Big(\stem\xi,\root{\stem\xi}\xi\Big)$
  is a bijection between the sets

  \fix{\Big\{\hbox{bounded elements of $\specTCK$}\Big\}}
  {\left\{
  \matrix{(\w,R) : & \w \hbox{ is a finite admissible word},\hfill\cr
                   & R\subseteq\Gen, \hbox{ and }
                     \last\w\in R \hbox{ if } |\w|>0\hfill\cr}
  \right\}.}

  \sysstate{Remark}{\rm}
  {\label \SmallestGuy
  It is interesting to observe that there are many elements $\xi$ in
$\specTCK$ for which $\stem\xi$ is the empty word.  In fact, there is
one such $\xi$ for each subset $R$ of $\Gen$, as seen above.  In
particular, if $R$ is the empty set then $\xi=\{e\}$.  We will denote
this element of $\specTCK$ by $\phi$.}

  \section{The topology of $\specTCK$}
  By definition, $\specTCK$ is a closed subset of $\OmegaE$ and hence a
compact topological space.  We would now like to describe the topology
of $\specTCK$ by exhibiting a fundamental system of neighborhoods of
each $\xi$ in $\specTCK$ in terms of the characterization of elements by
stems and roots.  The unbounded case is treated first.

  \state Proposition
  \label \NBDforUnbdd
  Let $\xi\in\specTCK$ and assume that $\w=\stem{\xi}$ is infinite.
  Define
  $$
  V_n= \{\eta\in\specTCK : \trunc{\w}{n}\in\eta\}.
  $$
  Then the collection $\{V_n\}_{n\in\N}$ forms a fundamental system of
open neighborhoods of $\xi$ in $\specTCK$.

  \proof
  Since $\specTCK$ is a topological subspace of $2^\F$, it is clear that
each $V_n$ is a relatively open set, which obviously includes $\xi$.

Let $U$ be a neighborhood of $\xi$.  Then by the definition of the
product topology, $U$ contains a neighborhood of $\xi$ of the form
  $$
  W = \{\eta\in\specTCK :
  t_1,\ldots,t_p \in \eta,~s_1,\ldots,s_q\notin\eta
  \},
  $$
  where $t_1,\ldots,t_p$ and $s_1,\ldots,s_q$ are elements of $\F$.
Since $W$ contains $\xi$, the $t_i$'s are in $\xi$, while the $s_j$'s
are not.

Let $n$ be an integer bigger than $|t_i|$ and $|s_j|$ for all $i$ and
$j$.
  We will check that $V_n\subseteq W$, from which the proof will be
concluded.  Suppose $\eta\in V_n$, that is,
  $\trunc{\w}{n}\in\eta$.  So
  $$
  \trunc{\w}{n} \in
  \eta\cap\Pos =
  \set{\stem\eta},
  $$
  from which it follows that $\stem\xi$ and $\stem\eta$ agree up to the
$n^{th}$ coordinate.  Let us use the membership criteria provided by
\lcite{\MembershipCriteria}, in order to decide whether or not the $t_i$
and the $s_j$ belong to $\eta$.  In doing so we observe that the $\a$
mentioned there is necessarily shorter than $\trunc{\w}{n}$, since we
took $n$ bigger than $|t_i|$ and $|s_j|$.  The outcome of this criteria
will therefore be the same as if we tested membership to $\xi$.  It
follows that $t_i\in\eta$ and that $s_j\notin\eta$.  Thus $\eta\in W$,
showing that $V_n \subseteq W \subseteq U$.
  \proofend

The characterization of neighborhoods of bounded elements of $\specTCK$
is a little more involved.

  \state Proposition
  \label \NBDforBounded
  Let $\xi\in\specTCK$ be bounded (as defined in \lcite{\BoundedDef})
and let $\w$ be the stem of $\xi$.  Given finite subsets $X$, $Y$, and
$Z$ of\/ $\Gen$, with
  $X\subseteq\root\w\xi$ and $Y\cap\root\w\xi=\0$,
  define
  $$
  \W XYZ =
  \left\{\matrix{
  \eta\in\specTCK : & \w\hfill\in\eta,\cr
  &\hfill \w x\inv\in\eta, & \hbox{ for $x$ in $X$,}\hfill\cr
  &\hfill \w y\inv\notin\eta, & \hbox{ for $y$ in $Y$,} \hfill\cr
  &\w z\hfill\notin\eta, & \hbox{ for $z$ in $Z$}\hfill\cr}
  \right\}.
  $$
  Then the collection $\{\W XYZ\}$ forms a fundamental system of open
neighborhoods of $\xi$ in $\specTCK$.

  \proof
  As before, since $\specTCK$ is a topological subspace of $2^\F$ with
the product topology, each $\W XYZ$ is an open set.  In addition, one
has that $\xi\in\W XYZ$ since
  \iBitem $\w\in\xi$, by definition of $\w=\stem\xi$,
  \Bitem $\w x\inv\in\xi$, for $x$ in $X$, because
$X\subseteq\root\w\xi$,
  \Bitem $\w y\inv\notin\xi$, for $y$ in $Y$, since $Y
\cap\root\w\xi=\0$, and
  \Bitem $\w z\notin\xi$, for $z$ in $Z$, because otherwise
$\xi\cap\Pos$ would be bigger than $\set\w$.

  \medskip
  Let $U$ be a neighborhood of $\xi$.  Then, again, $U$ contains a
neighborhood of $\xi$ of the form
  $$
  W = \{\eta\in\specTCK :
  t_1,\ldots,t_p \in \eta,~s_1,\ldots,s_q\notin\eta
  \},
  $$
  and our task is to exhibit a neighborhood of $\xi$, of the form $\W
XYZ$, contained in $W$.

Observe that the collection $\{\W XYZ\}$ is closed under finite
intersections.  Because of this we may restrict our attention to just
two cases, namely
  $$
  W = \{\eta\in\specTCK : t\in\eta\},
  $$
  on the one hand, and
  $$
  W = \{\eta\in\specTCK : s\notin\eta\}
  $$
  on the other.

\medskip\noindent{\tensc case 1:}
  $\xi\in W := \{\eta\in\specTCK : t\in\eta\}$.

  \medskip
  Then $t\in\xi$ and hence, by \lcite{\MembershipCriteria}, $t =
\a\b\inv$ with $|t|=|\a|+|\b|$, and such that $\a$ is a finite sub-word
of $\w$ and $\b$ is either the empty word or a finite admissible word
such that $\last\b\in\root\a\xi$.

Let us suppose first that that $|\a|<|\w|$.  We claim that
  $\W{\0}{\0}{\0}\subseteq W$.
  Observing that $\W{\0}{\0}{\0}$ is just the set of $\eta$'s in
$\specTCK$ such that $\w\in\eta$, we have to show that
  $$
  \w\in\eta \imply t\in\eta.
  $$
  So, given $\eta$ containing $\w$, we have, by convexity, that $\trunc
\w {n-1} \in \eta$, where we denote by $n$ the length of $\w$.  Clearly
both $\w$ and $\trunc \w {n-1}$ belong to $\xi$ as well and so we have,
by \lcite{\Edge}, that
  $$
  \xi\cap \STree{\trunc\w {n-1}}{\w} = \eta\cap \STree{\trunc\w
{n-1}}{\w}.
  $$
  Since $t\in\STree{\trunc\w {n-1}}{\w}$, we conclude that $t$ belongs
to $\eta$, as claimed.

Let us suppose now that $|\a|=|\w|$, which actually means that $\a=\w$.
If, in addition, $\b=e$ then $t = \a\b\inv = \a = \w$ and hence we have
that
  $\W{\0}{\0}{\0}=W$.
  If, on the other hand, $|\b|\geq 1$, let $x=\last\b$ and observe that
by
  \lcite{\MembershipCriteria},
  $x\in\root\w\xi$.  Setting $X=\{x\}$ we claim that
  $\W X{\0}{\0} \subseteq W$, i.e. that, for all $\eta\in\specTCK$ we
have
  $$
  (\w\in\eta) \wedge (\w x\inv\in\eta) \imply t\in\eta.
  $$
  To prove this, let $\eta$ be such that both $\w$ and $\w x\inv$ belong
to $\eta$.  Since $\w$ and $\w x\inv$ also belong to $\xi$ we have, by
\lcite{\Edge}, that
  $$
  \xi\cap \STree{\w x\inv}{\w} = \eta\cap \STree{\w x\inv}{\w}.
  $$
  Since $t\in\STree{\w x\inv}{\w}$, we conclude that $t\in\eta$ as
desired.

\medskip\noindent{\tensc case 2:}
  $\xi\in W := \{\eta\in\specTCK : s\notin\eta\}$.

  \medskip
  Again we must show that some $\W XYZ$ can be found inside $W$.
Observe that, unless the reduced form of $s$ looks like $\a\b\inv$,
where $\a$ and $\b$ are admissible words, we will have that
$W=\specTCK$, which allows for any choice of $\W XYZ$ whatsoever.  So we
will assume that $s=\a\b\inv$ as above.

Let us suppose first that $\w$ is not a sub-word of $\a$. We claim that
  $\W {\0}{\0}{\0}\subseteq W$, i.e. that
  $$
  \w\in\eta \imply s\notin\eta.
  $$
  Let $\eta$ be such that $\w\in\eta$.
  As in case 1, observe that both $\trunc\w {n-1}$ (where
  $n=|\w|$) and $\w$ belong to both $\eta$ and $\xi$.  Therefore, again
by \lcite{\Edge}, we conclude that
  $$
  \xi\cap \STree{\trunc\w {n-1}}{\w} = \eta\cap \STree{\trunc\w
{n-1}}{\w}.
  $$
  Because $\w$ is not a sub-word of $\a$ we have that
  $s\in\STree{\trunc\w {n-1}}{\w}$.  Since, in addition, we have that
$s\notin\xi$, we conclude that $s\notin\eta$, proving the last claim
stated, that is, that $\w\in\eta \imply s\notin\eta$.

Let us suppose now that $\a=\w$ and let $y=\last\b$.  Since, by
assumption, $s\notin\xi$ then at least one of the conditions of
  \lcite{\MembershipCriteria} must fail and the only possibility, at
this stage, is that $y\notin\root\w\xi$.
  Let $Y=\{y\}$, so that $Y\cap\root\w\xi = \0$. We claim that
  $\W {\0}Y{\0} \subseteq W$, i.e. that, for all $\eta$ in $\specTCK$,
  $$
  (\w\in\eta) \wedge (\w y\inv\notin\eta) \imply s\notin\eta,
  $$
  This follows from the fact that $\eta$ is convex.

Finally let us suppose that $\w$ is a proper sub-word of $\a$.  To deal
with this case let $z=\a_{|\w|+1}$ and $Z=\{z\}$.  Then, again by the
convexity of the elements of $\specTCK$, we have that
  $\W {\0}{\0}Z \subseteq W$.
  \proofend

  \section{The unital Cuntz--Krieger algebra}
  \label \UCKAlgSec
  We will consider from now on \cstar-algebras $\OA$ and $\uOA$
generated by partial isometries satisfying a set of conditions that
includes \lcite{(\ELCond)}.  Suppose, for a moment, that the $i^{th}$
row of $A$ is identically zero.  Then, taking $X = \{i\}$ and $Y=\0$
in \lcite{(\ELCond)}, the resulting condition reads ``$S_i=0$''.  In
order to avoid this triviality, we will assume from now on that the
matrix $A$ has no identically zero rows, but we remark that the
algebras that would result from considering matrices having zero rows
are included in our theory, since they are associated to matrices
without trivial rows, which are obtained by eliminating a (possibly
sizable) subset of $\Gen$.

  \definition
  \label \UnitalOADef
  Given an arbitrary set $\Gen$ and a 0--1 matrix
  $A = \{\A ij\}_{i,j\in\Gen}$ with no identically zero rows, we denote
by $\uOA$ the universal \stress{unital} \cstar-algebra generated by a
family of partial isometries
  $\{S_i\}_{i\in\Gen}$
  satisfying conditions \TCKCond1--\TCKCond3 of Section
  {\PCKPisos}, in addition to
  \lcite{(\ELCond)} for each pair of finite subsets $X$ and $Y$ of
$\Gen$ such that $\{j\in\Gen:A(X,Y,j) = 1\}$ is a finite set.

The notation is meant to stress that $\uOA$ is defined to be a universal
object in the category of \stress{unital} \cstar-algebras.  Later we
will study the subalgebra $\OA$ of $\uOA$ generated by the $S_i$'s.  We
will see that $\uOA$ coincides with $\OA$ if $1\in\OA$ and with its
unitization if $1\notin\OA$.  This will make our notation compatible
with the one often used to indicate the adjunction of a unit to an
algebra which does not already have one.

  \state Proposition
  \label \OARelations
  $\uOA$ is naturally isomorphic to
  $\UAlg{\F}{\RelCK}$, where $\RelCK$ is the following set of relations:
  \izitem
  $\e(ts) \e(t) = \e(ts)$,
  whenever $t,s\in\F$ are such that
  $|ts|=|t|+|s|$,
  \zitem $\e(x)\e(y) = 0$, for all $x\neq y$ in $\Gen$,
  \zitem $\e(x\inv) \e(y) = \A xy \e(y)$, for all $x$ and $y$ in $\Gen$,
  \zitem \vrule height0truept depth13truept width0truept
  $
  \prod_{x\in X} \e(x\inv) \prod_{y\in Y} (1 - \e(y\inv)) =
  \sum_{j\in\Gen} A(X,Y,j) \e(j),
  $
  whenever $X$ and $Y$ are finite subsets of $\Gen$ such that $A(X,Y,j)$
vanishes for all but a finite number of $j$'s.

  \proof
  Any non-degenerate representation of $\uOA$ will produce a family of
partial isometries satisfying \TCKCond1--\TCKCond3 of Section
{\PCKPisos} and hence, by \lcite{\Suffering}, also a semi-saturated
partial representation $\pr$ of $\F$.  This \prep will clearly satisfy
(ii) and (iii), as a consequence of \TCKCond2 and \TCKCond3,
respectively, as well as (i) by semi-saturatedness
(cf.~\lcite{\ParcRepAlg}).
  The fourth condition is also clear from \lcite{(\ELCond)}.

Conversely, any partial representation $\pr$ of $\F$ satisfying (i) --
(iv) will give rise to \pisos{} $S_i=\pr(i)$, for $i\in\Gen$, obeying
the requirements of \lcite{\UnitalOADef}.
  \proofend

We have thus realized $\uOA$ as the universal \cstar-algebra for partial
representations of $\F$ satisfying $\RelCK$.  We may then use
\scite{\ELQ}{4.4} to conclude that $\uOA$ is
naturally isomorphic to the covariance \cstar-algebra for the partial
dynamical system given by the partial action of $\F$ on the spectrum
$\GenSpec{\RelCK}$ of the relations $\RelCK$.

Shortening the notation from $\GenSpec{\RelCK}$ to $\specUCK$, recall
that
  $$
  \specUCK =
  \{\xi\in\OmegaE : f(t\inv\xi) = 0, \hbox{ for all } t\in\xi \hbox{ and
} f\in\RelCK'\},
  $$
  where $\RelCK'$ is the subset of $C(\OmegaE)$ consisting of the
functions obtained by replacing each occurrence of $\e(t)$ in $\RelCK$
by the corresponding $1_t$.  See \scite{\ELQ}{Section 4} for more
details.

Later we will study a subset $\specCK$ of $\specUCK$ which will
correspond, in a certain sense, to the subalgebra $\OA$ of $\uOA$, to be
introduced shortly.  We will see that $\specUCK=\specCK$ when $\specCK$
is compact and that $\specUCK$ coincides with the one-point
compactification of $\specCK$ otherwise.
  Again this will make our notation compatible with the one used to
indicate the addition of a point at infinity to a locally compact space
which is not already compact.

Since the relations defining $\uOA$ include those used in the definition
of $\LA$, we see that $\specUCK$ is a subset of $\specTCK$.  Our next
major goal is to to characterize $\specUCK$ by showing it to be a sort
of boundary of $\specTCK$, in a sense similar to that introduced in
  \cite{\Laca}
  for Toeplitz algebras of quasi-lattice ordered groups.
  This will be accomplished once we prove the following:

  \state Proposition
  \label \DeterminationSpecOA
  $\specUCK$ is the closure, within $\specTCK$, of the set of unbounded
elements.

In order to prove this result we must focus on an auxiliary topological
space, namely $2^\Gen$.  Observe that, according to the definition of
roots (Definition \lcite{\RootDef}), if $\xi\in\specTCK$ and $t\in\xi$
then $\root t\xi$ is a subset of $\Gen$.  So we may view $\root t\xi$ as
an element of $2^\Gen$.  Another collection of elements of this
topological space which is relevant to our study is given by the columns
of the matrix $A$, as explained below:

  \definition
  \label \ColumnDef
  For each $j$ in $\Gen$ we will denote by $\col j$ the subset of $\Gen$
given by
  $
  \col j = \{i\in \Gen : \A ij=1\}.
  $

It is usual to interpret the power set of $\Gen$, either as the set of
all subsets of $\Gen$, or as the Cartesian product $\{0,1\}^\Gen$.  If
the latter interpretation is adopted then a column of $A$, which, after
all, is a vector of zeros and ones, is actually an element of
$\{0,1\}^\Gen$.  Of course, these points of view are compatible with one
another and, in particular, the $j^{th}$ column vector corresponds to
$\col j$.

A crucial observation to be made at this point is that if
$\xi\in\specTCK$ and both $t$ and $tz$ belong to $\xi$, where $z\in\Gen$
(that is, $\xi$ contains the edge $(t,tz)$), then the root of $t$
relative to $\xi$ coincides with $\col z$.  This follows from Remark
\lcite{\RootIsColumn}.

Before we embark on the proof of
  \lcite{\DeterminationSpecOA} we need two preliminary results.

  \state Lemma
  \label \BoundedAccumulation
  Let $\xi$ be a bounded element of\/ $\specTCK$ and let $\w$ be the
stem of $\xi$.  If $\xi\in\specUCK$ then $\root\w\xi$ is an accumulation
point of
  $\{\col j\}_{j\in\Gen}$, within $2^\Gen$.

  \proof
  Suppose, by contradiction, that $\root\w\xi$ is not such an
accumulation point.  Then, for some neighborhood $U$ of $\root\w\xi$,
there is only a finite number of $j$'s for which $\col j \in U$.  By the
definition of the product topology on $2^\Gen$, we see that $U$ must
contain a \stress{basic} neighborhood of $\root\w\xi$ of the form
  $$
  V = \{c \in 2^\Gen : x_1,\ldots,x_p\in c,~y_1,\ldots,y_q\notin c \},
  $$
  where $x_1,\ldots,x_p$ and $y_1,\ldots,y_q$ are elements of $\Gen$.
Let
  $$
  X = \{x_1,\ldots,x_p\} \and Y = \{y_1,\ldots,y_q\},
  $$
  and consider the equation
  $$
  \prod_{x\in X} \A xj \prod_{y\in Y} (1-\A yj) = 1
  $$
  in the unknown $j$.  The solutions consist, of course, of those $j$'s
such that for all $x\in X$ and $y\in Y$, one has that
  $\A xj = 1$ and $\A yj = 0$.  Thus $j$ is a solution if and only if
  $x\in\col j$ and $y\notin\col j$
  for all $x\in X$ and $y\in Y$ which, in turn, is equivalent to the
assertion that $\col j\in V$.  Therefore our equation has as many
solutions as there are $\col j$'s in $V$ and hence
  $A(X,Y,j) = 0$ (see \lcite{\AXYj})
  for all but a finite number of $j$'s.
  So
  \lcite{\OARelations.(iv)} applies and, since $\xi\in\specUCK$, we have
that
  $$
  \prod_{x\in X} 1_{x\inv} \prod_{y\in Y} (1 - 1_{y\inv})
\calcat{t\inv\xi} =
  \sum_{j\in\Gen} A(X,Y,j) 1_j \calcat {t\inv\xi} \for{t\in\xi},
  $$
  which can be rewritten as
  $$
  \prod_{x\in X} [tx\INVIN\xi] \prod_{y\in Y} [ty\INVNOTIN\xi] =
  \sum_{j\in\Gen} A(X,Y,j) [tj\IN\xi] \for{t\in\xi}.
  \eqno{(\dagger)}
  $$

  Since $V$ is a neighborhood of $\root\w\xi$, we have that $\root\w\xi
\in V$ which translates into saying that
  $x\in\root\w\xi$ and $y\notin\root\w\xi$, for
  $x\in X$ and $y\in Y$.
  By definition of roots (see \lcite{\RootDef}), we conclude that
  $\w x\inv\in\xi$ and $\w y\inv\notin\xi$ which implies that the left
hand side of \lcite{($\dagger$)} equals 1 for $t=\w$.  Hence precisely
one of the summands on the right hand side equals 1.  In particular, for
some $j$ in $\Gen$, we must have $\w j\in\xi$.  This implies that $\w$
can be prolonged within $\xi$ which cannot happen since $\w$ is the stem
of $\xi$.  We have thus arrived at a contradiction, concluding the
proof.
  \proofend

The second result to be proved before we attack the proof of
\lcite{\DeterminationSpecOA} is in order:

  \state Lemma
  \label \AcumIsLim
  Let $\xi$ be a bounded element of\/ $\specTCK$ and let $\w$ be the
stem of $\xi$. Suppose that $\root\w\xi$ is an accumulation point of
  $\{\col j\}_{j\in\Gen}$, within $2^\Gen$.  Then $\xi$ is in the
closure of the set of unbounded elements of $\specTCK$.

  \proof
  Pick a generic neighborhood $V$ of $\xi$, which, in view of
\lcite{\NBDforBounded}, may be chosen to be of the form
  $\W XYZ$, where $X$, $Y$ and $Z$ are finite subsets of $\Gen$ such
that
  $X\subseteq\root\w\xi$ and $Y\cap\root\w\xi=\0$.  If $|\w|>0$ we may
suppose, without loss of generality, that $\last\w$ belongs to $X$.
  Now, consider the neighborhood $U$ of $\root\w\xi$, within $2^\Gen$,
given by
  $$
  U = \{c \in 2^\Gen :
  x\in c,\ y\notin c,\ \ \forall x\in X,\ \forall y\in Y\}.
  $$

Since, by hypothesis,
  $\root\w\xi$ is an accumulation point of
  $\{\col j\}_{j\in\Gen}$, there are infinitely many $j\in\Gen$ such
that $\col j\in U$.  We choose, among these, a $z_0$ lying outside $Z$.
Observe that the fact that $\col{z_0}$ belongs to $U$ implies that
  $\A x{z_0}=1$ and $\A y{z_0}=0$, for $x$ in $X$ and $y$ in $Y$.

Suppose, for now, that $|\w|>0$.  Then, as mentioned earlier,
  $\last\w\in X$ and hence $\A{\last\w}{z_0}=1$.  So, adding $z_0$ to
the end of $\w$ provides the admissible word $\w z_0$.  Since no row of
$A$ is identically zero, it is possible to extend $\w z_0$ to an
infinite admissible word.  Let $\w'$ be such a word and, using
  \lcite{\ExistenceInfiniteWord},
  let $\xi'$ be the unique element of $\specTCK$ such that
  $\stem{\xi'}=\w'$.  Obviously
  $\xi'$ is an unbounded element of $\specTCK$.

We now aim at proving that $\xi'$ belongs to $\W XYZ$.
  It is easy to see, based on \lcite{\TheWord}, that both $\w$ and
  $\w z_0$ belong to $\xi'$.  So, by \lcite{\CondOne}, there is no other
$z$ such that $\w z \in \xi'$.  In particular, for $z$ in $Z$ we have
that
  $\w z\notin\xi'$.
  Also, by \lcite{\CondTwo}, we have that
  $\w x\inv\in\xi'$, for all $x\in X$, since $\A x{z_0}=1$.  Similarly,
  $\w y\inv\notin\xi'$ for all $y$ in $Y$.
  So, in fact,
  $\xi'$ belongs to $\W XYZ$.
  Summarizing, we have shown that any neighborhood of $\xi$ contains an
unbounded element.  So $\xi$ is in the closure of the set of unbounded
elements of $\specTCK$, as desired.

In the not yet resolved case in which $|\w|=0$, take any infinite word
$\w'$ beginning in
  $z_0$ and choose, as above, $\xi'$ such that
  $\stem{\xi'}=\w'$.  A similar argument shows that $\xi'\in \W XYZ$ and
hence that $\xi$ is in the closure of the set of unbounded elements.
The proof is thus complete.
  \proofend

We are now ready to prove \lcite{\DeterminationSpecOA}.

  \proof
  (of \DeterminationSpecOA) Let us begin by showing that the closure of
the set of unbounded elements of $\specTCK$ is contained in $\specUCK$.
Since $\specUCK$ is closed, it suffices to show that it contains every
unbounded element.

So let $\xi\in\specTCK$ be unbounded and let $\w$ be the stem of $\xi$.
We claim that, for every
  $t$ in $\xi$, there exists $z$ in $\Gen$ such that $tz\in\xi$.  To see
this, write $t=\a\b\inv$ as explained in \lcite{\MembershipCriteria}.
In case $|\b|>0$ we clearly may take $z$ to be $\b_1$.
  In case $|\b| = 0$, recall that $\a$ is a finite sub-word of the
infinite word $\w$.  So, if we take $z$ to be
  $\w_{|\a|+1}$, we'll have $tz = \trunc{\w}{|\a|+1} \in \xi$, proving
our claim.

Now, showing that $\xi\in\specUCK$ amounts to proving that
  $$
  \prod_{x\in X} [tx\INVIN\xi] \prod_{y\in Y} [ty\INVNOTIN\xi] =
  \sum_{j\in\Gen} A(X,Y,j) [tj\IN\xi] \for{t\in\xi},
  \eqno{(\dagger)}
  $$
  whenever $X$ and $Y$ are finite subsets
 of $\Gen$ such that $A(X,Y,j)=0$ for all but a finite number of $j$'s.
  So, let $X$ and $Y$ be chosen accordingly and let $t\in\xi$.  We
denote by $z$ the (necessarily unique by \lcite{\CondOne}) element of
$\Gen$ such that $tz\in\xi$.
  The right hand side of
  \lcite{($\dagger$)} then takes the value $A(X,Y,z)$.

Observe that, by \lcite{\CondTwo}, we have that, for all $x\in\Gen$,
  $\A xz = 1$ if and only if $tx\inv\in\xi$.  Using the boolean valued
operator $[\cdot]$, we may rewrite this as
  $\A xz = [tx\INVIN \xi]$.  We then have, cf. \lcite{(\AXYj)}, that
  $$
  A(X,Y,z) =
  \prod_{x\in X} \A xz \prod_{y\in Y} (1-\A yz) =
  \prod_{x\in X} [tx\INVIN\xi] \prod_{y\in Y} [ty\INVNOTIN\xi],
  $$
  which is precisely the left hand side of \lcite{($\dagger$)}.  This
proves that $\xi\in\specUCK$.

We must now prove that every bounded element of $\specUCK$ is in the
closure of the unbounded elements of $\specTCK$.  But this follows
immediately from applying, successively, \lcite{\BoundedAccumulation}
and \lcite{\AcumIsLim}.
  \proofend

One of the side consequences of this is the following simple
characterization of the bounded elements of $\specUCK$:

  \state Corollary
  \label \AcumIfANdOnlyIfIsLim
  Let $\xi\in\specTCK$ be bounded and let $\w$ be the stem of $\xi$.
Then $\xi$ is in $\specUCK$ if and only if $\root\w\xi$ is an
accumulation point of $\{\col j\}_{j\in\Gen}$, within $2^\Gen$.

  \proof
  The \stress{only if} part is nothing but
\lcite{\BoundedAccumulation}. So let us deal with the converse and hence
we assume that $\root\w\xi$ is an accumulation point of
  $\{\col j\}_{j\in\Gen}$.  By \lcite{\AcumIsLim}, we see that $\xi$
is in closure of the set of unbounded elements of $\specTCK$ and hence
in $\specUCK$, by \lcite{\DeterminationSpecOA}.
  \proofend

The following is also a useful property of the elements of $\specUCK$:

  \state Proposition
  \label \MustContinue
  Let $\xi\in\specUCK$ and let $\a$ be a finite admissible word of
positive length such that $\a\in\xi$.  Suppose that the row of $A$
labeled by $\last\a$ is finite (in the sense that it contains only a
finite number of nonzero entries).  Then there exists $j$ in $\Gen$ such
that
  $\a j\in\xi$.

  \proof
  Let $x=\last\a$ and $X=\{x\}$. Observe that, since
  $A(X,\0,j) = \A xj$ then $A(X,\0,j)$ vanishes for all but a finite
number of $j$'s.  Hence
  \lcite{\OARelations.(iv)} applies and we have that
  $$
  [tx\INVIN\xi] =
  \sum_{j\in\Gen} A(X,\0,j) [tj\IN\xi] \for{t\in\xi}.
  $$
  With $t=\a$ we have that $tx\inv=\a(\last\a)\inv$ which is a sub-word
of $\a$ and hence represents an element belonging to $\xi$.  Hence the
right hand side above also equals 1 and so there is some $j\in\Gen$ such
that $\a j\in\xi$.
  \proofend

In the special case of a row-finite $A$ (see also \cite{\KPRR}), the
spectrum $\specUCK$ has a much simpler description.

  \state Proposition
  \label \specRowFinite
  Suppose that each row of $A$ has a finite number of nonzero entries.
Let $\Omega_\infty$ be the set of unbounded elements of\/ $\specTCK$ and
let $\phi=\{e\}$ be the smallest element of $\specTCK$.  Then
  $$
  \Omega_\infty \subseteq
  \specUCK \subseteq
  \Omega_\infty \cup \{\phi\}.
  $$
  In addition, if $\Gen$ is an infinite set then $\specUCK$ coincides
with the set on the right hand side. If, on the other hand, $\Gen$ is
finite then $\specUCK=\Omega_\infty$.

  \proof
  We have already seen that the unbounded elements of $\specTCK$ belong
to $\specUCK$.
  To prove the second inclusion in the statement we must show that if
$\xi\in\specUCK$ is bounded then $\xi=\phi$.  Given such a $\xi$ let
$\w$ be its stem and observe that, by
  \lcite{\MustContinue}, $\w$ is either infinite or the empty word.
Since the first possibility is ruled out, we have that $\w=e$.

On the other hand, by \lcite{\AcumIfANdOnlyIfIsLim}, we have that $\root
e\xi$ is an accumulation point of the columns of $A$.  It is easy to see
that, since $A$ is row-finite, the only possible accumulation point of
its columns is the zero vector.  It follows that $\root e\xi$ is the
empty set and hence that $\xi=\phi$ by \lcite{\UniqueWordRoot}.

Assume that $\Gen$ is infinite and let us prove that the second
inclusion in the statement is actually an equality.  For this we only
need to prove that $\phi\in\specUCK$.  Since $2^\Gen$ is a compact space
and there are infinitely many columns in $A$, these columns must have an
accumulation point which, by the reasoning above, must be the zero
vector.  That $\phi\in\specUCK$ then follows from
\lcite{\AcumIfANdOnlyIfIsLim}.

To conclude we must only show that,
  in case $\Gen$ is finite,
  $\phi$ is not in $\specUCK$,
  but this also follows from \lcite{\AcumIfANdOnlyIfIsLim}, since a
finite set of columns does not admit an accumulation point.
  \proofend

Since $\uOA$ is isomorphic to $\UAlg{\F}{\RelCK}$, by
\lcite{\OARelations}, we may conclude from
\scite{\ELQ}{4.4} that $\uOA$ is also the crossed
product of $C(\specUCK)$ by $\F$.  The above characterization of
$\specUCK$ gives this picture a much more concrete flavor.  Together
with that characterization, the following is one of our main results:

  \state Theorem
  \label \UnitalOAasCrossProd
  Let $A = \{\A ij\}_{i,j\in \Gen}$ be a 0--1 matrix with no
identically zero rows and let $\specUCK$ be the closure, within
$\specTCK$, of the set of the unbounded elements.  For each $t\in\F$ let
  $$
  \Xa_t = \{\xi\in\specUCK : t\in \xi\},
  $$
  and let
  $$
  \ha_t : \Xa_{t\inv} \arw \Xa_t
  $$
  be given by $\h_t(\xi)=t\xi$, for $\xi$ in $\Xa_{t\inv}$.
  Then the pair $\(\{\Xa_t\}_{t\in\F}, \{\ha_t\}_{t\in\F}\)$ is a
partial action of $\F$ on $\specUCK$ and
  $\uOA$ is isomorphic to the crossed product
  $C(\specUCK)\crossproduct \F$
  under an isomorphism that maps each generating \piso{} $S_x$ to
$1_x\delta_x$, where $1_x$ is the characteristic function of $\Xa_x$.

  \proof
  Follows from \lcite{\OARelations}, \lcite{\DeterminationSpecOA}, and
\scite{\ELQ}{4.4}.
  \proofend

Recall from Section \lcite{\Preliminaries} that the crossed product
algebra associated to a given a partial dynamical system is obtained by
taking the cross-sectional \cstar-algebra of the corresponding
semi-direct product bundle.

In the case of the partial action of a given group $G$ on the spectrum
$\GenSpec{\Rel}$, for a set $\Rel$ of relations, it is easy to see that
the semi-direct product bundle is isomorphic to the bundle obtained (as
in \scite{\Amena}{6.1}) from the partial representation $\pr$ of $G$
arising from any faithful representation of $\UAlg G\Rel$.  The crucial
observation supporting this statement is that the collection of all the
  $$
  \e(t) := \pr(t)\pr(t)^*
  $$
  generates the unit fiber of the bundle.

Since $\RelCK$ includes \lcite{\OARelations.(i--ii)} we have that any
partial representation of $\F$ satisfying $\RelCK$ is semi-saturated and
orthogonal.  From Theorem \lcite{4.1} of \cite{\Ortho} we then deduce
the following:

  \state Proposition
  \label\OurBundleAP
  The semi-direct product bundle corresponding to the partial action of
$\F$ on $\specUCK$ satisfies the approximation property of
\cite{\Amena}.  Therefore this bundle is amenable and hence
  the (full) crossed product algebra
  $C(\specUCK)\crossproduct \F$
  coincides with its reduced version
  $C(\specUCK)\crossproduct_r \F$.

  \state Proposition
  \label \Nuke
  Let $A = \{\A ij\}_{i,j\in\Gen}$ be a 0--1 matrix with no
identically zero rows.  Then $\uOA$ is nuclear.

  \proof
  As seen above $\uOA$ is the cross sectional \cstar-algebra of a Fell
bundle satisfying the approximation property.  The unit fiber of this
bundle is $C(\specUCK)$, which is a nuclear \cstar-algebra.  The result
then follows from \cite{\FAbadie}.
  \proofend

  \section{Cuntz--Krieger algebras for infinite matrices}
  \label \JASection
  Given an arbitrary 0--1 matrix
  $A = \{\A ij\}_{i,j\in\Gen}$,
  we will now define the \cstar-algebra $\OA$,
  generalizing the definition given in \cite{\CKbib} for finite matrices
and in \cite{\KPRR} for row-finite matrices.  As usual, we assume that
$A$ has no row consisting only of zeros.

  \definition
  \label \NonUnitalOADef
  We denote by $\OA$ the sub-\cstar-algebra of $\uOA$ generated by
  $\{S_i: i\in\Gen\}$.

Recall that $\uOA$ was defined via a universal property in the category
of \stress{unital} \cstar-algebras and so $\uOA$ is obviously unital.
On the other hand notice that $\OA$ need not include the unit of $\uOA$
as there is no obvious reason for $\OA$ to coincide with $\uOA$.
However, the addition of the unit to $\OA$ is clearly enough to generate
the whole of $\uOA$.  From this it is easy to see that:

  \state Proposition
  \label \OAIsIdeal
  $\OA$ is a two sided ideal in $\uOA$ of linear co-dimension at most
one.

As a consequence of \lcite{\Nuke} we then have:

  \state Corollary
  For every 0--1 matrix $A$ with no identically zero rows,
  $\OA$ is a nuclear \cstar-algebra.

  \proof
  $\OA$ is an ideal in the nuclear \cstar-algebra $\uOA$, and hence it
is nuclear by \cite{\ChoiEffros}.
  \proofend

The following result characterizes $\OA$ as the crossed product algebra
for a partial dynamical system.

  \state Theorem
  \label \NonUnitOAasCrossProd
  Let $A = \{\A ij\}_{i,j\in \Gen}$ be a 0--1 matrix with no
identically zero rows and view $\uOA$ as the crossed product
  $C(\specUCK)\crossproduct \F$, as in \lcite{\UnitalOAasCrossProd}.
Also let $\phi=\{e\}$ be the smallest element of $\specTCK$ and put
  $\specCK:=\specUCK\setminus\{\phi\}$.  Then
  \izitem $\specCK$ is an open invariant subset of $\specUCK$,
  \zitem $\OA \cap C(\specUCK) = C_0(\specCK)$, and
  \zitem $\OA \simeq C_0(\specCK)\crossproduct \F$.

  \proof
  Initially note that $\phi$ may or may not be in $\specUCK$, as seen in
\lcite{\specRowFinite}.  Therefore $\specCK$ may or may not coincide
with $\specUCK$.  If $\phi\in\specUCK$ then $\{\phi\}$ is a closed
invariant subset, hence (i).

By \lcite{\OAIsIdeal}, $\OA$ is also the closed ideal of $\uOA$
generated by
  $\{S_x S_x^*: x\in\Gen\}$.  Since $S_x S_x^* = 1_x$ we see that this
set is contained in $C(\specUCK)$.  Denoting by $\big\<C\big\>$ the
closed ideal of $\uOA$ generated by a subset
  $C\subseteq \uOA$, we have
  $$
  \OA =
  \big\<\{S_x S_x^*: x\in\Gen\}\big\> \subseteq
  \big\<\OA \cap C(\specUCK)\big\> \subseteq
  \OA,
  $$
  so equality holds throughout.  Let $I = \OA \cap C(\specUCK)$.  Then,
clearly, $I$ is an invariant ideal of $C(\specUCK)$ and, as seen above,
$\OA=\big\<I\big\>$.  Using \scite{\ELQ}{3.1} we have that
$\OA = I\crossproduct \F$.  To conclude the proof it is therefore enough
to prove (ii) or, in other words, that $I=C_0(\specCK)$.

Since $I$ is an invariant ideal of $C(\specUCK)$ we have that
  $I = C_0(U)$ for some open invariant set $U\subseteq\specUCK$.  For
each $x\in\Gen$ we have
  $$
  1_x = S_x S_x^* \in C_0(U)
  \and
  1_{x\inv} = S_x^* S_x \in C_0(U).
  $$
  It follows that both $\Xa_x$ and $\Xa_{x\inv}$ are contained in $U$.
We leave it for the reader to verify that
  $$
  \bigcup_{x\in\Gen} (\Xa_x \cup \Xa_{x\inv}) =
\specUCK\setminus\{\phi\},
  $$
  which implies that $\specCK\subseteq U$.  Hence
  $C_0(\specCK)\subseteq C_0(U) = I$.

  To prove the reverse inclusion observe that, for every $x\in\Gen$, we
have that $1_x$ vanishes at $\phi$.  Therefore
  $1_x\in C_0(\specCK)$ and thus
  $S_x \in \big\< C_0(\specCK) \big\>$.  It follows that $\OA\subseteq
\big\<C_0(\specCK)\big\>$ and hence that
  $$
  I =
  \OA\cap C(\specUCK) \subseteq \big\<C_0(\specCK)\big\>\cap C(\specUCK)
= C_0(\specCK),
  $$
  where the last equality follows from
  \scite{\ELQ}{3.1}.
  \proofend

The following is a useful criteria to determine when $\OA$ is unital.

  \state Proposition
  \label \AboutUnits
  Let $A = \{\A ij\}_{i,j\in\Gen}$ be a
  0--1 matrix with no identically zero rows.  Then the following
are equivalent
  \izitem $\OA$ is a unital \cstar-algebra,
  \zitem $\OA = \uOA$,
  \zitem the element $\phi=\{e\}$ does not belong to $\specUCK$,
  \zitem the identically zero vector in $2^\Gen$ is not an accumulation
point for the columns of $A$,
  \zitem there exists a finite set $Y\subseteq\Gen$ such that
$A(\0,Y,j)$ is finitely supported in $j$.

  \proof
  \def\part#1#2#3{(#1$#2$#3):}
  \medskip\part {iii}{\bimply}{iv} Observe that $\stem\phi=e$ and that
$\root e\phi$ is the empty set, which should be interpreted as the
  zero vector in $2^\Gen$.  So, employing \lcite{\AcumIfANdOnlyIfIsLim},
we have that $\phi$ is in $\specUCK$ if and only if the zero vector is
an accumulation point of the set of columns of $A$.

  \medskip\part{iii}{\imply}{ii}
  If $\phi\notin\specUCK$ then, again by \lcite{\NonUnitOAasCrossProd},
we have that
  $1\in C(\specUCK) = C_0(\specCK) \subseteq \OA$.  So $\OA$ contains
the unit of $\uOA$ and hence $\OA=\uOA$.

  \medskip\part{ii}{\imply}{i} Obvious.

  \medskip\part {i}{\imply}{iii}
  Suppose that $f$ is a unit for $\OA$.  We will first show that
  $f\in C_0(\specCK)$.
  Since $1_x\in\OA$ for all $x\in\Gen$, we have that $f1_x=1_x$.
  If we apply the conditional expectation $E$ of \scite{\Amena}{2.9} to
this last identity we conclude that
  $E(f) 1_x = 1_x$ and hence that
  $$
  E(f) S_x = E(f) S_x S_x^* S_x = E(f) 1_x S_x = 1_x S_x = S_x.
  $$
  Similarly, one can show that $E(f) S_x^* = S_x^*$ and hence that
$E(f)$ is a unit for $\OA$.  Thus $E(f)=f$.

It then follows that
  $
  f\in \OA\cap C(\specUCK) = C_0(\specUCK \backslash \{\phi\}),
  $
  and hence that $f$ is the characteristic function of the complement of
$\{\phi\}$ within $\specUCK$.  But then, if $\phi\in\specUCK$, we would
have that $\phi$ is an isolated point of $\specUCK$ which is impossible
in view of the fact that the unbounded elements are dense in $\specUCK$
(see \lcite{\DeterminationSpecOA}).  Thus $\phi\notin\specUCK$.

  \medskip\part{v}{\bimply}{iii}
  A fundamental system of neighborhoods of the zero vector in $2^\Gen$
is obtained by taking neighborhoods of the form
  $$
  U = \{c \in 2^\Gen : y\notin c,\ \forall y\in Y\},
  $$
  where $Y$ ranges in the collection of finite subsets of $\Gen$.  By
\lcite{\AcumIfANdOnlyIfIsLim}, we have that $\phi$ does not belong to
$\specUCK$ if and only if there exists $Y$, as above, such that the set
  $\{j: \col j\in U\}$ is finite.  Since $\col j\in U$ if and only if
  $A(\0,Y,j)=1$, the proof is concluded.
  \proofend

Observe that if the equivalent conditions of \lcite{\AboutUnits} hold
then
  $
  \specCK =
  \specUCK\setminus\{\phi\} =
  \specUCK
  $
  and hence $\specCK$ is a compact space.
  On the other hand, if the conditions of \lcite{\AboutUnits} fail, that
is, $\phi\in\specUCK$, then $\phi$ is not an isolated point of
$\specUCK$ because of the density of the unbounded elements.  In this
case $\specCK$ is not compact and its one-point compactification is
clearly $\specUCK$.  Therefore our notation turns out to be compatible
with the one used to indicate the addition of a point at infinity to a
locally compact space which is not already compact.

We would now like to characterize $\OA$ as the universal (not
necessarily unital) \cstar-algebra generated by the $S_x$'s.

  \state Theorem
  \label \JAUniversal
  Let $\Gen$ be any set,
  $A = \{\A ij\}_{i,j\in \Gen}$ be a 0--1 matrix with no
identically zero rows and $\{T_x\}_{x\in\Gen}$ be a family of partial
isometries on a Hilbert space $H$ satisfying
  \izitem \commute{$T_i^*T_i$}{$T_j^*T_j$},
  \zitem $T_i^*T_j=0$, for $i\neq j$,
  \zitem $T_i^*T_iT_j = \A ij T_j$,
  \zitem
  whenever $X$ and $Y$ are finite subsets of $\Gen$ such that $A(X,Y,j)$
is zero for all but a finite number of $j$'s, one has
  $$
  \prod_{x\in X} T_x^* T_x \prod_{y\in Y} (1 - T_y^* T_y) =
  \sum_{j\in\Gen} A(X,Y,j) T_j T_j^*.
  $$
  \medskip\noindent
  Then there exists a unique representation $\pi$ of $\OA$ on $H$ such
that
  $\pi(S_x)=T_x$.

  \proof
  Recall from \lcite{\UnitalOADef} that $\uOA$ is the universal unital
\cstar-algebra generated by \pisos satisfying the above conditions.
  So, the given set of partial isometries will provide for a
non-degenerate representation $\rho$ of $\uOA$ on $H$ such that
$\rho(S_x)=T_x$.  The hypothesis do not impose an explicit restriction
on the value of $\rho$ on the unit of $\uOA$ and so we cannot say that
$\rho$ is unique.  However, the restriction $\pi$ of $\rho$ to $\OA$ is
uniquely determined and will clearly satisfy the requirements above.
  \proofend

This results expresses that $\OA$ is the universal not necessarily
unital \cstar-algebra generated by the $S_x$'s under the conditions
above; but the reference to ``1'' in the fourth condition needs
clarification.  Observe that if $X$ is not the empty set
then the left hand side of (iv), once expanded, may be rewritten without
reference to the identity.  So, we can get rid of all references to the
identity as long as the matrix $A$ is such that (iv) never occurs for an
empty set $X$.  That is, as long as the condition that $A(X,Y,j)$
vanishes for all but a finite number of $j$'s never holds for $X=\0$.

On the other hand if, for some $Y$, one has that $A(\0,Y,j)=0$ for all
but a finite number of $j$'s then the corresponding occurence of (iv)
cannot be rewritten without explicit mention of ``1''.  In this case
our relations force us back into the category of unital algebras.
Indeed, by \lcite{\AboutUnits}, the existence of such a $Y$ implies
that $\OA$ is unital!

One of the important consequences of the description of $\OA$ as a
crossed product is:

  \state Proposition
  \label \ItExists
  Let $A = \{\A ij\}_{i,j\in \Gen}$ be a 0--1 matrix with no
identically zero rows.  Then:
  \izitem The canonical partial isometries $S_x\in\OA$ are non-zero.
  \zitem More generally, if $\a$ is a finite admissible word then
  $S_{\a_1} \ldots S_{\last\a} \neq 0$.
  \zitem There exists a collection $\{T_x\}_{x\in\Gen}$ of
  \stress{non-zero} partial isometries on some Hilbert space
satisfying conditions \lcite{\JAUniversal.(i--iv)}.

  \proof
  Under the isomorphism $\OA \simeq C_0(\specCK)\crossproduct \F$ one
has that $S_{\a_1} \ldots S_{\last\a}$ corresponds to $1_\a\delta_\a$,
where $1_\a$ is the characteristic function of $\Xa_\a$.  To prove
(ii) one simply has to note that $\Xa_\a$ is nonempty since it
contains any $\xi\in\specCK$ whose stem is an infinite admissible word
extending $\a$.  This proves (ii) and hence also (i).  To prove (iii)
it is enough to consider the images of the canonical $S_x$'s in any
faithful representation of $\OA$.
  \proofend

Even though we have avoided using the strong operator topology in our
study, the following fact should be noted:

  \state Proposition
  \label \SOT
  Let $A = \{\A ij\}_{i,j\in \Gen}$ be a 0--1 matrix with no
identically zero rows and suppose that $\{T_i\}_{i\in\Gen}$ is a
family of partial isometries on a Hilbert space $H$ satisfying
conditions \CKcond1 and \CKcond2 (from Section \lcite{\Introduction})
with respect to the strong operator topology.  Then
  \izitem the $T_i$ satisfy conditions
  \lcite{\JAUniversal.(i--iv)}, and
  \zitem there exists a representation $\pi$ of\/ $\OA$ on $H$ sending
each generating partial isometry
   $S_x\in\OA$ to the corresponding $T_x$.

  \proof
  The first statement follows easily from \CKcond1 and \CKcond2, and
the second is a consequence of \lcite{\JAUniversal}.
  \proofend

Let us now present some examples in order to relate our study to
various versions of Cuntz--Krieger algebras that have appeared in the
literature.  The proofs of the following statements are elementary and
left to the reader.

  \sysstate{Examples}{\rm}
  {Let $A = \{\A ij\}_{i,j\in \Gen}$ be a 0--1 matrix with no
identically zero rows.
  \iItem{(a)} If $\Gen$ is finite and  $A$ satisfies condition (I)
of \cite{\CKbib} then our $\OA$ coincides with the algebra $\OA$
introduced by Cuntz and Krieger in \cite{\CKbib}.
  \Item{(b)} More generally, if $\Gen$ is finite then our $\OA$
coincides with the algebra ${\cal A}\OA$ of an Huef and Raeburn
\cite{\HR}.
  \Item{(c)} If $\Gen$ is infinite but there are only finitely many
ones in each row of $A$ then our $\OA$ coincides with the algebra $\OA$
studied by Kumjian, Pask, Raeburn, and Renault in \cite{\KPRR}.
  \Item{(d)} If $\Gen$ is countably infinite and $\A ij=1$ for all $i$
and $j$ in $\Gen$ then our $\OA$ coincides with the algebra ${\cal
O}_\infty$ of \cite{\CuntzOInfinite} (cf.~\scite{\CKbib}{Remark
2.15}).}

  \section{A faithful representation}
  The density of the unbounded elements in the spectrum $\specUCK$ of
the generalized Cuntz-Krieger relations $\RelCK$ can be used to
construct a natural faithful representation of $\OA$.  As before, we let
$A = \{\A ij\}_{i,j\in\Gen}$ be a fixed 0--1 matrix having no
identically zero rows.

We will denote by $P_A$ the set of all infinite paths on $\Graph$.  As
usual, the elements of $P_A$ will also be seen as infinite admissible
words.

Considering the Hilbert space
  $ \ell^2(P_A) $
  and its canonical orthonormal basis
  $\{\bvec_\w: \w \in P_A \}$,
  define, for each $x\in\Gen$, the partial isometry
  $L_x \in\B(\ell^2(P_A))$
  by
  $$
  L_x (\bvec_\w) =
  \left\{\matrix{
  \bvec_{x\w} & \text{ if } A(x,\w_1) = 1,\cr
  0 \hfill & \text{ if } A(x,\w_1) = 0.}
  \right. 
  $$
  Observe that the conditon that $A(x,\w_1) = 1$ is equivalent to the
infinite word $x\w$ being admissible.

  Denote by
  $\lambda$
  the left regular representation of
  $\F$
  on
  $\ell^2(\F)$.


  \state Proposition
  Let $A = \{\A ij\}_{i,j\in\Gen}$ be a 0--1 matrix with no
identically zero rows. Then
  \izitem There exists a unique representation
  $$
  \pi : \uOA \arw \B(\ell^2(P_A))
  $$
  such that $\pi(S_x)=L_x$ for all $x\in\Gen$.
  \zitem There exists a unique representation
  $$
  \rho : \uOA \arw \B(\ell^2(P_A)\otimes\ell^2(\F))
  $$
  such that $\rho(S_x)=L_x\otimes\lambda_x$ for all $x\in\Gen$.
  \zitem $\rho$ is faithful.

  \proof
  With respect to (i) it suffices to show that the family
$\{L_x\}_{x\in \Gen}$ satisfies conditions
  \TCKCond1--\TCKCond3 of Section {\PCKPisos},
  as well as condition
  \lcite{(\ELCond)}
  for each pair of finite subsets $X$ and $Y$ of $\Gen$ such that
$A(X,Y,j)$ is finitely supported in $j$.

  Let $Q_x$ be the projection onto the closed linear span of
  $\{\bvec_\w: A(x,\w_1) = 1\}$
  and let $P_x $ be the projection onto the closed linear span of
  $\{\bvec_\w: \w_1 = x\}$.
  It is obvious that $L_x$ has initial projection $Q_x$ and final
projection $P_x$.

Since $Q_x$ and $P_x$ are defined by subsets of the canonical basis,
their products are determined by the intersection of the corresponding
subsets, so \TCKCond1 holds, and since the subsets corresponding to
$P_i$ and $P_j$ are disjoint for $i \neq j$, \TCKCond2 holds as well. To
prove \TCKCond3 simply observe that the intersection of the set
  $\{\bvec_\w: A(i,\w_1) =1\}$,
  corresponding to $Q_i$,
  and the set
  $\{\bvec_{\w} : \w_1 = j\}$,
  corresponding to $P_j$,
  equals
  $\{\bvec_{\w} : \w_1 = j \}$ when $A(i,j) = 1$ and reduces to the
empty set otherwise.

  We must now verify condition \lcite{(\ELCond)}, for which we assume
that $X$ and $Y$ are two finite subsets of $\Gen$ such that $A(X,Y,j)$
is finitely supported on the variable $j\in \Gen$.  A simple argument
shows that
  $$
  \big(\bigcap_{x\in X}\{\bvec_\w: A(x,\w_1) = 1\}\big)
  \cap
  \big(\bigcap_{y\in Y}\{\bvec_\w: A(y,\w_1) = 0\}\big)
  =
  \bigcup_{A(X,Y,j) =1} \{\bvec_\w : \w_1 = j\},
  $$
  from where condition
  \lcite{(\ELCond)} follows without difficulty.

Since $\uOA$ is the universal unital \cstar-algebra generated by partial
isometries
  $\{S_i\}_{i\in\Gen}$
  satisfying the conditions above, the existence and uniqueness of $\pi$
follows.

To prove (ii) one needs only note that the argument used above also
applies to the family of partial isometries
  $\{L_x\otimes\lambda_x\}_{x\in\Gen}$.

Let $L$ denote the partial representation of $\F$ given by
  $$
  L(t) = \pi(\pr(t)) \for t\in\F,
  $$
  where $\pr : \F \to \uOA = \UAlg{\F}{\RelCK}$ is the universal partial
representation of
  $\F$ satisfying $\RelCK$ (see Section \lcite{\Preliminaries}).
Obviously
  $L(x) = L_x$.

We leave it for the reader to verify that
  if $\a$ is a finite admissible word and $\w\in P_A$ then
  $$
  L(\a) (\bvec_\w) =
  \left\{\matrix{
  \bvec_{\a\w} & \text{ if $\a\w$ is admissible,} \cr
  0 \hfill & \text{ otherwise,}\hfill}
  \right. 
  $$
  and that
  $$
  L(\a)^* (\bvec_\w) =
  \left\{\matrix{
  \bvec_{\a\inv\w} & \text{ if $\a$ is a sub-word of $\w$,} \cr
  0 \hfill & \text{ otherwise.}\hfill}
  \right. 
  $$
  We should point out that $\a\inv\w$ is supposed to mean the infinite
word obtained by deleting the initial segment $\a$ from $\w$.

In order to prove that $\rho$ is faithful we will use the conditional
expectation
  $$
  E: C(\specUCK) \crossproduct \F \to C(\specUCK)
  $$
  associated to the semi-direct product bundle as in
\scite{\Amena}{2.9}.  We will also use the map
  $$
  \Delta: \B(\ell^2(P_A)\otimes \ell^2(\F)) \to
    \B(\ell^2(P_A)\otimes \ell^2(\F))
  $$
  which assigns, to every operator on $\ell^2(P_A)\otimes \ell^2(\F)$,
its (block) diagonal relative to the canonical basis of $\ell^2(\F)$.
In symbols
  $$
  \Delta(T) = \sum_{t\in\F} (1\otimes p_t) T (1\otimes p_t)
  \for T \in \B(\ell^2(P_A)\otimes \ell^2(\F)),
  $$
  where the $p_t$ are the rank one projections associated to the standard
basis of $\ell^2(\F)$.

  Identifying $\uOA$ and $C(\specUCK)\crossproduct \F$ as permited by
  \lcite{\UnitalOAasCrossProd},
  it is easy to see that the following diagram commutes:
  $$
  \matrix{
  C(\specUCK)\crossproduct \F &
    \labar{\rho} &
    \B(\ell^2(P_A)\otimes \ell^2(\F))\cr\cr
  E \Big\downarrow &&
  \Big\downarrow \Delta \cr\cr
  C(\specUCK) &
    \labar{\rho} &
    \B(\ell^2(P_A)\otimes \ell^2(\F))
  }
  $$

We initially show that $\rho$ is faithful on $C(\specUCK)$.   But, since
  $\rho = \pi\otimes 1$ on $C(\specUCK)$, it is enough to show
faithfulness of $\pi$, instead.
In order to
do this let $f\in C(\specUCK)$.  We claim that if the unbounded element
$\xi\in\specUCK$ and $\w\in P_A$ are related by $\w=\stem\xi$ (see
\lcite{\TheWord}) then
  $$
  \pi(f)( \bvec_{\w} )
  =
  f(\xi) \; \bvec_{\w}.
  \eqno{(\dagger)}
  $$
  Since $C(\specUCK)$ is generated as a \cstar-algebra by the
projections $1_t$ with $t\in \F$, it suffices to prove
\lcite{($\dagger$)} for $f=1_t$.  In addition we may assume that the
reduced form of $t$ is $\a \b\inv$ for admissible words $\a$ and $\b$
since, otherwise, $1_t = 0$ in $C(\specUCK)$ by
  \lcite{\MembershipCriteria}.  Observe that
  $$
  \pi(1_t) =
  \pi\big(\pr(t)\pr(t)^*\big) =
  \pi\big(\pr(\a)\pr(\b)^*\pr(\b)\pr(\a)^*\big) =
  L(\a)L(\b)^*L(\b)L(\a)^*.
  $$

Suppose first that $t\in\xi$.  Then, by \lcite{\MembershipCriteria}, we
have that $\a$ is a sub-word of $\w$.  Hence, as seen above,
  $L(\a)^*(\bvec_\w) = \bvec_{\a\inv\w}$.  Again by
\lcite{\MembershipCriteria}, it follows that $\b\a\inv\w$ is an
admissible word and hence
  $L(\b)L(\a)^*(\bvec_\w) = \bvec_{\b\a\inv\w}$.
  We then have that $\bvec_\w$ belongs to the initial space of
  $L(\b)L(\a)^*$ and so
  $$
  \pi(1_t)(\bvec_\w) =
  L(\a) L(\b)^* L(\b) L(\a)^* (\bvec_\w) =
  \bvec_\w.
  $$
  On the other hand we have that $1_t(\xi) = [t\IN\xi] = 1$, verifying
\lcite{($\dagger$)} under the assumption that $t\in\xi$.

  Suppose now that $t\notin\xi$.  Then at least one of the conditions in
\lcite{\MembershipCriteria} must fail and then one may verify that
  $L(\b)L(\a)^*(\bvec_\w) = 0$. It follows that
  $$
  \pi(1_t)(\bvec_\w) =
  L(\a) L(\b)^* L(\b) L(\a)^* (\bvec_\w) =
  0 =
  1_t(\xi)\bvec_\w,
  $$
  concluding the proof of \lcite{($\dagger$)}.

  This immediately implies our claim that $\pi$, and hence also
$\rho$, is faithful on $C(\specUCK)$: if $\pi(f)=0$ then
\lcite{($\dagger$)} says that $f$ vanishes on every unbounded
$\xi\in\specUCK$.  Since these form a dense subset of $\specUCK$ by
\lcite{\DeterminationSpecOA}, one must have $f=0$.

The proof that $\rho$ is faithful on all of $\uOA$ is now routine
(e.g.~\scite{\newpim}{2.9} and \cite{\LR}): if
  $a\in\uOA$ is such that $\rho(a)=0$ then, by the commutativity of the
diagram above, we have
  $$
  0 = \Delta(\rho(a^*a)) =
  \rho(E(a^*a)),
  $$
  and hence $E(a^*a)=0$.  Since $E$ is a faithful conditional
expectation by
  \lcite{\OurBundleAP} and \scite{\Amena}{3.6}, we have that $a=0$.
  \proofend

  \section{Invariant sets}
  \label \InvarSection
  This section is dedicated to the study of the first important
dynamical aspect of the partial action of $\F$ on $\specUCK$, namely
open invariant subsets.  Afterwards we will also look into fixed points
and topological freeness.

Recall that a subset $C\subseteq\specUCK$ is said to be invariant if
$\ha_t( C \cap \Xa_{t\inv} ) \subseteq C$, for all $t$ in $\F$.  Even if
$C$ is not invariant, there always exists a smallest invariant set $\orb
C$ containing $C$, called the \stress{orbit} of $C$ and given by
  $$
  \orb C = \bigcup_{t\in\F} \ha_t( C \cap \Xa_{t\inv} ).
  $$

Observe that the fundamental systems of open neighborhoods obtained in
\lcite{\NBDforUnbdd} and \lcite{\NBDforBounded}, if suitably intersected
with $\specUCK$, will collectively provide a basis for the topology of
$\specUCK$.  Therefore, taking the orbits of these sets, we will get a
collection of open invariant sets whose arbitrary unions consist of the
collection of all open invariant sets.

So, in order to characterize the open invariant sets, it is crucial to
understand the orbits of these \stress{basic} open sets.  Observe that
the sets $\W XYZ$ of \lcite{\NBDforBounded} are more general than the
sets $V_n$ of \lcite{\NBDforUnbdd}, since the latter set is among the
former, for suitable choices of $X$, $Y$, $Z$ and $\w$. In other
words, the collection of all $\W XYZ$, alone, gives a basis for the
topology of $\specUCK$.  Note, however, that we must allow the
``$\w$'' appearing in \lcite{\NBDforBounded} to vary.  To be precise,
let
  $$
  \W XYZ^\w =
  \left\{\matrix{
  \eta\in\specUCK : & \w\hfill\in\eta,\cr
  &\hfill \w x\inv\in\eta, & \hbox{ for $x$ in $X$,}\hfill\cr
  &\hfill \w y\inv\notin\eta, & \hbox{ for $y$ in $Y$,} \hfill\cr
  &\w z\hfill\notin\eta, & \hbox{ for $z$ in $Z$}\hfill\cr}
  \right\},
  $$
  where $X$, $Y$ and $Z$ are finite subsets of $\Gen$ and $\w$ is any
finite admissible word.  As already observed, the collection of all $\W
XYZ^\w$ will then form a basis of open sets for $\specUCK$.

Suppose $\w$ is an admissible word with $|\w|=n>0$ and let $\xi\in\W
XYZ^\w$.  Then $\w\in\xi$ and, since $\xi$ is convex, we must have that
$\trunc{\w}{n-1}=\w\last\w\inv \in \xi$.  So we would get a
contradiction if $\last\w\in Y$, which, in turn would imply that $\W XYZ
^\w$ is the empty set.  It is therefore advisable to require that
$\last\w\notin Y$, when
  $\w$ is of positive length.  In fact, it does not hurt if we, instead,
require that $\last\w\in X$.

  \state Proposition
  \label \SameOrbits
  Let $X$, $Y$ and $Z$ be finite subsets of $\Gen$ and let $\w$ be an
admissible word of finite positive length such that $\last\w\in X$.
Then the orbit of $\W XYZ^\w$ and that of $\W XYZ^e$ coincide.

  \proof
  Let $\xi\in\W XYZ^e$.  We claim that $\w\inv\in\xi$.  In fact, by
\lcite{\MembershipCriteria}, it suffices to show that
  $\last\w\in\root e\xi$, which is an immediate consequence of the fact
that $\last\w\in X$.  This implies that $\xi\in\Xa_{\w\inv}$ and hence,
since $\xi$ is arbitrary, that $\W XYZ^e \subseteq \Xa_{\w\inv}$.  One
now easily verifies that
  $\ha_\w(\W XYZ^e) = \W XYZ^\w$, from which the conclusion follows.
  \proofend

If we are to classify the orbits of all $\W XYZ^\w$, we may now restrict
attention to the case $\w=e$.

  \state Proposition
  \label \OrbWXYZ
  The orbit of each $\W XYZ^e$ consists of the set of all $\xi$ in
$\specUCK$ such that, for some $t\in\xi$, one has
  \iBitem $t x\inv\in\xi$, for all $x$ in $X$,
  \Bitem $t y\inv\notin\xi$, for all $y$ in $Y$, and
  \Bitem $t z \notin\xi$, for all $z$ in $Z$.

  \proof
  If $\xi$ satisfies the conditions stated for some $t$ in $\xi$ then
one clearly has that
  $$
  t\inv\xi\in\W XYZ^e\cap\Xa_{t\inv}
  $$ and that
  $\xi=\ha_t(t\inv\xi)$.  Conversely, it is evident that any point in
the orbit of
  $\W XYZ^e$ satisfies these conditions.
  \proofend

At this point it is convenient to expose a somewhat pictorial view of
our situation.  First of all, given a $\xi$ in $\specUCK$, one should
imagine that the points of $\F$ belonging to $\xi$ have been painted
black.  So $\xi$ is represented by a painted (necessarily convex) region
of the Cayley graph of $\F$.

With this in mind, the last result can be interpreted as saying that our
\stress{basic} open invariant sets, that is, the orbits of the
  $\W XYZ^\w$, consist of the configurations which show, at some vertex
and its adjacent ones, the pattern determined by $X$, $Y$ and $Z$.

  \state Proposition
  \label \ClassifyInvar
  Every open invariant subset of\/ $\specUCK$ is the union of a
collection of sets of the form $\orb{\W XYZ^e}$, as described in
\lcite{\OrbWXYZ}.

  \proof
  Obvious.
  \proofend

Of special relevance is the case where $Y$ and $Z$ are empty and $X =
\{x\}$ is a singleton.  In this case note that $\W XYZ^e =
\{\xi\in\specUCK: x\inv\in\xi\}$ which turns out to be nothing but
$\Xa_{x\inv}$.
  The orbit of $\Xa_{x\inv}$ will henceforth be denoted by $\E_x$ and,
since $\ha_x$ maps $\Xa_{x\inv}$ onto $\Xa_x$, we see that $\E_x$ is
also the orbit of $\Xa_x$.

According to \lcite{\OrbWXYZ}, note that in order for $\xi$ to belong to
$\E_x$ there must exist some $t$ in $\xi$ such that $tx\inv$ is also in
$\xi$.  With the change of variables $s=tx\inv$, this is the same as
saying that, for some $s$ in $\xi$, one has that $sx\in\xi$.  In other
words $\xi$ must contain (the vertices of) an edge of the Cayley graph
of $\F$ labeled $x$.  In symbols
  $$
  \E_x = \{\xi\in\specUCK : \exists s\in\xi,\ sx\in\xi\}.
  $$

  \state Proposition
  \label \ExInAnyInvarSet
  Every nonempty open invariant subset of\/ $\specUCK$ contains at least
one $\E_x$.

  \proof
  Let $U$ be a nonempty open invariant set and observe that, since the
unbounded elements are dense in $\specUCK$ by
\lcite{\DeterminationSpecOA}, $U$ must contain some unbounded element.
Applying \lcite{\NBDforUnbdd}, we conclude that there exists a finite
admissible word $\w$ such that the set
  $V = \{\eta\in\specUCK : \w\in\eta\}$ is contained in $U$.  We clearly
may suppose that $|\w|\geq1$ and so we let
  $x = \last\w$.
  It is easy to see that $V=\W XYZ^\w$, where $X=\{x\}$ and $Y=Z=\0$.
Now, since $U$ is invariant, it must contain the orbit of $\W XYZ^\w$
which, by
  \lcite{\SameOrbits}, coincides with the orbit of
  $\W XYZ^e$.  To conclude if suffices to note that
  $\W XYZ^e = \Xa_{x\inv}$ and that the orbit of the latter set is
precisely $\E_x$.
  \proofend

In order to discuss our next result we need to introduce a directed
graph associated to $A$, cf.~\cite{\KPRR,\KPR}.

  \definition
  \label \Grafone
  Given a 0--1 matrix $A = \{\A ij\}_{i,j\in\Gen}$ we will let
$\Graph$ be the directed graph whose vertex set is $\Gen$ and such that
the number of oriented edges from a vertex $i$ to the vertex $j$ is
precisely $\A ij$. By a \stress{path} in $\Graph$ we mean a finite or
infinite sequence $(x_1,x_2,\ldots)$ of elements of $\Gen$ such that
$\A{x_k}{x_{k+1}}=1$ for all $k$.

Notice that there is no conceptual difference between \stress{paths} and
\stress{admissible words} (see Definition \lcite{\AdmissibleWord}).  The
reason for our present interest in $\Graph$ stems from an increased
importance of the combinatorial properties of $A$, as we investigate the
dynamical properties of our partial action.

  \state Proposition
  \label \Hierarchy
  If $x,y\in\Gen$ and there exists a path from $x$ to $y$ in $\Graph$
then $\E_y\subseteq \E_x$.

  \proof
  Let $\xi\in\E_y$.  Then, by definition, $\xi$ contains an edge of the
form $(ty\inv,t)$ and hence, letting $\eta=t\inv\xi$, we have that
$\eta$ contains $y\inv$.  In the language of roots (see
\lcite{\RootDef}), this means that $y\in\root e\eta$.  Let $\b$ be a
finite admissible word such that $\b_1=x$ and $\last\b=y$.  We then
have, by \lcite{\MembershipCriteria}, that $\b\inv\in\eta$.  It follows
that $\eta$ contains the edge $(\b\inv,\b\inv x)$ and hence that
$\eta\in\E_x$.  Since $\E_x$ is invariant, it follows that $\xi=t\eta$
is in $\E_x$ as well.
  \proofend

If we combine the last two results, we see that an open invariant set
must contain not only one but several of the $\E_x$, depending, of
course, on the matrix $A$.  In particular, if $A$ is transitive, meaning
that there exists a path joining any two vertices in $\Gen$, one
concludes that any open invariant set must contain everything in
$\specUCK$ except, possibly, for the element $\phi = \{e\}$ which may or
may not be in $\specUCK$, as we saw in \lcite{\AboutUnits}.

There is another situation in which the $\E_x$'s are related:

  \state Proposition
  \label \Garfo
  Suppose that $x$ is in $\Gen$ and that the $x^{th}$ row of $A$ is
finite (meaning that $\{y\in\Gen:\A xy=1\}$ is a finite set).  Then
  $$
  \E_x=\bigcup_{{y\in\Gen} \atop {\A xy=1}}
  \E_y.
  $$

  \proof
  We already know, by \lcite{\Hierarchy}, that the $\E_y$ above are
contained in $\E_x$.  Conversely, let $\xi\in\E_x$, so that there exists
$t\in\xi$ such that $tx\inv\in\xi$ as well.  By
  \lcite{\MembershipCriteria} and \lcite{\MustContinue}
  we have that $ty\in\xi$ for some $y\in\Gen$, which implies that $\xi$
is in $\E_y$.  In view of \lcite{\CondTwo}, we have that $\A xy=1$.
  \proofend

We remark that the finiteness hypothesis in the above result cannot be
removed.  In fact, consider the following:

  \sysstate{Example}{\rm}
  {\label \Example
  Let $\Gen$ be the disjoint union of three copies of the set $\N$ of
natural numbers, indexed by
  $$
  \Gen = \{x_n: n\geq 1\} \cup \{y_n: n\geq 1\} \cup \{z_n: n\geq 1\}
  $$
  and suppose that the nonzero entries of $A$ are exactly the following
ones:
  \iBitem $\A{x_{n+1}}{x_n}$,
  \Bitem $\A{x_1}{y_n}$,
  \Bitem $\A{y_n}{z_1}$, and
  \Bitem $\A{z_n}{z_{n+1}}$,

  \medskip\noindent
  where, in all cases, $n$ takes all integer
  values greater than zero.}

In particular, $x_1$ is the only vertex of $\Graph$ which is the source
of infinitely many edges.  Observe that the column $\col{y_n}$ has
exactly one nonzero entry, namely $\A{x_1}{y_n}$.  It follows that
$\lim_n\col{y_n}$ is the column vector possessing a single nonzero entry
in the coordinate labeled by $x_1$.  This corresponds to the singleton
$\{x_1\}$ under the usual identification between 0--1 vectors and
subsets of $\Gen$.

Let $\xi=\{e,x_1\inv,x_1\inv x_2\inv, x_1\inv x_2\inv x_3\inv, \ldots\}$
and observe that $\xi$ is in $\specTCK$ as one may easily verify, with
the help of the description of $\specTCK$ given in \lcite{\TheoremOnLA}.
Now, observe that $\stem\xi=e$ and that $\root{e}{\xi} = \{x_1\}$.
Hence, by the discussion above and \lcite{\AcumIfANdOnlyIfIsLim} we have
that $\xi\in\specUCK$.

Note that $\xi$ is in $\E_{x_1}$, since it contains the edge
  $(x_1\inv,e)$, but $\xi$ is not in any of the $\E_{y_n}$ for lack of
an edge of the form $(t,ty_n)$.  Therefore the union of the $\E_{y_n}$
is strictly contained in $\E_{x_1}$, showing that our last result cannot
be extended without the finiteness hypothesis.

  \section{Fixed points}
  As promised earlier, we will now discuss fixed points.  In order to do
so we must begin by studying how roots (see Definition \lcite{\RootDef})
and stems (see Proposition \lcite{\TheWord}) change under our partial
action.

  \state Lemma
  \label \DynamicOnRootsAndWords
  Let $t\in\F$ and $\xi\in\Xa_{t\inv}$.
  Then
  $\stem{t\xi}=t\stem\xi$ and
  $\root{ts}{t\xi} =\root{s}{\xi}$,
  for every $s\in\xi$.

  \proof
  An element $x$ in $\Gen$ belongs to $\root{ts}{t\xi}$ if and only if
  $tsx\inv\in t\xi$ and this is the same as saying that
  $sx\inv\in\xi$, proving that $\root{ts}{t\xi} =\root{s}{\xi}$.

If $\xi$ is bounded then $\stem\xi$ is the maximum, among the elements
of $\xi$ for the order structure on $\F$ according to which
  $t\leq s$ if and only if $t\inv s\in\Pos$.  This can be proved easily
{}from
  \lcite{\MembershipCriteria}.  This order is clearly invariant under
left multiplication and hence
  $
  \max t\xi = t\max\xi,
  $
  establishing the validity of the first statement for bounded elements.

As for unbounded $\xi$'s we should first point out that the reference to
$t\stem\xi$ means the concatenation of the reduced form of $t$ with the
infinite word $\stem\xi$ after the possible (formal) cancelations are
performed.  It should be stressed that it is not obvious that the
resulting infinite string will contain only elements from $\Gen$ (as
opposed to elements from $\Gen\inv$).

Since $\xi\in\Xa_{t\inv}$ then $t\inv\in\xi$.  By
\lcite{\MembershipCriteria}, we have that $t\inv=\a\b\inv$, where
$|t|=|\a|+|\b|$, $\a$ is a finite sub-word of $\stem\xi$, and $\b$ is a
finite admissible word.  Therefore $t=\b\a\inv$ and then $t\stem\xi$ is
the result of replacing the initial segment $\a$ of $\stem\xi$ by $\b$.
In particular, this shows that $t\stem\xi$ can, in fact, be interpreted
as an infinite word.

Using \lcite{\MembershipCriteria} it is easy to see that $\stem\xi$ is
the only infinite word such that all of its finite sub-words belong to
$\xi$.  This implies that $t\stem\xi$ coincides with $\stem{t\xi}$
because every sub-word of $t\stem\xi$ belongs to $t\xi$.
  \proofend

We may use this to describe the fixed points for each $\ha_t$.

  \state Proposition
  \label \FixedPoints
  Suppose $t\neq e$ and let $\xi\in\Xa_{t\inv}$ be a fixed point for
$\ha_t$.  Then
  \iBitem there is a unique pair of finite words $\a$ and $\gamma$ such
that $t$ is equal to $\a\gamma^{\pm 1}\a\inv$ and $|t|=2|\a|+|\gamma|$,
  \Bitem $\stem\xi$ coincides with the eventually periodic word
$\a\gamma\gamma\gamma\ldots$, so that, in particular, $\xi$ is
unbounded, and
  \Bitem $\ha_t$ has no other fixed points.

  \proof
  Given that $\xi\in\Xa_{t\inv}$ we have that $t\inv\in\xi$ and hence
$t\inv=\a\b\inv$ where
  $t=|\a|+|\b|$, $\a$ is a finite sub-word of $\stem\xi$ and the
remaining properties of \lcite{\MembershipCriteria} hold.  So, write
$\stem\xi=\a\delta$ where $\delta$ is another word.
  From
  \lcite{\DynamicOnRootsAndWords} we have that
  $t\stem\xi = \stem{t\xi} = \stem\xi$ which means that
  $\b\a\inv\a\delta = \a\delta$, or simply $\a\delta=\b\delta$.

This says, in particular, that either $\a$ is a sub-word of $\b$ or
vice-versa.  Let us take up the case in which $\a$ is a sub-word of $\b$
and hence we can write $\b=\a\gamma$, from which
  $t=\b\a\inv= \a\gamma\a\inv$ and we will clearly have that
  $|t|=2|\a|+|\gamma|$.

Equation $\a\delta=\b\delta$ then becomes
  $\a\delta = \a\gamma\delta$ from which it follows that
$\gamma\delta=\delta$.
  Now, since $t\neq e$ we must also have $\gamma\neq e$ and so $\delta$
has to be infinite.  In fact we must have that
$\delta=\gamma\gamma\gamma\ldots$ and so
$\stem\xi=\a\delta=\a\gamma\gamma\gamma\ldots$.  The case in which $\b$
is a sub-word of $\a$ may be treated similarly.

With respect to uniqueness, observe that, if
  $t=\a\gamma^{\pm 1}\a\inv$, with
  $|t|=2|\a|+|\gamma|$, then $\a$ and $\gamma$ are uniquely determined
by $t$.  Therefore $\stem\xi$ is also uniquely determined by $t$.  The
uniqueness of the fixed point then follows from \lcite{\UniqueWordRoot}.
  \proofend

  \section{Topological freeness}
  A partial action of a group on a topological space is said to be
topologically free \scite{\ELQ}{Definition 2.1} (see also
\cite{\ArchSpiel}) if the set of fixed points $F_t$, for each partial
homeomorphism $\h_t$ with $t\neq e$, has empty interior.  We now want
to discuss the conditions under which the partial action of $\F$ on
$\specUCK$ possesses this property.

Recall from \lcite{\Grafone} that $\Graph$ is the directed graph whose
vertices correspond to the elements of $\Gen$ and such that there is an
oriented edge from vertex $i$ to vertex $j$ if and only if $\A ij=1$.
  The following Definition includes some well known concepts from graph
theory and is intended to pinpoint the aspects of $\Graph$ which we will
be looking at, cf.~\cite{\KPRR,\KPR}.

  \definition
  \label \Graphtwo
  \iBitem A \stress{circuit} in $\Graph$ is a finite path
$(x_1,\ldots,x_n)$ such that $\A{x_n}{x_1}=1$.
  \Bitem A circuit $(x_1,\ldots,x_n)$ is said to have an \stress{exit}
if, for some $k$, there exists $y\in\Gen$ with $\A{x_k}y=1$ and $y\neq
x_{k+1}$ (where $x_{k+1}$ is to be understood equal to $x_1$ if $k=n$).
  \Bitem A circuit is said to be \stress{terminal} when it has no exit.
  \Bitem A circuit $(x_1,\ldots,x_n)$ is said to be \stress{transitory}
if there is no path $(y_1,\ldots,y_m)$ with $m\geq3$ such that $y_1$ and
$y_m$ belong to $\{x_1,\ldots,x_n\}$ but $y_2$ does not.  In other
words, a transitory circuit is one to which no exiting path can return.

We would like to call the reader's attention to the fact that the term
\stress{loop} is sometimes used to refer to what we call
\stress{circuits}.  We have decided to adopt the terminology above
because it seems to be the one of choice among graph theorists and
also because in the established graph theory terminology the word
\stress{loop} refers to a single edge whose source and range coincide.

  \state Proposition
  \label \Topofree
  The partial action of $\F$ on $\specUCK$ is topologically free if and
only if $\Graph$ has no terminal circuits.

  \proof
  Let $e\neq t\in\F$ and let $F_t=\{\xi\in\Xa_{t\inv}:\ha_t(\xi)=\xi\}$
be the set of fixed points for $\ha_t$.  We already know that $F_t$ has
at most one element.  Suppose $F_t\neq\0$ and let $\xi$ be the unique
element in $F_t$.  Then the question about the interior of $F_t$ hinges
upon whether or not $\xi$ is an isolated point in $\specUCK$.  Write
$\stem\xi=\a\gamma\gamma\ldots$ as in \lcite{\FixedPoints}.  Since
$\stem\xi$ is admissible, $\gamma$ represents a circuit in $\Graph$.

Recall from \lcite{\NBDforUnbdd} that $\xi$ has a basis of neighborhoods
$\{V_n\}_{n\in\N}$ where each $V_n$ consists of all $\eta$'s in
$\specUCK$ containing $\trunc{\stem\xi}{n}$ (we are, of course,
referring to the topology of $\specUCK$ as a subspace of $\specTCK$
which, in turn, is the topological space under consideration in
\lcite{\NBDforUnbdd}).  So, to say that $\xi$ is an isolated point is to
say that there is an integer $n$ such that $V_n=\{\xi\}$.

If $\Graph$ has no terminal circuits then obviously $\gamma$ cannot be a
terminal circuit and it must therefore admit an exit.  Beginning with
this exit, say $\gamma_k$, one can construct an alternative infinite
path, hence an infinite word
  $\zeta$, such that
  $\trunc{\gamma}{k-1}{\zeta}$ is admissible.  If $\xi'$ is the unique
element of $\specUCK$ such that
  $\stem{\xi'}=\a\gamma\gamma\ldots\gamma\trunc{\gamma}{k-1}\zeta$,
where $\gamma$ is repeated sufficiently many times, we will have that
$\xi'$ contains
  $\trunc{\stem\xi}{n}$ and hence that $\xi'\in V_n$.  So $V_n$ cannot
be the singleton $\{\xi\}$.  That is, $\xi$ is not an isolated point and
hence $F_t$ has empty interior.  Since this holds for all $t\neq e$, we
see that our action is topologically free.

Conversely, if $\Graph$ admits a terminal circuit, say $\gamma$, then
let $\xi\in\specUCK$ be such that $\stem\xi=\gamma\gamma\ldots$.  Also
let $t=\gamma$.  We leave it for the reader to verify, using
\lcite{\MembershipCriteria}, that $t\inv\in\xi$, so that
$\xi\in\Xa_{t\inv}$.  In addition
  $\stem{t\xi} = t\stem\xi = \stem\xi$ and hence, by
\lcite{\UniqueWordRoot}, we have that $t\xi=\xi$.  That is, $\xi$ is a
fixed point for $\ha_t$.  Consider the neighborhood $V$ of $\xi$ given
by
  $
  V = \{\eta\in\specUCK : \gamma\in\eta\}.
  $
  We claim that $V=\{\xi\}$.  To prove this, suppose that $\eta$ is in
$V$, that is, $\gamma\in\eta$.  By \lcite{\MembershipCriteria} we have
that $\gamma$ is a sub-word of $\stem\eta$ and hence, since $\gamma$ is
a terminal circuit, $\stem\eta$ must be a sub-word of
$\gamma\gamma\ldots$.

Now, if $\stem\eta$ is infinite, we will have
  $\stem\eta = \gamma\gamma\ldots = \stem\xi$ and hence $\eta=\xi$, by
\lcite{\UniqueWordRoot}, as claimed.  To conclude the proof, we must
rule out the possibility that $\stem\eta$ be finite, but this is
precisely what \lcite{\MustContinue} does.
  \proofend

A topologically free partial action need not remain topologically free
when restricted to an invariant subset.  This is because, even if
$F_t$ has no interior, its intersection with an invariant subset may
have a nonempty interior in the relative topology.  For future
applications we would now like to study the conditions under which our
partial action is topologically free on closed invariant subsets.

  \state Proposition
  \label \TopoFreeOnsubsets
  The partial action of $\F$ on $\specUCK$ is topologically free on
every closed invariant subset of\/ $\specUCK$ if and only if\/ $\Graph$
has no transitory circuits.

  \proof
  Suppose that there exists a closed invariant subset
  $C\subseteq \specUCK$ where our partial action is not topologically
free.  Then there exists $t\neq e$ in $\F$ such that $F_t\cap C$ has
nonempty interior relative to $C$.  Let $\xi$ be a point in that
interior.  By \lcite{\FixedPoints}, fixed points are unique, so we
conclude that $\xi$ is an isolated point of $C$.  Obviously $\xi$ is
also an isolated point relative to its orbit $\orb \xi$.

We have thus shown that the existence of a closed invariant set where
our partial action is not topologically free implies the existence of a
fixed point that is isolated in its orbit.

We claim that the converse of this implication also holds.  To prove it
suppose that $\xi$ is fixed for some $\ha_t$, with $t\neq e$ and that
$\xi$ is an isolated point of $\orb\xi$.  Then there exists an open
subset $V$ of $\specUCK$ such that
  $$
  \orb\xi \cap V = \{\xi\}.
  $$
  It is easy to see that also $\overline{\orb\xi} \cap V = \{\xi\}$.  So
$\overline{\orb\xi}$ is a closed invariant subset of $\specUCK$ where
our action is not topologically free.  This proves our claim.

The conclusion then follows from the next Proposition.
  \proofend

  \state Proposition Let $t\neq e$ and let $\xi$ be a fixed point for
$\ha_t$. As in \lcite{\FixedPoints}, let $\a$ and $\g$ be finite words
such that
  $t = \a\gamma^{\pm 1}\a\inv$,
  $|t|=2|\a|+|\gamma|$, and
  $\stem\xi = \a\gamma\gamma\gamma\ldots$.
  Then $\xi$ is isolated in $\orb\xi$ if and only if $\g$ is a
transitory circuit in $\Graph$.

  \proof
  Suppose that $\g$ is a transitory circuit and let $V$ be the
neighborhood of $\xi$ given by
  $$
  V = \{\eta\in\specUCK : \a\g\in\eta\}.
  $$
  We claim that $V\cap\orb{\xi}=\{\xi\}$.  To see this, suppose that
$\eta$ is in $V\cap\orb{\xi}$, say $\eta=s\xi$ where $s\inv\in\xi$.
Then, by \lcite{\DynamicOnRootsAndWords}, we have that
  $
  \stem\eta =
  s\stem\xi =
  s\a\g\g\g\ldots,
  $
  where the convention adopted in the proof of
  \lcite{\DynamicOnRootsAndWords}, for multiplying group elements by
infinite words, is in effect.  In particular, $\stem\eta$ is eventually
periodic with period $\g$.

On the other hand, since $\eta\in V$, we have that $\a\g\in\eta$ and
hence, by \lcite{\MembershipCriteria}, $\a\g$ is a sub-word of
$\stem\eta$. So we see that $\stem\eta$ represents an infinite path in
$\Graph$ which, after following $\a$, loops around $\g$ once and
eventually returns to $\g$.  Now, since $\g$ is transitory, we conclude
that $\stem\eta$ must be equal to $\a\g\g\g\ldots$ So, in particular,
$\stem\eta=\stem\xi$ and hence $\eta=\xi$ by \lcite{\UniqueWordRoot}.

Conversely, supposed that $\xi$ is isolated in $\orb\xi$.  So there
exists a neighborhood $U$ of $\xi$ such that $U\cap\orb\xi=\{\xi\}$.
Then, by \lcite{\NBDforUnbdd}, $U$ must contain a set of the form
  $$
  V = \{\eta\in\specUCK : \a\g^n \in\eta\},
  $$
  where $n$ is an integer.  It follows that $V\cap\orb\xi=\{\xi\}$, as
well.

Let us argue by contradiction and hence suppose that $\g$ is not a
transitory circuit.  Therefore we may construct an infinite admissible
path of the form
  $\a\g^n\b\g\g\g\ldots$ which differs from $\a\g\g\g\ldots$.
  This is because from a circuit that is not transitory, one may leave
for a while and then return.

Let $\eta\in\specUCK$ be such that $\stem\eta= \a\g^n\b\g\g\g\ldots$ and
observe that $\eta$ is in $V$.  So, if we manage to prove that $\eta$ is
also in the orbit of $\xi$, we will arrive at a contradiction and hence
will complete the proof.  In order to do so all we must show is that
$\xi$ is in the domain of the partial homeomorphism associated to
  $t=\a\g^n\b\a\inv$, since we will then have that
  $$
  \stem{t\xi} = t\stem\xi = \a\g^n\b\a\inv \a\g\g\g\ldots =
  \a\g^n\b\g\g\g\ldots = \stem\eta,
  $$
  from which we will conclude that $\ha_t(\xi)=\eta$ and hence that
  $\eta$ is in $\orb\xi$, as desired.

To show that $\xi$ is in the domain $\Xa_{t\inv}$ of $\ha_t$ is
equivalent to showing that $t\inv=\a(\a\g^n\b)\inv\in\xi$, which we do
by using \lcite{\MembershipCriteria}.  In fact, $\a$ is a sub-word of
$\stem\xi$, $\a\g^n\b$ is an admissible word, whose last letter, namely
$\last\b$, belongs to $\root\a\xi$, precisely because $\b\g$ is
admissible.
  \proofend

  \section{Uniqueness of $\OA$}
  We now wish to address the question of \stress{uniqueness of} $\OA$.
This phrase has been used in the literature to refer to the fact that
any \cstar-algebra generated by a finite family $\{S_i\}_{i\in\Gen}$ of
\stress{nonzero} partial isometries satisfying the Cuntz--Krieger
conditions is isomorphic to $\OA$, provided $A$ satisfies a certain
condition (I)
  \scite{\CKbib}{Theorem 2.13}.

  \state Theorem
  \label \OAUnique
  Let $A = \{\A ij\}_{i,j\in\Gen}$ be a
  0--1 matrix with no identically zero rows and suppose that
$\Graph$ has no terminal circuits.  Then
  \iBitem any non-trivial ideal of $\uOA$ contains at least one of the
$S_i$, and
  \Bitem any non-trivial ideal of $\OA$ contains at least one of the
$S_i$.

  \proof
  Since $\OA$ is an ideal in $\uOA$ it suffices to prove the first
assertion.  The absence of terminal circuits guarantees, by
\lcite{\Topofree}, that the action of $\F$ on $\specUCK$ is
topologically free.  We may then use
\scite{\ELQ}{2.6} to conclude that any
non-trivial ideal in the \stress{reduced} crossed product
  $C(\specUCK)\crossproduct_r\F$
  has a non-trivial intersection with $C(\specUCK)$.

Since, by \lcite{\OurBundleAP}, we have that
  $C(\specUCK)\crossproduct \F
  \simeq
  C(\specUCK)\crossproduct_r \F$,
  we conclude that any nontrivial ideal $I$ in $\uOA$ satisfies
  $I\cap C(\specUCK) \neq \{0\}$.  This intersection, seen as an ideal
of $C(\specUCK)$, is then invariant under the partial action of $\F$ on
$C(\specUCK)$.  Let $U$ be the open subset of $\specUCK$ such that
  $I\cap C(\specUCK) = C_0(U)$, so that $U$ is invariant under the
partial action on $\specUCK$.  By \lcite{\ExInAnyInvarSet} we conclude
that $U$ contains some $\E_x$ and hence $\Xa_x\subseteq U$.

Recall from Section {\Preliminaries} that $1_x$ is the characteristic
function of $\Xa_x$.  So $1_x\in C_0(U) \subseteq I$ and since $1_x=S_x
S_x^*$, we conclude that $S_x$ is in $I$.
  \proofend

We therefore obtain the following generalization of the well known
  Cuntz--Krieger uniqueness Theorem \scite{\CKbib}{Theorem 2.13}.

  \state Corollary
  \label \UniquenessInOriginalForm
  Let $A = \{\A ij\}_{i,j\in\Gen}$ be a
  0--1 matrix with no identically zero rows and suppose that
$\Graph$ has no terminal circuits.  Let
  $\{S_i\}_{i\in\Gen}$
  and
  $\{T_i\}_{i\in\Gen}$
  be two sets of \stress{nonzero} \pisos on Hilbert space, both of
which satisfy \lcite{\JAUniversal.(i--iv)}.  Then the
  \cstar-algebras generated by $\{S_i\}_{i\in\Gen}$ and by
  $\{T_i\}_{i\in\Gen}$ are isomorphic to each other under an isomorphism
$\psi$ such that
  $\psi(S_i)=T_i$ for all $i\in\Gen$.

  \proof
  The representations $\pi$ and $\rho$ of $\OA$ given by
\lcite{\JAUniversal} relative to
  $\{S_i\}_{i\in\Gen}$
  and
  $\{T_i\}_{i\in\Gen}$, respectively, will both be faithful by
\lcite{\OAUnique}.  The isomorphism required is then obtained by taking
the composition $\psi = \rho\pi\inv$.
  \proofend

  \sysstate{Remark}{\rm}
  {Given the natural difficulties in generalizing the Cuntz--Krieger
conditions to the infinite case, one could ask whether or not our
definition of $\OA$ is the \stress{correct} generalization of the
definition given, in the finite case, by Cuntz and Krieger
\cite{\CKbib}.  Since all of the relations required by Definition
\lcite{\UnitalOADef} of the partial isometries $S_i$ hold in the
finite case, one is led to agree that these relations are justified.
So, the main question is whether or not we have imposed
\stress{sufficiently many} conditions, as one may argue that,
according to \lcite{\UnitalOADef}, condition \lcite{(\ELCond)} only
occurs in the perhaps few cases when finite sets $X$ and $Y$ can be
found such that $A(X,Y,j)=0$ for all but a finite number of $j$'s.  In
fact, it is easy to construct matrices $A$ for which \lcite{(\ELCond)}
does not interfere at all.  An example of this situation is given by
any matrix $A$ whose columns form a dense subset of $2^\Gen$, as the
reader may easily verify.  The answer to the above metamathematical
question is given by \lcite{\OAUnique} since, at least in the absence
of terminal circuits, we see that adding any genuinely new condition
will force at least one of the $S_x$ to be zero, which is obviously
undesirable.}

Corollary \lcite{\UniquenessInOriginalForm} can be used to prove the
result stated in \scite{\CKbib}{Remark 2.15} (cf.~also \cite{\KPR}).
Indeed, assuming that $\Graph$ has no terminal circuits (which is a
weaker hypothesis than condition (I) of \cite{\CKbib}, as already
observed in \scite{\KPR}{Lemma 3.3}), let
  $\{S_i\}_{i\in\Gen}$
  and
  $\{T_i\}_{i\in\Gen}$
  be sets of partial isometries with pairwise orthogonal ranges and
satisfying \CKcond2 (from Section \lcite{\Introduction}) with respect
to the strong operator topology.

Let $P = \sum_{j\in\Gen} S_j S_j^*$.  It is easy to show that, for all
$i$, one has that $S_i = S_i P = P S_i$.  We may therefore assume
without loss of generality that $P=1$ (cf.~\scite{\CKbib}{Section 2}).
That is, we may assume that our partial isometries satisfy \CKcond1,
as well.  Then, by \lcite{\SOT.(i)} and
\lcite{\UniquenessInOriginalForm}, the algebras generated by
  $\{S_i\}_{i\in\Gen}$
  and
  $\{T_i\}_{i\in\Gen}$
  are isomorphic.

However we point out that condition \CKcond2 in the strong topology is
unnecessarily restrictive since, for instance, it is not satisfied by
the Fock representation of ${\cal O}_\infty$ \cite{\Evans}, and, more
seriously, not being a C*-algebraic condition, it is not preserved by
quotients.

  \section{The simplicity criteria}
  We will now present a generalization of
  \scite{\CKbib}{Theorem 2.14},
  giving a sufficient condition for the simplicity of $\OA$.

  \state Theorem
  \label \Simplicity
  Let $A = \{\A ij\}_{i,j\in\Gen}$ be a
  0--1 matrix with no identically zero rows.  Suppose that $\Graph$
is transitive, in the sense that for every $x$ and $y$ in $\Gen$ there
is a path from $x$ to $y$.  If $\Gen$ is finite we assume, in addition,
that $A$ is not a permutation matrix.  Then $\OA$ is simple.

  \proof
  Suppose that $\Gen$ is infinite.  If $\gamma$ is a circuit in $\Graph$
then, given that $\Gen$ is infinite and a circuit is necessarily finite,
there must be an $x$ in $\Gen$ not in $\gamma$.  Using the transitivity,
one may find a path from anywhere in $\gamma$ to $x$.  Therefore there
are no terminal circuits in $\Graph$.  If, on the other hand, $\Gen$ is
finite and $A$ is not a permutation matrix, it is also easy to prove the
inexistence of terminal circuits.

Let $I$ be a non-trivial ideal in $\OA$.  Then $I$ is also an ideal in
$\uOA$ and hence, by \lcite{\OAUnique}, $I$ contains some $S_y$.  Now,
let $U$ be the open invariant subset of $\specUCK$ such that $I\cap
C(\specUCK) = C_0(U)$.  Since $S_y S_y^*=1_y\in I\cap C(\specUCK) =
C_0(U)$, we have that $\Xa_y\subseteq U$ and, since $U$ is invariant, we
see that $U$ also contains $\E_y$.  In addition, $U$ must contain all
$\E_z$ such that there exists a path from $y$ to $z$, by
\lcite{\Hierarchy}.  But then transitivity implies that $U$ contains the
union of all $\E_z$ which turns out to be precisely $\specUCK \backslash
\{\phi\}$.  It follows that $U$ contains
  $\Xa_z$ for all $z$ in $\Gen$ and hence that $C_0(U)$ contains the
corresponding $1_z$.  Finally, this implies that $I$ contains all of the
$S_z$ and hence also $\OA$.
  \proofend

In the row-finite case
  Theorem \lcite{6.8} of \cite{\KPRR}
  and
  Corollary \lcite{3.11} of \cite{\KPR}
  give conditions for the simplicity of $\OA$ based on the co-finality
of $\Graph$.  A graph is said to be \stress{co-final} if, given any
vertex $x$ and any infinite path $\b$, there exists a path starting in
$x$ and ending in some vertex of $\b$.
  However, we would like to argue that these results do not generalize
beyond the row-finite case. In fact, in the context of Example
\lcite{\Example}, the underlying graph is co-final and has no terminal
circuits and hence satisfies the hypotheses of both
\scite{\KPRR}{Theorem 6.8} and \scite{\KPR}{Corollary 3.11}.  However,
$\OA$ is not simple.  In fact, the union of all the $\E_{y_n}$ is a
proper open invariant set which therefore gives rise to a non-trivial
ideal of $\OA$ by \scite{\ELQ}{3.1}.

  \section{Classification of ideals}
  We are now going to use Theorem \lcite{3.5} from \cite{\ELQ} in
order to classify the ideals of $\uOA$ and therefore also of $\OA$.
According to the statement of that Theorem, if a partial action of a
group $G$ on a locally compact space $\Omega$ is topologically free on
every closed invariant subset and moreover satisfies the approximation
property then the correspondence
  $$
  U\mapsto \<C_0(U)\>,
  $$
  where $\<C_0(U)\>$ refers to the ideal generated by
  $C_0(U)$ in $C_0(\Omega)\crossproduct G$, is a bijection between the
collection of open invariant subsets of $\Omega$ and the closed
two-sided ideals of the crossed product.

  \state Theorem
  \label \ClassIdeals
  Let $A = \{\A ij\}_{i,j\in\Gen}$ be a 0--1 matrix with no
identically zero rows and suppose that $\Graph$ has no transitory
circuits.  Then the correspondence
  $
  U\mapsto \<C_0(U)\>
  $ is a bijection between the collection of open invariant subsets of\/
$\specUCK$ and the closed two-sided ideals of $\uOA$.

  \proof
  By \scite{\ELQ}{Theorem 3.5} it suffices to verify that
the partial action of $\F$ on $\specUCK$ satisfies the approximation
property and is topologically free on every closed invariant subset of
$\specUCK$.

Both these verifications are now immediate: the approximation property
follows from \lcite{\OurBundleAP} and topological freeness is a
consequence of \lcite{\TopoFreeOnsubsets}.
  \proofend

At this point we should recall that a classification of open invariant
sets was given in \lcite{\ClassifyInvar}.  We would also like to point
out that the situation for non-row-finite matrices is more complicated
than in the finite case studied in
  \scite{\Cuntz}{2.5} and the row-finite case studied in
  \scite{\KPRR}{6.6}.
  In fact, in the context of Example \lcite{\Example}, the results
mentioned do not apply since they would predict that $\OA$ is simple,
which is not the case.

  \section{Pure infiniteness}
  In this final section we exploit the crossed product structure of
$\OA$ to give a necessary condition for it to be purely infinite.  The
key idea is to use the circuits of the graph of $A$ to produce infinite
projections in hereditary subalgebras via an argument borrowed from
  \cite{\LS},
  cf.~\scite{\KPR}{Theorem 3.9}.
  As before, we let
  $\pr : \F \to \uOA$ be the universal partial representation of
  $\F$ satisfying $\RelCK$, and $\e(t)=\pr(t)\pr(t)^*$.

  \state Lemma
  \label\VariousProjections
  Let $n\geq1$,
  $x_1,\ldots,x_n\in\Gen$ and put
  $\g=x_1\ldots x_n\in\Pos$.
  \izitem If $1\leq m\leq n$ then $\e(\g) \leq \e(x_1\ldots x_m)$,
  \zitem If $\g$ is an admissible word then $\e(\g\inv) = \e(x_n\inv)$.
Othervise $\e(\g\inv) = 0$.
  \zitem If $\g$ is a circuit in $\Graph$ then
  $\e(\g) \leq \e(\g\inv)$,
  \zitem If, in addition to the assumption in (iii), $\g$ admits an exit
then
  $\e(\g) \neq \e(\g\inv)$.
  \zitem If $\a$ is an admissible word, $\g$ is a circuit with
$\last\a=\last\g$, and $v=\pr(\a\g)\pr(\a)^*$  then  $vv^*\leq v^*v$.
If, moreover, $\g$ has an exit then $vv^*\neq v^*v$.

  \proof
  Let $\g'=x_1\ldots x_m$ and $\g''=x_{m+1}\ldots x_n$.  Then $\g=\g'\g''$ and
  $$
  \e(\g) =
  \pr(\g')\e(\g'')\pr(\g')^* =
  \pr(\g')\pr(\g')^*\pr(\g')\e(\g'')\pr(\g')^* =
  \e(\g')\e(\g),
  $$
  proving (i).
  As for (ii), it follows from the equation
  $$
  \e(\g\inv) =
  \big(\prod_{k=1}^{n-1}\A{x_k}{x_{k+1}}\big)\;
  \e(x_n\inv),
  $$
  appearing as claim (1) in the proof of \lcite{\Suffering}.

Assuming $\g$ to be a circuit we have that $\A{x_n}{x_1} = 1$ and hence,
{}from
  \lcite{\OARelations.(iii)}, it follows that
  $\e(x_1) \leq \e(x_n\inv)$.  Therefore
  $$
  \e(\g) \leq
  \e(x_1) \leq
  \e(x_n\inv) =
  \e(\g\inv),
  $$
  proving (iii).

Suppose now that $\g$ has an exit, that is, there exists $j=1,\ldots,n$
and $y\in\Gen$ such that
  $\A{x_j}{y} = 1$ but $y\neq x_{j+1}$ (where $x_{j+1}$ is to be
understood equal to $x_1$ if $j=n$).  Define $\mu$ in $\F$ by
  $$
  \mu := \left\{
  \matrix{
  x_1 x_2 \cdots x_j {y}
    &\text{ if } 1 \leq j < n\cr
  {y}
    &\text{ if } j = n.\hfill
  }
  \right. 
  $$

  \noindent
  We will show that $\e(\g) \neq \e(\g\inv)$ by showing that
  \iBitem $\e(\mu) \perp \e(\g)$,
  \Bitem $\e(\mu) \leq \e(\g\inv)$, and
  \Bitem $\e(\mu)\neq 0$.

  \medskip\noindent
  With respect to the first point above let $\b=\mu y\inv$ so that $\b$
is a (possibly empty) subword of $\g$ and $\mu=\b y$.  We have
  $$
  \pr(\mu)^*\pr(\g) =
  \pr(y)^*\pr(\b)^*\pr(\b)\pr(\b\inv\g) =
  \pr(y)^*\e(\b\inv)\pr(\b\inv\g).
  $$
  If $j<n$ then $|\b|>0$ and, by (ii), $\e(\b\inv)=\e(x_j\inv)$.  Also,
since the positive word $\b\inv\g$ starts with $x_{j+1}$ we have that
  $\pr(\b\inv\g) = \e(x_{j+1}) \pr(\b\inv\g)$.  Using
  \lcite{\OARelations.(iii)} we conclude that the above equals
  $$
  \pr(y)^*\e(x_j\inv)\e(x_{j+1})\pr(\b\inv\g) =
  \A{x_j}{x_{j+1}}
  \pr(y)^*\e(x_{j+1})\pr(\b\inv\g) \$=
  \pr(y)^*\e(y)\e(x_{j+1})\pr(\b\inv\g) =
  0,
  $$
  where the last step follows from
  \lcite{\OARelations.(ii)} and the fact that $y\neq x_{j+1}$.
  Having proved that $\pr(\mu)^*\pr(\g) = 0$ we then have that
  $
  \e(\mu)\e(\g) =
  \pr(\mu)\pr(\mu)^*
  \pr(\g)\pr(\g)^* =
  0
  $,
  and hence that $\e(\mu) \perp \e(\g)$.

  If, on the other hand, $j=n$ then $\mu=y$ and hence
  $$
  \e(\mu)\e(\g) =
  \e(y)\e(\g) \leq
  \e(y)\e(x_1) =
  0,
  $$
  where, again, the last step follows from
  \lcite{\OARelations.(ii)} and $y\neq x_1$.

  Let us now prove that $\e(\mu) \leq \e(\g\inv)$.  If $j<n$ we have
  $$
  \e(\mu) \leq
  \e(x_1) \leq
  \e(x_n\inv) =
  \e(\g\inv),
  $$
  where the first step follows from (i), the second from
  \lcite{\OARelations.(iii)}, and the last one from (ii).

  If $j=n$ then
  $$
  \e(\mu) =
  \e(y) \leq
  \e(x_n\inv) =
  \e(\g\inv),
  $$
  where the crucial second step follows from
  the fact that $\A{x_n}y=1$, and
  \lcite{\OARelations.(iii)}.

  That  $\e(\mu)\neq 0$ follows from  \lcite{\ItExists.(ii)}.

  Let $\a$ be an admissible word and $\g$ be a circuit with
$\last\a=\last\g$.  Then
  $$
  v v^* =
  \pr(\a\g)\pr(\a)^*\pr(\a)\pr(\a\g)^* =
  \pr(\a) (\pr(\g)\pr(\a)^* \pr(\a) \pr(\g)^*)\pr(\a)^* \$\leq
  \pr(\a) (\pr(\g)\pr(\g)^*)\pr(\a)^* \leq
  \pr(\a) (\pr(\g)^*\pr(\g))\pr(\a)^* =
  \pr(\a)(\pr(\a\g)^*\pr(\a\g))\pr(\a)^* =
  v^* v
  $$
To see that $v v^*\neq v^* v$ when $\g$ has an exit, let $\mu$ be
defined as above.  One can then prove, as before, that
  \iBitem $\e(\a\mu) \perp \e(\a\g) \geq vv^*$,
  \Bitem $\e(\a\mu) \leq \pr(\a)\e((\a\g)\inv) \pr(\a)^* = v^* v$, and
  \Bitem $\e(\a\mu)\neq 0$.
  \proofend

  \state Theorem
  Let $A = \{\A ij\}_{i,j\in\Gen}$ be a
  0--1 matrix with no identically zero rows and suppose that
  \iBitem $\Graph$ has no terminal circuits, and
  \Bitem for every $x\in\Graph$, there exists a path starting at $x$ and
ending in a circuit.
  \medskip\noindent
  Then $\uOA$ and $\OA$ are purely infinite.

  \proof
  Since $\OA$ is an ideal in $\uOA$ it suffices to prove the statement
for $\uOA$.

  We adapt the argument used in
  \scite{\LS}{Theorem 9}
  to deal with local boundary actions to the present situation of $\uOA$
viewed as a crossed product; the proof here is technically simpler
because of the abundance of partial isometries.

Let $a$ be a nonzero positive element in $\uOA$. We will construct an
infinite projection in the hereditary subalgebra ${\cal H} (a)$
generated by $a$. Since rescaling $a$ does not affect ${\cal H} (a)$ we
may assume, without loss of generality, that $\|E(a)\| = 1$.

By \lcite{\Topofree} the partial action of $\F$ is topologically free.
Viewing $\uOA$ as $C(\specUCK) \crossproduct \F$, there exists, by
  \scite{\ELQ}{Proposition 2.4}, a function $h \in C(\specUCK)$ such
that
  \izitem $\| h E(a) h \| \geq \tfrac{7}{8}$,
  \zitem $\| h E(a) h - h a h \| \leq \tfrac{1}{8}$, and
  \zitem $ 0 \leq h \leq 1$.
  \medskip

We may assume further that the support of $h$ is contained in the set
  $
  \{\xi\in\specUCK: E(a)\calcat\xi > \frac{3}{4}\}.
  $
  Otherwise, multiply the original $h$ by a suitable positive function.

The subset $\{\xi\in \specUCK: h(\xi) > \tfrac{1}{2}\}$ is open in
$\specUCK$ so it contains $\Xa_\w$, for some $\w\in\Pos$, by the density
of the unbounded elements and
  \lcite{\NBDforUnbdd}.

Prolonging $\w$ until we reach a circuit $\g$, and relabelling $\g$ if
necessary, we get a path $\a$ such that $\last\a = \last\g$ and such
that $\w$ is an initial segment of $\a$.
  Therefore, using \lcite{\VariousProjections.(i)},
  $$
  \e(\a) \leq \e(\w) = 1_{\Xa_\w} ,
  $$
  from which it follows that the support of $\e(\a)$ is contained in the
set $\{\xi\in \specUCK: h(\xi) > \tfrac{1}{2}\}$ and hence that
  $\e(\a) \leq 2h$.

Let $v = \pr(\a\g) \pr(\a)^*$.  By
\lcite{\VariousProjections.(v)} we have that
  $v v^*$ is a proper subprojection of $v^*v.$
  Furthermore, $v^* v \in {\cal H} (h)$ because
  $v^* v \leq \pr(\a) \pr(\a)^* = \e(\a) \leq 2h$.
  By (iii),
  $v^* v \geq h v^*v \geq \tfrac{1}{2} v^*v$,
  so $hv^* v$ is invertible in ${\cal H} (v^*v)$, and its inverse, which
we denote by $g$, is bounded by 2.  Notice that
  $h v^* v \in {\cal H} (v^* v) \cap C(\specUCK)$,
  so
  $g \in {\cal H} (v^* v) \cap C(\specUCK)$.

Since $v v^* gh = v v^* v^* v = v v^*$, we have that
  $$
  \|v^*a v - v^* E(a) v\| = \|v^* v v^* (a - E(a)) v v^* v\|
   = \|v^* v v^* g h (a - E(a)) h g v v^* v\|,
  $$
  and using (ii) we obtain
  $$
  \|v^*a v - v^* E(a) v\| \leq \|g\|^2 \, \| h (a - E(a)) h \| \leq
  4 \frac{1}{8} = \frac{1}{2}.
  $$

Since the support of $vv^* \in C(\specUCK)$ is contained on the set on
which $h > \frac{1}{2}$, and since $\frac{3}{4} < E(a) \leq 1$ on the
support of $h$, we have that
  $$
  \frac{3}{4} v v^* \leq v v^* E(a) \leq v v^*
  $$
  and, multiplying on the left by $v^*$ and on the right by $v$,
  $$
  \frac{3}{4} v^* v \leq v^* E(a) v \leq v^*v.
  $$
  Thus
  $$
  \| v^*a v - v^* v\| \leq
  \|v^* a v - v^* E(a) v\| + \|v^* E(a) v - v^* v\| <
  \frac{1}{2} + \frac{1}{4},
  $$
  so $v^*a v$ is invertible in ${\cal H} (v^*v)$.

The rest of the proof is verbatim from
  \scite{\LS}{Theorem 9}
  and consists of initially showing that
  $w := a^{\frac{1}{2}} v (v^* a v)^{-\frac{1}{2}}$
  satisfies $w^* w = v^* v$, and then showing that $u:= w v w^*$ is a
partial isometry in ${\cal H} (a^{\frac{1}{2}}) = {\cal H} (a) $ such
that $uu^*$ is a proper subprojection of $ww^* = u^*u$.
  \proofend

  \vfill\eject 
  \centerline{\tensc Index of Symbols}

  \bigskip\noindent
  The symbol in the left hand side column is briefly described in the
middle column and its definition, or first occurrence, is to be found in
the section displayed in the right hand side column.

  \def\margin{\hbox to 10pt{\hfill}}
  \def\subtit#1{\bigskip \goodbreak \noindent \margin{\bf #1}\medskip}
  \def\symb#1;#2;#3;{\noindent
    \hbox to \hsize{\margin\hbox to 11ex{#1\hfill}
    #2\dotfill\hbox to 6ex{\hfill#3}
    \margin
    }}
  
  \subtit{Matrices}

  \symb
  $A$;
  Matrix of zeros and ones;
  {\Introduction};

  \symb
  $\Gen$;
  Index set for $A$;
  {\Introduction};

  \symb
  $\A i j$;
  $(i,j)$ entry of $A$;
  {\Introduction};

  \symb
  $\col j$;
  $j^{th}$ column of $A$;
  {\ColumnDef};

  \symb
  $A(X,Y,j)$;
  Same as $\prod_{x\in X} \A xj \prod_{y\in Y} (1-\A yj)$, for finite
subsets $X,Y\subseteq\Gen$;
  {\AXYj};

  \subtit{Free group}

  \symb
  $\F$;
  Free group on $\Gen$;
  {\Introduction};

  \symb
  $\Pos$;
  Positive cone of $\F$;
  {\PCKPisos};

  \subtit{Relations}

  \symb
  \PR i;
  Axioms for partial representations;
  {\Preliminaries};

  \symb
  \CKcond i;
  Classical Cuntz--Krieger relations;
  {\Introduction};

  \symb
  \TCKCond i;
  Toeplitz--Cuntz--Krieger relations;
  {\PCKPisos};

  \symb
  $\Rel$;
  Generic set of relations in terms of partial representations;
  {\Preliminaries};

  \symb
  $\Rel'$;
  Set of functions on $\OmegaE$ associated to $\Rel$;
  {\Preliminaries};

  \symb
  $\RelCK$;
  Cuntz--Krieger relations in terms of partial representations;
  {\OARelations};

  \symb
  $\RelTCK$;
  Toeplitz--Cuntz--Krieger relations in terms of partial
representations;
  {\TCKAlgSec};

  \subtit{Spectra}

  \symb
  $\OmegaE$;
  Subset of the power set of a group $G$ given by $\{\xi\in 2^G: e\in
\xi\}$;
  {\Preliminaries};

  \symb
  $\phi$;
  Smallest element of $\OmegaE$, namely $\{e\}$;
  {\SmallestGuy};

  \symb
  $\GenSpec\Rel$;
  Spectrum of a set $\Rel$ of relations (a subset of $\OmegaE$);
  {\Preliminaries};

  \symb
  $\specCK$;
  Same as $\specUCK\setminus \{\phi\}$;
  {\NonUnitOAasCrossProd};

  \symb
  $\specUCK$;
  Spectrum of Cuntz--Krieger relations = $\GenSpec{\RelCK}$;
  {\UCKAlgSec};

  \symb
  $\specTCK$;
  Spectrum of Toeplitz--Cuntz--Krieger relations = $\GenSpec{\RelTCK}$;
  {\TCKAlgSec};

  \symb
  $\W XYZ$;
  Fundamental neighborhoods for bounded elements of $\specTCK$;
  {\NBDforBounded};

  \symb
  $\E_x$;
  Certain invariant open subsets of $\specUCK$;
  {\InvarSection};

  \subtit{Algebras}

  \symb
  $\UAlg G\Rel$;
  Universal algebra for partial representations of $G$ subject to
$\Rel$;
  {\Preliminaries};

  \symb
  $\OA$;
  Cuntz--Krieger algebra;
  {\NonUnitalOADef};

  \symb
  $\uOA$;
  Unital Cuntz--Krieger algebra;
  {\UnitalOADef};

  \symb
  $\LA$;
  Toeplitz--Cuntz--Krieger algebra;
  {\TCKAlgSec};

  \subtit{Graphs and trees}

  \symb
  $\Graph$;
  Graph associated to a matrix $A$;
  {\Grafone};

  \symb
  $\STree ts$;
  Connected component, containing $t$, of the Cayley tree of $\F$ minus
the edge $(t,s)$;
  {\Subtree};

  \subtit{Words, stems, and roots}

  \symb
  $\w$;
  Generic finite or infinite word on $\Gen$;
  {\WordDef};

  \symb
  $|\w|$;
  Length of $\w$;
  {\WordDef};

  \symb
  $\w_i$;
  $i^{th}$ coordinate of $\w$;
  {\WordDef};

  \symb
  $\trunc \w n$;
  Initial segment of $\w$ of length $n$;
  {\WordDef};

  \symb
  $\set \w$;
  Set of group elements formed by the sub-words of $\w$;
  {\SetWDef};

  \symb
  $\stem \xi$;
  Stem of $\xi$, for $\xi\in\specTCK$ (a word);
  {\BoundedDef};

  \symb
  $\root t\xi$;
  Root of $t$ relative to $\xi$, for $\xi\in\specTCK$ and $t\in\xi$ (a
subset of $\Gen$);
  {\RootDef};

  \subtit{Dynamics}

  \symb
  $\h_t$;
  Generic partial homeomorphism associated to a group element $t$;
  {\Preliminaries};

  \symb
  $\X_t$;
  Range of $\h_t$;
  {\Preliminaries};

  \symb
  $\ha_t$;
  Canonical partial homeomorphism on $\specUCK$ associated to $t\in\F$;
  {\UnitalOAasCrossProd};

  \symb
  $\Xa_t$;
  Range of $\ha_t$;
  {\UnitalOAasCrossProd};

  \symb
  $\orb C$;
  Orbit of the set $C$ under a partial action;
  {\InvarSection};

  \References

  \def\www{Available from http://www.ime.usp.br/$\sim$exel/}

\Unpublished{\FAbadie,
  auth = {F. Abadie},
  title = {Tensor products of Fell bundles over discrete groups},
  institution = {Universidade de S\~ao Paulo},
  year = {1997},
  type = {preprint (funct-an/9712006)},
  note = {},
  NULL = {},
  }

\Article{\ArchSpiel,
  auth = {R. J. Archbold and J. Spielberg},
  title = {Topologically free actions and ideals in discrete dynamical
systems},
  journal = {Proc. Edinburgh Math. Soc.},
  year = {1993},
  volume = {37},
  pages = {119--124},
  NULL = {},
  }

\Article{\Blackadar,
  auth = {B. Blackadar},
  title = {Shape theory for $C^*$-algebras},
  journal = {Math. Scand.},
  year = {1985},
  volume = {56},
  pages = {249--275},
  NULL = {},
  author = {Bruce Blackadar},
  }

\Article{\ChoiEffros,
  auth = {M. D. Choi and E. G. Effros},
  title = {Nuclear $C^*$-algebras and the approximation property},
  journal = {Amer. J. Math.},
  year = {1978},
  volume = {100},
  pages = {61--79},
  NULL = {},
  }

\Article{\CuntzOInfinite,
  auth = {J. Cuntz},
  title = {Simple $C^*$-algebras generated by isometries},
  journal = {Comm. Math. Phys.},
  year = {1977},
  volume = {57},
  pages = {173--185},
  NULL = {},
  }

\Article{\Cuntz,
  auth = {J. Cuntz},
  title = {A class of $C^*$-algebras and topological Markov chains II:
reducible chains and the Ext-functor for $C^*$-algebras},
  journal = {Inventiones Math.},
  year = {1981},
  volume = {63},
  pages = {25--40},
  NULL = {},
  author = {Joachim Cuntz},
  }

\Article{\CKbib,
  auth = {J. Cuntz and W. Krieger},
  title = {A class of $C^*$-algebras and topological Markov chains},
  journal = {Inventiones Math.},
  year = {1980},
  volume = {56},
  pages = {251--268},
  NULL = {},
  author = {Joachim Cuntz and W. Krieger},
  }

\Article{\Evans,
  auth = {D. Evans},
  title = {On ${\cal O}_n$},
  journal = {Publ. Res. Inst. Math. Sci.},
  year = {1980},
  volume = {16},
  pages = {915--927},
  NULL = {},
  }

\Article{\newpim,
  auth = {R. Exel},
  title = {Circle actions on {$C^*$}-algebras, partial automorphisms and
a generalized {P}imsner--{V}oiculescu exact sequence},
  journal = {J. Funct. Analysis},
  year = {1994},
  volume = {122},
  pages = {361--401},
  NULL = {},
  akey = {Exel1994b},
  author = {Ruy Exel},
  number = {2},
  MR = {95g:46122},
  amsclass = {46L55 (46L80 47B35)},
  atrib = {IR},
  }

\Unpublished{\Inverse,
  auth = {R. Exel},
  title = {Partial actions of groups and actions of inverse semigroups},
  institution = {Universidade de S\~ao Paulo},
  year = {1995},
  type = {preprint},
  note = {to appear in {\sl Proc. Amer. Math. Soc\/}. \www},
  NULL = {},
  akey = {Exel1995b},
  author = {Ruy Exel},
  atrib = {S},
  }

\Article{\TPA,
  auth = {R. Exel},
  title = {Twisted partial actions, a classification of regular
{$C^*$}-algebraic bundles},
  journal = {Proc. London Math. Soc.},
  year = {1997},
  volume = {74},
  pages = {417-443},
  NULL = {},
  akey = {Exel1997b},
  author = {Ruy Exel},
  number = {3},
  atrib = {IR},
  }

\Unpublished{\Amena,
  auth = {R. Exel},
  title = {Amenability for {F}ell bundles},
  institution = {Universidade de S\~ao Paulo},
  year = {1996},
  type = {preprint},
  note = {to appear in {\sl J. reine angew. Math\/}. \www},
  NULL = {},
  akey = {Exel1996b},
  author = {Ruy Exel},
  atrib = {A},
  annote = {Journal f\"ur die Reine und Angewandte Mathematik},
  }

\Unpublished{\Ortho,
  auth = {R. Exel},
  title = {Partial representations and amenable Fell bundles over free
groups},
  institution = {Universidade de S\~ao Paulo},
  year = {1997. \www},
  type = {preprint},
  note = {},
  NULL = {},
  akey = {Exel1997d},
  author = {Ruy Exel},
  atrib = {N},
  annote = {},
  }

\Unpublished{\ELQ,
  auth = {R. Exel, M. Laca, and J. C. Quigg},
  title = {Partial dynamical systems and $C^*$-algebras generated by
partial isometries},
  institution = {University of Newcastle},
  year = {1997},
  type = {preprint},
  note = {},
  NULL = {},
  akey = {ExelLacaQuigg1997a},
  author = {Ruy Exel and Marcelo Laca and John C. Quigg},
  atrib = {N},
  annote = {},
  }

\Book{\FD,
  auth = {J. M. G. Fell and R. S. Doran},
  title = {Representations of *-algebras, locally compact groups, and
Banach *-algebraic bundles},
  publisher = {Academic Press},
  year = {1988},
  volume = {125 and 126},
  series = {Pure and Applied Mathematics},
  NULL = {},
  author = {J. M. G. Fell and R. S. Doran},
  }

\Book{\Glf,
  auth = {F. P. Greenleaf},
  title = {Invariant means on topological groups},
  publisher = {van Nostrand-Reinhold},
  year = {1969},
  volume = {16},
  series = {Mathematical Studies},
  NULL = {},
  author = {F. P. Greenleaf},
  }

\Article{\HR,
  auth = {A. an Huef and I. Raeburn},
  title = {The ideal structure of Cuntz--Krieger algebras},
  journal = {Ergod. Th. and Dynam. Sys.},
  year = {1997},
  volume = {17},
  pages = {611--624},
  NULL = {},
  }

\Article{\KPRR,
  auth = {A. Kumjian, D. Pask, I. Raeburn, and J. Renault},
  title = {Graphs, groupoids, and Cuntz-Krieger algebras},
  journal = {J. Funct. Anal.},
  year = {1997},
  volume = {144},
  pages = {505--541},
  NULL = {},
  author = {Alex Kumjian and David Pask and Iain Raeburn and Jean
Renault},
  }

\Unpublished{\KPR,
  auth = {A. Kumjian, D. Pask, and I. Raeburn},
  title = {Cuntz--Krieger algebras of directed graphs},
  institution = {University of Newcastle},
  year = {1996},
  type = {preprint},
  note = {},
  NULL = {},
  }

\Unpublished{\Laca,
  auth = {M. Laca},
  title = {Purely infinite simple Toeplitz algebras},
  institution = {University of Newcastle},
  year = {1997},
  type = {preprint},
  note = {},
  NULL = {},
  }

\Article{\LR,
  auth = {M. Laca and I. Raeburn},
  title = {Semigroup crossed products and Toeplitz algebras of
nonabelian groups},
  journal = {J. Funct. Analysis},
  year = {1996},
  volume = {139},
  pages = {415--440},
  NULL = {},
  }

\Article{\LS,
  auth = {M. Laca and J. Spielberg},
  title = {Purely infinite C*-algebras from boundary actions of discrete
groups},
  journal = {J. reine angew. Math.},
  year = {1996},
  volume = {480},
  pages = {125--139},
  NULL = {},
  }

\Article{\McCl,
  auth = {K. McClanahan},
  title = {$K$-theory for partial crossed products by discrete groups},
  journal = {J. Funct. Analysis},
  year = {1995},
  volume = {130},
  pages = {77--117},
  NULL = {},
  author = {Kevin McClanahan},
  }

\Article{\Nica,
  auth = {A. Nica},
  title = {$C^*$-algebras generated by isometries and Wiener-Hopf
operators},
  journal = {J. Operator Theory},
  year = {1991},
  volume = {27},
  pages = {1--37},
  NULL = {},
  author = {Andu Nica},
  }

\Article{\QR,
  auth = {J. C. Quigg and I. Raeburn},
  title = {Characterizations of crossed products by partial actions},
  journal = {J. Operator Theory},
  year = {1997},
  volume = {37},
  pages = {311--340},
  NULL = {},
  }

\Book{\Renault,
  auth = {J. Renault},
  title = {A groupoid approach to $C^*$-algebras},
  publisher = {Springer},
  year = {1980},
  volume = {793},
  series = {Lecture Notes in Mathematics},
  NULL = {},
  author = {Jean Renault},
  }

  \endgroup

  \bigskip
  \bigskip
  \noindent
  S\~ao Paulo, January 1998.

  \bye